\documentclass[a4pitaper,11pt]{amsart}
\usepackage{amsmath,amsthm,amssymb}
\usepackage[mathscr]{eucal}
 \usepackage{cite}
\usepackage{upgreek}
\usepackage[bookmarks,bookmarksnumbered]{hyperref}
\usepackage{enumerate}

\usepackage[dvipsnames]{xcolor}

\numberwithin{equation}{section}

\newcommand{\car}{\curvearrowright}

\newcommand{\G}{\Gamma}

\newcommand{\g}{\gamma}
\newcommand{\sg}{\sigma}

\newcommand{\ca}{\curvearrowright}

\newcommand{\emm}{\mathcal{M}}

\newcommand{\enn}{\mathcal N}

\newcommand{\la}{\lambda}

\newcommand{\El}{\mathcal{L}}

\newcommand{\rar}{\rightarrow}
\newcommand{\La}{\Lambda}

\newcommand{\cA}{\mathcal A}
\newcommand{\cB}{\mathcal B}
\newcommand{\cC}{\mathcal C}
\newcommand{\cD}{\mathcal D}
\newcommand{\cE}{\mathcal E}
\newcommand{\cF}{\mathcal F}
\newcommand{\cG}{\mathcal G}

\newcommand{\cK}{\mathcal K}
\newcommand{\cL}{\mathcal L}
\newcommand{\cM}{\mathcal M}
\newcommand{\cN}{\mathcal N}

\newcommand{\cP}{\mathcal P}
\newcommand{\cQ}{\mathcal Q}

\newcommand{\cU}{\mathcal U}

\newcommand{\pee}{\mathcal{P}}

\newcommand{\Sp}{\operatorname{Sp}}

\theoremstyle{plain}
\newtheorem{main}{Theorem}
\newtheorem{mcor}[main]{Corollary}

\newtheorem{theorem}{Theorem}[section]
\newtheorem{claim}{Claim}
\newtheorem{lemma}[theorem]{Lemma}
\newtheorem{proposition}[theorem]{Proposition}
\newtheorem{corollary}[theorem]{Corollary}
\theoremstyle{definition}

\newtheorem{definition}[theorem]{Definition}

\newtheorem{remark}[theorem]{Remark}

\begin{document}

\title[New examples of W$^*$ and C$^*$-superrigid groups]
{New examples of W$^*$ and C$^*$-superrigid groups}



\author{Ionu\c t Chifan}
\address{14 MacLean Hall, Department of Mathematics, The University of Iowa, IA, 52242, U.S.A.}
\email{ionut-chifan@uiowa.edu}
\thanks{I.C. has been supported in part by the NSF grants DMS-1600688, FRG DMS-1854194 and a CDA award from the University of Iowa}

\author{Alec Diaz-Arias}
\email{alec-diaz-arias@uiowa.edu}
\thanks{A.D-A. was supported in part by the AGEP- supplemental grant}

\author{Daniel Drimbe}
\address{Department of Mathematics, KU Leuven, Celestijnenlaan 200b, B-3001 Leuven, Belgium}
\email{daniel.drimbe@kuleuven.be}
\thanks {D.D. was partially supported by a PIMS fellowship and by a FWO fellowship.}


\begin{abstract} A countable discrete group $G$ is called $W^*$-superrigid (resp.\ $C^*$-superrigid) if it is completely recognizable from its von Neumann algebra $\cL(G)$ (resp. reduced $C^*$-algebra $C_r^*(G)$). Developing new technical aspects in Popa's deformation/rigidity theory we introduce several new classes of $W^*$-superrigid groups which appear as direct products, semidirect products with non-amenable core and iterations of  amalgamated free products and HNN-extensions. As a byproduct we obtain new rigidity results in $C^*$-algebra theory including additional examples of $C^*$-superrigid groups and explicit computations of symmetries of reduced group $C^*$-algebras.

\end{abstract}

\maketitle

\section{Introduction}

\subsection{Background}

\noindent The group von Neumann algebra $\cL(\G)$ of a countable discrete group $\G$ is defined as the weak operator
closure of complex group algebra $\mathbb C[\G]$ acting by left convolution on the Hilbert space $\ell^2 \G$ of square-summable functions on $\G$, \cite{MvN43}.
The classification of group von Neumann algebras has since been a main
theme in operator algebras centered around the following question: 
\emph{How does $\cL(\G)$ depend on $\G$?} The question
is particularly interesting when all non-trivial conjugacy classes of $\G$ are infinite  (abbrev.\ \emph{$\G$ is icc}), which corresponds to $\cL(\G)$ being a II$_1$ factor. While this is a fairly broad thematic, the main interest is to identify what purely algebraic aspects of $\G$ could be recovered from $\cL(\G)$ which in essence is a highly analytic environment. This is a rather complicated task as in general $\cL(\G)$ tends to have only a ``faded memory'' of $\G$. Perhaps the best illustration in this direction is A. Connes' celebrated result \cite{Co76} which asserts that for \emph{any} icc amenable groups $\G$ and $\La$ the corresponding von Neumann algebras are isomorphic, $\cL(\G)\cong \cL(\La)$. Thus besides amenability, a representation theoretic property of the group, $\cL(\G)$ has no recollection of any algebraic structure of the underlying group $\G$.

In the non-amenable case the situation is even more complex. For instance, methods in free probability show that for any collections of infinite amenable groups $\G_1,\G_2 ,...,\G_n$ and $\La_1,\La_2,...,\La_n$ where $n\geq 2$, the potentially non-isomorphic free products $\G=\G_1\ast \G_2\ast...\ast \G_n$ and $\La= \La_1\ast\La_2\ast ...\ast \La_n$ always give rise to isomorphic von Neumann algebras, $\cL(\G)\cong \cL(\La)$ \cite{Dy93}. Other examples of isomorphic von Neumann algebras arising from non-isomorphic non-amenable groups can be constructed using wreath products, see for instance \cite{Io06, Bo10, IPV10}.

\noindent Over the years however there have been discovered a variety of instances where $\cL(\G)$ is sensitive to various algebraic, analytic and representation theoretic properties of $\G$. A significant part of this progress was achieved through the emergence of Popa's deformation/rigidity theory in the early $2000$. This is a remarkably powerful conceptual framework that allows for various algebraic, dynamical, geometric and cohomological information of groups to be
completely recovered from their von Neumann algebras.

\noindent In this paper we are interested in an absolute form of reconstruction, namely when, up to isomorphism, $\G$ is entirely recoverable from $\cL(\G)$. Specifically, a group $\G$ is called \emph{W$^*$-superrigid} if the following holds: \emph{given any group $\La$ and any $\ast$-isomorphism  $\theta: \cL(\G)\rar \cL(\La)$,  then one can find a group isomorphism $\delta\in {\rm Isom}(\G,\La)$, a multiplicative character $\omega\in {\rm Char} (\G)$, and a unitary $w\in \cL(\La)$ such that $\theta = {\rm ad} (w) \circ \Psi_{\omega,\delta}$; here $\Psi_{\omega,\delta}$ denotes the canonical $\ast$-isomorphism given by $\Psi_{\omega,\delta}(u_g)=\omega(g) v_{\delta(g)}$, for all $g\in \G$, where $\{u_g \,:\,g\in \G\}$ and $\{v_h\,:\,h\in \La\}$ are the canonical group unitaries of $\cL(\G)$ and $\cL(\La)$, respectively.}

\noindent 
The first examples of $W^*$-superrigid groups were discoverd by Ioana, Popa, and Vaes in their groundbreaking work \cite{IPV10}. Moreover, their paper is particularly important as it introduces a conceptual approach towards the study of $W^*$-superrigidity through several novel techniques such as the analysis of comultiplication and height arguments. Developing several new technological aspects of these methods, a few other classes of examples of $W^*$-superrigid groups were found subsequently in \cite{BV12,Be14,CI17}.

 \subsection{Statements of the main results}

 \noindent In this paper we introduce several methods of constructing $W^*$-superrigid groups. Some on them are entirely new while others rely on prior constructions, mainly from \cite{IPV10,CI17}. We start by introducing our classes of groups, highlighting their features and importance along with the main results in the von Neumann algebric setting and various applications for the study of $C^*$-algebras.

 \noindent {\bf Class} $\mathcal{IPV}$. As already mentioned, the first examples of $W^*$-superrigid groups were discovered by Ioana, Popa and Vaes in \cite{IPV10}. Their groups arise via a specific generalized wreath product construction and some of its algebraic features play a key role in our work too. Therefore, to properly introduce our results, we briefly recall some of these examples below. Throughout the paper we denote by class $\mathcal {IPV}$ the collection of all generalized wreath product groups $\G=A \wr_{I} G$ satisfying the following conditions: 
\begin{enumerate}
    \item The base group $A \cong \mathbb{Z}_2$ \text{ or } $\mathbb{Z}_3$; 
    \item The acting group $G$ is any icc,  non-amenable  bi-exact group that contains an infinite property (T) normal subgroup;  
    \item The underlying set $I = G/K$ on which $G$ acts is the set of left cosets with respect to an infinite amenable malnormal\footnote{A subgroup $K<G$ is called malnormal if $gKg^{-1}\cap K=1$, for any $g\in G\setminus K$.} subgroup $K < G$. 
\end{enumerate}

\noindent There are many natural examples of groups in class $\mathcal {IPV}$, e.g.\ one can let $G$ be any uniform lattice in $Sp(n,1)$ $n\geq 2$ and $K$ be any maximal amenable subgroup of $G$ (see Section \ref{ipvclass} for other concrete examples of groups that belong to the class $\mathcal{IPV}$).

\noindent Besides being $W^*$-superrigid another important feature for us is that all groups in $\mathcal {IPV}$ are bi-exact 
\cite{BO08,De18}. Recall that a countable group $\G$ is called {\it bi-exact} (in the sense of \cite{Oz04}) if it is exact and admits a map $\mu:\G\to$Prob$(\G)$ such that lim$_{h\to\infty}\|\mu(ghk)-g\mu(h)\|\to 0$, for all $g,k\in\G.$ Other examples of bi-exact groups include all hyperbolic groups \cite{Oz03} and $\mathbb Z^2\rtimes$SL$_2(\mathbb Z)$ \cite{Oz08}.

\noindent The first main result of our paper is establishing product rigidity for von Neumann algebras of bi-exact groups, in the same vein with \cite[Theorem A]{CdSS15}. Namely, we have the following result:

\begin{main}\label{mainthm1}
	 Let $\G=\G_1\times\dots\times \G_n$ be a product of $n\ge 1$ icc, non-amenable, bi-exact groups and denote $\emm=\El(\G)$. Let $\La$ be any countable group and $t>0$ such that $\emm^t=\El(\La)$. 
	 
	 \noindent Then there exist a product decomposition $\La=\La_1\times\dots\times \La_n$, scalars $t_1,\dots,t_n>0$ with $t_1\cdots t_n=t$, and a unitary $u\in \emm^t$ such that $u\El( \La_i)u^*=\El(\G_i)^{t_i}$ for all $1\leq i\leq n$.	
	 \end{main}

\noindent In fact this generalizes the main results from \cite{CdSS15} in two ways. Firstly, it allows one to remove the weak amenability condition on the underlying groups in \cite[Theorem A]{CdSS15}. Secondly, it extends \cite[Theorem B]{CdSS15} from direct products of two to any finite number of groups. 

\noindent Consequently, when Theorem \ref{mainthm1} is combined with \cite[Theorem 8.4]{IPV10} we obtain the following examples of $W^*$-superrigid of direct product type.

\begin{mcor}\label{Cor:product}Let $\G_1,\G_2,...,\G_n \in \mathcal{IPV}$ and the product group $\G=\G_1 \times \G_2 \times ...\times \G_n$. Assume that $t>0$ is a scalar,  $\La$ is an arbitrary group and let $\theta: \El(\G) \to \El(\La)^t$ be an arbitrary $*$-isomorphism. 

\noindent Then $t=1$ and there exist $\delta \in {\rm Isom}(\G,\La)$, $ \omega \in {\rm Char}(\G)$ and  $u \in \mathscr{U}(\cL(\La))$ such that $\theta = {\rm ad}(u) \circ \Psi_{\omega,\delta}$. 

\end{mcor}


\noindent The class of all non-trivial finite product groups in class $\mathcal{IPV}$ can be used in conjunction with other canonical constructions in group theory to provide new examples of $W^*$-superrigid groups. Therefore, we introduce a family of groups that is constructed iteratively from the class $\mathcal {IPV}$ via a mix of two operations: amalgamated free product and HNN-extensions of groups belonging to a certain class of groups.

\noindent {\bf Class $\mathscr D$.} Let $\mathscr D_0$ be the class of all groups $G=\G_1\times \G_2\times \cdots \times \G_n$,  where $\G_i\in \mathcal{IPV}$ and $n\geq 2$. Also, for any $G\in \mathscr D_0$ we consider the set of factor subgroups $f_0(G)=\{G\}$ and the amalgamated subgroup $a_0(G)=1$. Then for every $i\geq 1$ we define inductively $\mathscr D_i$ 
 and for every $G\in \mathscr D_i$ its factor subgroups set $f_i(G)$ and its amalgamated subgroup $a_i(G)$ as follows. Assuming these are constructed, let $\mathscr D_{i+1}$ be the collection of all the groups that appear via one of the following constructions:
 \begin{itemize}
     \item [i)] amalgamated free product groups $G=G_1\ast_\Sigma G_2$ with $G_1,G_2\in \mathscr D_i$ and $\Sigma$ is an infinite, icc, amenable group such that ${\rm QN}^{(1)}_{G_1}(\Sigma)\footnote{If $\Sigma<G$ is a subgroup, then ${\rm QN}^{(1)}_G(\Sigma)$ denotes the one-sided quasi-normalizer of $\Sigma$ in $G$ and it is defined as the semigroup  of all $g\in G$ for which there exists a finite set $F\subset G$ such that $\Sigma g\subset F\Sigma$; see Section \ref{subsection.quasi.normalizer} for a more detailed account on the properties of quasi-normalizers.}={\rm QN}^{(1)}_{G_2}(\Sigma)=\Sigma$ and whenever $i\ge 1$ we require that $a_i(G_1)=a_i(G_2)=\Sigma$; 
     \item [ii)] HNN-extension groups $G={\rm HNN}(H,\Sigma, \phi)$ with $\phi: \Sigma \rightarrow H $ is a monomorphism such that $H\in \mathscr D_i$ and ${\rm QN}^{(1)}_H(\Sigma)=\Sigma$, ${\rm QN}^{(1)}_H(\varphi(\Sigma))=\varphi(\Sigma)$ , $[\Sigma: \Sigma\cap g\varphi(\Sigma)g^{-1}]=\infty$ and  $[\varphi(\Sigma) : \varphi(\Sigma) \cap g \Sigma g^{-1}]=\infty$, for all $g\in H$. Also, whenever $i\ge 1$ we require that $a_i(H)=\Sigma$. 
     \end{itemize}

In case i) we define $f_{i+1}(G)= f_i(G_1)\cup f_i(G_2)$ for the factor subgroups and $a_{i+1}(G)=\Sigma$ for the amalgamated subgroup set. In case ii) we define $f_{i+1}(G)= f_i(H)$ and $a_{i+1}(G)=\Sigma$, respectively. When is no confusion we will drop the $i$-subscript from the definitions of factors and amalgamated subgroups sets. 
For $i\ge 1$, we denote by $\mathscr D_i^m$ the subclass of groups $G\in \mathscr D_i$ for which there
exist $g_1,\dots, g_k\in G$ such that $\cap _{i=1}^k g_i\Sigma g_i^{-1}$ is finite, where $a(G)=\Sigma$.
Finally, denote by $\mathscr D=\cup^\infty_{i=0} \mathscr D_i$.

\begin{main}\label{mainthmclassd}Every group $G \in \mathscr D$ is $W^*$-superrigid. \end{main}

We refer the reader to the second part of Section 6 for concrete examples of groups that belong to $\mathscr D$. 
We also note that whenever $G\in \mathscr D_i$ for some $i\ge 1$ and $a_i(G)=\{\Sigma\}$, then  ${\rm QN}^{(1)}_{G}(\Sigma)=\Sigma$ (see Theorem \ref{theorem.classD.qntrivial}).

\noindent All known examples of semidirect product groups that are $W^*$-superrigid are generalized wreath product groups $A\wr_I G$ for which the base group $A$ is a
 finite abelian group \cite{IPV10,BV12,Be14}. Therefore, it is natural to investigate whether there are other kinds of $W^*$-superrigid semidirect groups beyond these family of examples which could potentially lead to new technological advancements. 
The following class of groups contains the first examples of $W^*$-superrigid groups that are semidirect product groups arising from actions on non-amenable groups.

\noindent {\bf Class $\mathscr {A}$}. Let $\G$ be a non-trivial, icc, bi-exact, torsion free, property (T) group. Let $n\geq 2$ be an integer and let  $\G_1$, $\G_2$, ..., $\G_n$ be isomorphic copies of $\G$. For every $1\leq i\leq n$ consider the action $\G \curvearrowright^{\rho^i} \G_i$ by conjugation, i.e.\ $\rho^i_\g (\la)=\g\la \g^{-1}$ for all $\g\in \G, \la\in \G_i$. Then let $\G\ca^{\rho} \G_1\ast \G_2\ast ...\ast \G_n$ be the action of $\G$ on the free product group $\G_1 \ast \G_2\ast...\ast \G_n$ induced by the canonical free product automorphism $\rho_\g = \rho^1_\g \ast \rho^2_\g \ast ...\ast \rho^n_\g$ for all $\g\in \G$ and denote by $G= (\G_1\ast \G_2 \ast... \ast \G_n) \rtimes_{\rho} \G$, the corresponding semidirect product. 

	\noindent Developing several new techniques in deformation/rigidity theory we were able to show that all groups in class $\mathscr A$ are completely recognizable from their von Neumann algebras. 
	
	\begin{main} \label{mainthm2} Every semidirect product group $G\in \mathscr A$ is $W^*$-superrigid.

	\end{main}

\noindent Many concrete examples of groups $G \in \mathscr A$ can be obtained appealing to methods in geometric groups theory, see Section \ref{semidirectwsuperrigid}. For example, one can start with any group $\G$ in any of the following categories: any uniform lattice in ${\rm Sp}(k,1)$ with $k\geq 2$; any torsion free property (T) group that is hyperbolic relative to any given finitely generated amenable subgroup via the method described in \cite[Theorem 1.1]{AMO}. In fact, it is worth noting that using methods in \cite{Os06,AMO} one can show that class $\mathscr A$ is uncountable (see item 3) in Proposition \ref{P:classA}) and that Theorem \ref{mainthm2} provides new  examples of residually finite $W^*$-superrigid groups (see item 2) in Proposition \ref{P:classA}) which add to the ones discovered previously in \cite[Theorem B]{BV12}. 

\noindent Another problem, closely related to $W^*$-superrigidity, is to investigate groups $G$ which are completely recognizable from their reduced $C^*$-algebra, $C^*_r(G)$; these are termed in the literature as \emph{$C^*$-superrigid groups}. Philosophically speaking, since $C^*_r(G)$ is a much ``smaller'' object than $\cL(G)$ it is reasonable to expect there should exist many $C^*$-superrigid groups. Despite this, unfortunately, very few examples are known in this direction.

\noindent The only examples of $C^*$-superrigid amenable groups are: all torsion free abelian groups by a classic result \cite{Sc74}, certain Bieberbach groups \cite{KRTW17}, some families of $2$-step nilpotent groups \cite{ER18}, and all free nilpotent groups of finite class and rank \cite{Om18}. The only known examples of non-amenable $C^*$-superrigid groups are the
amalgamated free product groups from \cite{CI17}.

\noindent There are a couple of points of contrast between these results. In the amenable case all groups are torsion free while the ones in the non-amenable setting could contain any type of torsion. In the amenable case all the results rely more or less on the $C^*$-superrigidity of abelian groups \cite{Sc74} as these are the building blocks for the groups considered. In the non-amenable case on the other hand the methods rely on deformation/rigidity arguments and their von Neumann algebraic superrigid behavior. Finally, in the amenable case the results always proceed by showing that an $\ast$-isomorphism between the $C^*$-algebras entail an abstract isomorphism between the underlying groups without explicitly connecting the two. By contrast, in the non-amenable case the methods used allow one to explicitly relate the two, essentially classifying all such $\ast$-isomorphisms.

\noindent In the same vein, our aforementioned rigidity results in the von Neumann algebraic setting shed new light towards the $C^*$-superrigidity problem  by providing many new examples of such non-amenable groups. Specifically, many of our groups have trivial amenable radical (see Propositions \ref{P:amenable radical} and \ref{P:classA}). Therefore, their $C^*$-algebras have unique trace by \cite[Theorem 1.3]{BKKO14} and hence, any $\ast$-isomorphism of these algebras will ``lift'' to the corresponding von Neumann algebras and using Theorems \ref{mainthmclassd} and \ref{mainthm2} we get the following:

\begin{mcor}  Let  $G\in \cup_{i\geq 1}\mathscr D_i^m \cup \mathscr A$. Let $\La$ be an arbitrary group and let $\theta: C^*_r(G) \to C^*_r(\La)$ be an arbitrary $*$-isomorphism. Then  there exist $\delta \in {\rm Isom}(G,\La)$, $ \omega \in {\rm Char}(G)$ and  $u \in \mathscr{U}(\cL(\La))$ such that $\theta = {\rm ad}(u) \circ \Psi_{\omega,\delta}$.  
\end{mcor}

	\noindent This result automatically enables us to describe all symmetries (automorphisms) of these algebras. 
	
\begin{mcor}\label{full} Let $G\in \cup_{i\geq 1}\mathscr D_i^m \cup \mathscr A$. Then for any $\theta\in {\rm Aut}( C^*_r(G))$ there exist $\delta \in {\rm Aut}(G)$, $ \omega \in {\rm Char}(G)$ and  $u \in \mathscr{U}(\cL(G))$ such that $\theta = {\rm ad}(u) \circ \Psi_{\omega,\delta}$.
\end{mcor}

\noindent A similar statement for the amalgamated free products considered in \cite{CI17} follows directly from \cite[Corollary C]{CI17}. To our knowledge, besides the Corollary \ref{full} above these are the only cases known of icc groups $G$ with $\cL(G)$  full factor for which the symmetries of $C^*_r(G)$ can be described entirely.     
	
\subsection{Organization of the paper.} \noindent Besides the introduction there are ten other sections  and an appendix in the paper. In Section 2 we recall some preliminaries and prove a few useful lemmas needed in the remainder of the paper. In Section 3 we use a new augmentation technique to prove
an intertwining result in von Neumann algebras that arise from products of bi-exact groups. We then continue in Section 4 with recalling some useful properties for groups that belong to the class $\mathcal {IPV}$. In Section 5 we use the result from Section 3 to prove our first main result, Theorem \ref{mainthm1}, and derive Corollary \ref{Cor:product} from it. We then continue in Sections 6 and 7 by presenting several properties for groups that belong to the classes $\mathscr D$ and $\mathscr A$, respectively. In Section 8 we provide a new situation where we can control the lower bound for height of certain unitary elements (Theorem \ref{thingroupheight}) and two technical results that provide ``discretization" results (Theorems \ref{discretizationgeneratinggroups} and \ref{discretizationhnngroups}). In Section 9 we present several results that allow us to reconstruct at the von Neumann algebra level the ``peripheral structure" of groups that belong to the classes $\mathscr D$ and $\mathscr A$. Finally, by using the machinery established in the previous sections, we present in Sections 10 and 11 the proofs of the remaining main results that are stated in the introduction.

\subsection{Acknowledgments} We are grateful to Stefaan Vaes for helpful comments and for kindly bringing to our attention that the W$^*$-superrigid groups from \cite{BV12} are residually finite. We also want to thank the anonymous referee for their numerous comments and suggestions which greatly improved the exposition and the overall quality of this paper. In particular, we thank the referee for pointing out a computational error in our original proof of Theorem  \ref{thingroupheight}.

	\section{Preliminaries}
	\subsection{Notations and Terminology}\label{terminology}

	\noindent Throughout this document all von Neumann algebras are denoted by calligraphic letters e.g.\ $\mathcal A$, $\mathcal B$, $\mathcal M$, $\mathcal N$, etc. Given a von Neumann algebra $\mathcal M$, we will denote by $\mathscr U(\mathcal M)$ its unitary group, by $\mathcal Z(M)$ its center, by $\mathscr P(\mathcal M)$ the set of all its nonzero projections and by $(\mathcal M)_1$ its unit ball.  Given a unital inclusion $\mathcal N\subseteq \mathcal M$ of von Neumann algebras we denote by $\mathcal N'\cap \mathcal M =\{ x\in \emm \,:\, [x, \enn]=0\}$ the relative commmutant of $\enn$ inside $\emm$ and by $\mathscr N_\emm(\enn)=\{ u\in \mathscr U(\emm)\,:\, u\enn u^*=\enn\}$ the normalizer of $\enn$ inside $\emm$. We say that $\cN$ is regular in $\cM$ if $\mathscr N_{\cM}(\cN)''=\cM$. We also denote by $W^*(S)$ the von Neumann algebra generated by a subset $S\subset \cM$.

	\noindent All von Neumann algebras $\emm$ considered in this document will be tracial, i.e.\ endowed with a unital, faithful, normal linear functional $\tau:M\rightarrow \mathbb C$  satisfying $\tau(xy)=\tau(yx)$ for all $x,y\in \emm$. This induces a norm on $\emm$ by the formula $\|x\|_2=\tau(x^*x)^{1/2}$ for all $x\in \emm$. The $\|\cdot\|_2$-completion of $\emm$ will be denoted by $L^2(\emm)$.  For any von Neumann subalgebra $\mathcal N\subseteq \mathcal M$ we denote by $E_{\mathcal N}:\mathcal M\rightarrow \mathcal N$ the $\tau$-preserving condition expectation onto $\mathcal N$. We denote the orthogonal projection from $L^2(\emm) \rightarrow L^2(\enn)$ by $e_{\enn}$. The Jones' basic construction \cite[Section 3]{Jo83} for $\enn \subseteq \emm$ will be denoted by $\langle \emm, e_{\enn} \rangle$.

	\noindent For any group $G$ we denote by $(u_g)_{g\in G} \subset \mathscr U(\ell^2G)$ its left regular representation, i.e.\ $u_g(\delta_h ) = \delta_{gh}$ where $\delta_h:G\rightarrow \mathbb C$ is the Dirac function at $\{h\}$. The weak operator closure of the linear span of $\{ u_g\,:\, g\in G \}$ in $\mathscr B(\ell^2 G)$ is called the group von Neumann algebra of $G$ and will be denoted by $\El(G)$; this is a II$_1$ factor precisely when $G$ has infinite non-trivial conjugacy classes (icc). Throughout this paper, for every subset $K\subseteq G$ we denote by $P_K$ the orthogonal projection from $\ell^2(G)$ onto the Hilbert subspace generated by the linear span of $\{\delta_g \,:\, g\in K \}$.

	\noindent All groups considered in this paper are countable and will be denoted by capital letters $A$, $B$, $G$, $H$, $Q$, $N$, $M$,  etc. Given groups $Q$, $N$ and an action $Q\curvearrowright^{\sg} N$ by automorphisms we denote by $N\rtimes_\sigma Q$ the corresponding semidirect product group. A group inclusion $H\leqslant G$ of finite index will be denoted by $H\leqslant_f G$. For any subgroup $H\leqslant G$ we denote by $C_G(H)=\{ g\in G\,|\, [g,H]=1\}$ its {\it centralizer} in $G$ and by $vC_G(H)=\{ g\in G\,|\, |g^H|<\infty\}$ its {\it virtual centralizer}. Note that $vC_G(G)=1$ precisely when $G$ is icc. Throughout the paper, we will also use the following observation: if $H<G$ is a subgroup satisfying $vC_G(H)=1$ (e.g. if $G$ is icc and $H<G$ has finite index), then $\cL(H)'\cap \cL(G)=1$.

\noindent Let $G$ be a group together with a family of subgroups $\mathcal F$. A set $K \subset G$ is called {\it small over $\mathcal F$} if there exist finite subsets $R,T\subset G$ and $\mathcal G\subseteq \mathcal F$ such that $K\subseteq \cup_{\Sigma \in \mathcal F} R \Sigma T$.  We denote by Sub$(G)$ the set of all the subgroups of $G$. If $G\car I$ is an action and $i\in I$, we denote by ${\rm Stab}_G(i)=\{g\in G|\; g\cdot i=i\}$ the {\it stabilizer of $i$ inside $G$}.

\noindent Finally, for any subset $S\subset \{1,\dots,n\}$ we denote its complement by $\hat S=\{1,\dots,n\}\setminus S$. If $S=\{i\},$ we will simply write $\hat i$ instead of $\widehat {\{i\}}$. Also, given
any product group $G=G_1\times \dots\times G_n$ we will denote the subproduct supported on $S$ by $G_S=\times_{i\in S}G_i$.

	\subsection{Popa's Intertwining Techniques}

\noindent Over fifteen years ago, S. Popa  introduced  in \cite[Theorem 2.1 and Corollary 2.3]{Po03} powerful analytic methods for identifying intertwiners between arbitrary subalgebras of tracial von Neumann algebras. These tools are now termed in the literature  as \emph{Popa's intertwining-by-bimodules techniques} and were highly instrumental to the classification of von Neumann algebras program via Popa's deformation/rigidity theory.

	\begin {theorem}[\cite{Po03}] \label{corner} Let $(\mathcal M,\tau)$ be a separable tracial von Neumann algebra and let $\mathcal P\subseteq p \emm p , \mathcal Q\subseteq  q\mathcal M q$ be von Neumann subalgebras. Let $\mathcal G\subset \mathscr U(\mathcal P)$ be a group such that $\mathcal G''= \mathcal P$. Then the following are equivalent:
	\begin{enumerate}
		\item There exist $ p_0\in  \mathscr P(\mathcal P), q_0\in  \mathscr P(\mathcal Q)$, a $\ast$-homomorphism $\theta:p_0 \mathcal P p_0\rightarrow q_0\mathcal Q q_0$  and a partial isometry $0\neq v\in q \mathcal M p$ such that $\theta(x)v=vx$, for all $x\in p_0 \mathcal P p_0$.
		\item There is no sequence $(u_n)_n\subset \mathcal G$ satisfying $\|E_{ \mathcal Q}(xu_ny)\|_2\rightarrow 0$, for all $x,y\in \mathcal  M$.
		\item There exist finitely many $x_i, y_i \in \mathcal M$ and $C>0$ such that  $\sum_i\|E_{ \mathcal Q}(x_i u y_i)\|^2_2\geq C$ for all $u\in \mathcal U(\mathcal P)$.
	\end{enumerate}
\end{theorem} 

\noindent If one of these equivalent conditions holds true,  then one writes $\pee\prec_{\emm} Q$, and says that {\it a corner of $\pee$ embeds into $\mathcal Q$ inside $\emm$.}
Furthermore, if $\mathcal P p'\prec_{\mathcal M} \mathcal Q$ for any non-zero projection $p'\in \pee '\cap p\emm p$ (equivalently, for any projection $0\neq p'\in\mathscr Z(\mathcal P'\cap  p  \mathcal M p )$), then we write $\pee \prec^{s}_{\emm }\mathcal Q$. We refer the readers to the survey papers \cite{Po07,Va10icm,Io12,Io18icm} for recent progress in von Neumann algebras using deformation/rigidity theory.

\noindent In the remaining part of the section we highlight a few technical intertwining results that will be used in an essential way to derive the main results of the paper. Some of them are either direct generalizations or follow from existent results in which case we only include some succinct proofs. For the new results we include more elaborated explanations.

\noindent The first lemma is a consequence of \cite[Lemma 2.4]{DHI16} and we omit its proof.

\begin{lemma}\label{L:bigger}
Let $(\cM,\tau)$ be a tracial von Neumann algebra and let $\cP\subset p\cM p$ and $\cQ\subset q\cM q$ be von Neumann subalgebras. Assume $\cP p'\prec_\cM^s \cQ$ for some non-zero projection $p'\in \cP'\cap p\cM p$.  Then there exists a non-zero projection $z\in \mathscr Z (\cP'\cap p\cM p)$ with $p'\leq z$ such that $\cP z\prec_\cM^s \cQ$.
\end{lemma}

\noindent In the proof of Theorem \ref{prodbiexactrig} we will need the following result that is essentially contained in \cite{HPV11}. Its proof is similar to the proof of \cite[Lemma 6.2]{Io11}, and we include it only for the reader's  convenience.

\begin{lemma}[\!\!\cite{HPV11}]\label{HPV}
Let $\Sigma<\Gamma$ be countable groups and denote $\cM=\cL(\Gamma)$. Let $\cB\subset \cM$ be a von Neumann subalgebra for which the quasi-normalizer of $\cB$ in $\cM$ is dense and $\cB\prec_\cM \cL(\Sigma)$. Let $\Omega$ be the subgroup of $\Gamma$ generated by all $\g\in\Gamma$ such that $\cB\prec_\cM \cL(\g\Sigma \g^{-1}\cap\Sigma)$. Then $\Omega$ has finite index in $\Gamma$.

\end{lemma}

We refer the reader to Section \ref{subsection.quasi.normalizer} for the definition of a quasi-normalizer of a subalgebra.

\noindent {\it Proof.}
 Let $\{u_\g\}_{\g\in\Gamma}$ be the canonical unitaries that generate $\cL(\Gamma).$ Following \cite[Section 4]{HPV11}, one can associate a projection $z(\Sigma_1)\in \cM$ to any subgroup $\Sigma_1<\Gamma$ such that $z(\Sigma_1)\neq 0$ if and only if $\cB\prec_\cM \cL(\Sigma_1)$. Moreover, $z(\g\Sigma_1 \g^{-1})=u_\g z(\Sigma_1)u_\g^*$, for any $\g\in\Gamma$ and $z(\Sigma_1\cap \Sigma_2)=z(\Sigma_1)z(\Sigma_2)$, for any subgroup $\Sigma_2<\Gamma$.

\noindent If $\Omega$ does not have finite index in $\Gamma$, then there exists a sequence of elements $(\g_n)_n\subset\Gamma$ such that $\cB\nprec_\cM \cL(\g_i^{-1}\g_j\Sigma \g_j^{-1}\g_i\cap \Sigma)$, for all $i\neq j.$ This is equivalent to $z(\g_j\Sigma \g_j\cap \g_i\Sigma \g_i)= 0$, for all $i\neq j.$ Hence, the projections $u_{\g_i}z(\Sigma)u_{\g_i}^*$, $i\ge 1$, are mutually orthogonal. Therefore, we deduce that $z(\Sigma)=0$, which implies that $B\nprec_M \cL(\Sigma)$, contradiction.
\hfill$\blacksquare$

\subsection{Quasinormalizers of groups and von Neumann algebras}\label{subsection.quasi.normalizer} 

\noindent Given a group inclusion $H<G$, the one-sided quasi-normalizer ${\rm QN}^{(1)}_G(H)$ is the semigroup of all $g\in G$ for which there exists a finite set $F\subset G$ such that $Hg\subset FH$ \cite[Section 5]{FGS10}; equivalently, $g\in {\rm QN}^{(1)}_G(H)$ if and only if $[H: gHg^{-1}\cap H]<\infty$. The quasi-normalizer ${\rm QN}_G(H)$ is the group of all $g\in G$ for which exists a finite set $F\subset G$ such that $Hg\subset FH$ and $gH\subset HF$.

\noindent Given an inclusion $\cN \subseteq \cM$ of finite von Neumann algebra we define the quasi-normalizer  $\mathscr {QN}_{\mathcal M}(\mathcal N)$ as the set of all elements $x\in \cM$ for which there exist $x_1,...,x_n\in \cM$ such that $\cN x\subseteq \sum x_i \cN$ and $x\cN \subseteq \sum \cN x_i$ (see \cite[Definition 4.8]{Po99}). Also the one-sided quasi-normalizer $\mathscr {QN}^{(1)}_{\mathcal M}(\mathcal N)$ is defined as the set of all elements $x\in \cM$ for which there exist $x_1,...,x_n\in \cM$ such that $\cN x\subseteq \sum x_i \cN$ \cite{FGS10}.

\noindent We record now some formulas for the quasi-normalizer of corners.

\begin{lemma}\emph{\cite{Po03,FGS10}}\label{QN1}
Let $\cP\subset \cM$ be tracial von Neumann algebras. For any projection $p\in\cP$, the following hold:
\begin{enumerate}
    \item $W^*({\rm \mathscr{QN}^{(1)}_{p\cM p}}(p\cP p))=pW^*({\rm \mathscr{QN}^{(1)}_{\cM}}(\cP))p$. 
    
    \item $W^*({\rm \mathscr{QN}_{p\cM p}}(p\cP p))=pW^*({\rm \mathscr{QN}_{\cM}}(\cP))p$. 
\end{enumerate}

\end{lemma}

We also mention the following remark which can be deduced directly from the definition.

\begin{remark}\label{QN.remark}
    Let $\cP\subset \cM$ be tracial von Neumann algebras. For any projection $p\in\cP'\cap\cM$, we have $W^*({\rm \mathscr{QN}_{p'\cM p'}}(\cP p'))=p'W^*({\rm \mathscr{QN}_{\cM}}(\cP))p'$.
\end{remark}

\noindent The following result provides a relation between the group theoretical quasi-normalizer and the von Neumann algebraic one.

\begin{lemma}\emph{\cite[Corollary 5.2]{FGS10}}\label{QN2}
Let $H<G$ be countable groups. Then the following hold:
\begin{enumerate}
    \item $W^*(\mathscr{QN}^{(1)}_{\cL(G)}(\cL(H)))=\cL(K)$,  where $K<G$ is the subgroup generated by ${\rm QN}^{(1)}_G(H)$. In particular, if ${\rm QN}^{(1)}_G(H)=H$, then $\mathscr{QN}^{(1)}_{\cL(G)}(\cL(H))=\cL(H)$.
    
    \item $W^*(\mathscr{QN}_{\cL(G)}(\cL(H)))=\cL({\rm QN}_G(H))$.
\end{enumerate}

\end{lemma}

\noindent We continue by emphasizing a few technical results regarding the control of quasinormalizers of von Neumann algebras subalgebras in various constructions including crossed products which are inspired by \cite[Theorem 3.1]{Po03}. We present a brief proof explaining how the same arguments from \cite{Po03} can be used.  
\begin{theorem}\label{controlquasi1} Let $\La$, $\Sigma$ be countable groups, let  $\La \curvearrowright^\rho \Sigma$ be an action by automorphisms and consider the corresponding semidirect product $\G=\Sigma \rtimes_\rho \La$. Denote by $\cM=\cL(\G)$ and $\cP =\cL(\La)$ assume that $\mathcal N \subseteq  \cP\subset \cM$ is a von Neumann subalgebra such that $\cN \nprec_{\cP} \cL({\rm Stab}_\La(\sigma))$ for all $\sigma\in \Sigma\setminus \{1\}$. Then we have that $\mathscr {QN}^{(1)}_\cM (\cN)''\subseteq \cP$. 
\end{theorem}

\noindent {\it Proof.} The conclusion follows immediately using the same arguments from \cite[Theorem 3.1]{Po03} once we show the following property:
given any sequence $(x_n)\subset \cN$ satisfying $\|E_{\cL({\rm Stab}_\La(\sigma))}(ax_n b)\|_2\rightarrow 0$ for all $a,b\in \cP$ and $\sigma\in \Sigma\setminus \{1\}$ we have that \begin{equation}\label{controlquasinorm1}
\|E_{\cP}(c x_n d)\|_2\rightarrow 0 \text{ for all } c,d\in \cM \ominus \cP.\end{equation}    
Using basic $\|\cdot \|_2$-approximations of $c$ and $d$ together with the $\cP$-bimodularity of $E_\cP$ one can easily see that it suffices to show \eqref{controlquasinorm1} only for $c=u_\mu$, $d=u_\sigma$ for $\mu,\sigma\in \Sigma\setminus \{1\}$.  Under these assumptions if we denote by $A_{\mu,\sigma} =\{ \la\in \Lambda \,:\, \rho_\la(\sigma)=\mu^{-1}\}$ basic computations show that $E_{\cP}(c x_n d)=E_{\cP}(u_\mu x_n u_\sigma)=\sum_\la \tau(x_n u_{\la^{-1}} ) \tau(u_{\mu\sigma_\la(\sigma)})u_\la= \sum_{\la\in A_{\mu,\sigma}} \tau(x_n u_{\la^{-1}} ) u_\la$. Since $A_{\mu,\sigma}= \nu {\rm Stab}_\La(\sigma)$ for some $\nu\in A_{\mu,\sigma}$ the above equation shows that $E_{\cP}(c x_n d)=u_\nu E_{\cL( {\rm Stab}_\La(\sigma))}(u_{\nu^{-1}} x_n)$ and using the hypothesis we get $\|E_{\cP}(c x_n d)\|_2=\| E_{\cL( {\rm Stab}_\La(\sigma))}(u_{\nu^{-1}} x_n)\|_2\rightarrow 0$ as $n\rightarrow \infty$. 
\hfill$\blacksquare$

For the following result, recall that a II$_1$ factor $\cM$ is called {\it solid} if for any diffuse subalgebra $\cA\subset \cM$, the relative commutant $\cA'\cap \cM$ is amenable. We refer the reader to Section \ref{section.amenable} for the notion of an amenable von Neumann algebra.

\begin{corollary} Let $\G$ be an icc, torsion free group such that $\cL(\G)$ is a solid von Neumann algebra.  Consider the product group $G= \G\times \G$ together with its diagonal subgroup $d(\G)=\{(\g,\g)\in G \,:\, \g\in\G  \}<G$. Let $p\in \cL(G)=\cM$ be a projection and assume that $\cA, \cB \subseteq p\cL(G)p$ are diffuse commuting von Neumann subalgebras such that $\cB$ has no amenable direct summand. Then $\cB \nprec_\cM \cL(d(\G))$.  

\end{corollary}

\noindent {\it Proof.} Assume by contradiction that $\cB \prec_\cM \cL(d(\G))$. Thus, one can find projections $b \in \cB$, $c\in \cL(d(\G))$, a non-zero partial isometry $v\in c\cM b$ and a $\ast$-isomorphism onto its image $\phi: b\cB b\rightarrow \cQ:=\phi(b \cB b))\subseteq c\cL(d(\G)) c$  such that $\phi(x)v=vx$ for all $x\in b \cB b$. Also note that $vv^*\in \cQ' \cap c \cM c$ and $v^*v\in b\cB b'\cap b\cM b$.

\noindent Next, we observe that the group $G= \G\times \G$ can be written alternatively as a semidirect product  $G= (\G\times 1) \rtimes_\rho d(\G)$ with respect to the action by conjugation of $d(\G)\curvearrowright^\rho \G\times 1$, i.e.\ $\rho_{(\g,\g)}(\la,1)= (\g \la \g^{-1},1)$,  for all $(\g,\g)\in d(\G)$ and $(\la,1)\in \G\times 1$. 
Then one can see that the stabilizers satisfy that ${\rm Stab}_{d(\G)}(\la,1)= d(C_\G(\la))$, where $C_\G(\la)$ is the centralizer of $\la$ in $\G$. Since $\G$ is torsion free and $\cL(\G)$ is solid it follows that the centralizer $C_\G(\la)$ and, hence, ${\rm Stab}_{d(\G)}(\la,1)$ is amenable for all $\la \neq 1$. Since $\cQ$ has no amenable direct summand, we have that $\cQ \nprec \cL({\rm Stab}_{d(\G)}(\la,1))$ and by Theorem \ref{controlquasi1} we get that $vv^*\in \mathscr{ QN}_{c\cM c}(\cQ)''\subseteq \cL(d(\G))$. Thus $vb\cB b v^*= \cQ vv^*\subseteq \cL(d(\G))$ and after extending $v
$ to a unitary $u$  we get $u \cB v^*v u^*\subseteq \cL(d(\G))$. Using Theorem \ref{controlquasi1} again we have that $u v^*v (\cB  \vee \cB '\cap p\cM p )v^*v u^*\subseteq \cL(d(\G))$. As $\G$ is icc after perturbing $u$ to a new unitary the previous relations imply that $u  (\cB  \vee \cB '\cap p\cM p )z u^*\subseteq \cL(d(\G))$, where $z$ is the central support of $vv^* \in \cB  \vee \cB '\cap p\cM p$. As $\cA \subseteq \cB'\cap p\cM p$ is diffuse this contradicts the solidity of $\cL(\G)$. \hfill$\blacksquare$

\noindent For further use, we record the following result which controls the intertwiners in algebras arising form certain subgroups. Its proof is essentially contained in \cite[Theorem 3.1]{Po03} (see also \cite[Lemma 2.7]{CI17}) so it will be left to the reader. 	  

\begin{lemma}\emph{\cite{Po03}}\label{malnormalcontrol} Let $H\leqslant G$ be countable groups and let $G \ca \mathcal N$ be a trace preserving action. Let $\mathcal P \subseteq p(\mathcal N \rtimes H)p$ be a von Neumann subalgebra such that $\mathcal P\nprec_{\mathcal N\rtimes H}  \cN\rtimes (g H g^{-1}\cap H)$ for all $g\in G\setminus H$. 

\noindent Then for all elements $x,x_1,x_2,...,x_l \in \mathcal N\rtimes G$ satisfying $\mathcal P x\subseteq \sum^l_{i=1} x_i \mathcal P$, we must have that $x p\in \mathcal N\rtimes H$. 
\end{lemma}

\noindent We also record the following result concerning von Neumann algebras of amalgamated free products and HNN-extension groups.

\begin{lemma}\emph{\cite[Theorem 1.1]{IPP05}}\label{L:ipp} Let $G=H*_\Sigma K$ be an amalgamated free product group or $G= {\rm HNN}(H,\Sigma,\varphi)$ is a HNN-extension group such that ${\rm QN}^{(1)}_G(\Sigma)=\Sigma$. Let $\cP\subset p\cL(H)p$ be a von Neumann subalgebra such that $P\nprec_{\cL(G)} \cL(\Sigma)$.

Then $\mathscr {QN}^{(1)}_{p\mathcal L(G)p} (\mathcal P)''\subseteq p\cL(H)p$.

\end{lemma}

\noindent {\it Proof.} Firstly, notice that $gHg^{-1}\cap H\subset \Sigma$, for any $g\in G\setminus H$. If $G$ is an amalgamated free product group, this is always true. On the other hand, if $G$ is an HNN-extension group as in the assumption, this follows from Lemma \ref{controlquasinorminhnn}. The lemma follows now from Lemma \ref{malnormalcontrol}. 
\hfill$\blacksquare$

\begin{lemma}\label{quasinormalizerdiagonal}
Let $\G$ be a countable non-amenable group such that for every $a\in \G\setminus \{1\}$ its centralizer $C_\G(a)$ is amenable. Then the diagonal subgroup $d(\G)<\G\times \G$  satisfies ${\rm QN}^{(1)}_{\G\times \G}(d(\G))=d(\G)$.
\end{lemma}

\noindent {\it Proof.} Let $(g,k)\in {\rm QN}^{(1)}_{\G\times \G}(d(\G))$. Thus, one can find $(g_i,k_i)\in \G\times \G$ with $1\leq i\leq n$ such that $d(\G)(g,k)\subseteq \bigcup^n_{i=1} (g_i,k_i)d(\G)$. Thus, for every $(\la,\la)\in d(\G)$ there exist an $i$ and 
$(\delta,\delta) \in d(\G)$ so that $(\la,\la)(g,k)=(g_i,k_i)(\delta,\delta)$. Basic calculations further imply that $g_i^{-1}\la g=\delta = k_i^{-1}\la k$; in particular, we have  $\la gk^{-1}= g_ik_i^{-1}\la$. Thus, if we denote $A_{i}=\{\la \in \La \, :\,g_i^{-1}\la g=\delta = k_i^{-1}\la k\}$, the above relations entail that $\G= \bigcup_{i=1}^n A_i$. However,  a simple calculation shows that $A_i$ is either empty or $A_i= \la_i C_{\G}(gk^{-1})$ for some $\la_i \in \G$. Combining with the previous relation we get that $\G= \bigcup_i \la_i C_{\G}(gk^{-1})$. In particular, we have  $[\G:C_\G(gk^{-1})]<\infty$ and as $\G$ is non-amenable we get that $C_\G(gk^{-1})$ is non-amenable as well. Then the hypothesis assumption implies that $gk^{-1}=1$ and, hence, $(g,k)\in d(\G)$, as desired. \hfill$\blacksquare$


\noindent We end this section by highlighting a result that allows us to obtain a genuine unitary conjugacy from some intertwining relations. The proof is essentially contained in the proof of \cite[Theorem A]{CI17} and we provide it for the reader's convenience.

\begin{theorem}\emph{\cite{CI17}} Let $A< G$ be icc groups such that ${\rm QN}^{(1)}_G(A)=A$ and denote by $\cM =\cL(G)$ the corresponding von Neumann algebra. Assume that $B<H$ are any groups satisfying $\cM =\cL(H)$, $\cL(A)\prec_\cM \cL(B)$ and $\cL(B)\prec^s_\cM \cL(A)$. Let $C<H$ be the subgroup generated by ${\rm QN}^{(1)}_H(B)$.

\noindent Then $[C:B]<\infty$ and there exists a unitary $w\in \mathscr U(\cM)$ such that $w\cL(A)w^*=\cL(C)$.  
\end{theorem}

\noindent
{\it Proof.} Since $\cL(A)\prec_{\cM}\cL(B)$, we can apply \cite[Lemma 2.4(4)]{DHI16} and obtain a non-zero projection $z\in \mathcal Z(\cL(B)'\cap \cM)\subset\cL(C)$ such that $\cL(A)\prec_{\cM} \cL(B)q'$, for any non-zero projection $q'\in \mathcal Z(\cL(B)'\cap \cM)z$. We continue by showing that 
\begin{equation}\label{c1}
\cL(B)z\nprec_{\cM}\cL(gAg^{-1}\cap A), \text{ for any }g\in G\setminus A.    
\end{equation}
Assume there exists $g\in G$ such that $\cL(B)z\prec_{\cM}\cL(gAg^{-1}\cap A)$. By \cite[Lemma 3.7]{Va08}, we have that $\cL(A)\prec_{\cM}\cL(gAg^{-1}\cap A)$. Since ${\rm QN}^{(1)}_G(A)=A$, it follows that ${\rm \mathscr {QN}^{(1)}_{\cM}}(\cL(A))=\cL(A)$, and hence, $\cL(A)\prec_{\cL(A)}\cL(gAg^{-1}\cap A)$. This implies by \cite[Lemma 2.5]{DHI16} that $[A: gAg^{-1}\cap A]<\infty$. Hence, $g\in {\rm QN}^{(1)}_G(A)=A$, which proves \eqref{c1}.

\noindent
{\bf Claim 1.} There exists a unitary $u\in\mathscr U(\cM)$ such that $uz\cL(C)zu^*\subset \cL(A)$. 

\noindent
{\it Proof of Claim 1.}
\noindent We first show that for any non-zero projection $q'\in (\cL(B)'\cap\cM)z$, there exists a non-zero projection $q''\in q'(\cL(B)'\cap\cM)q'$ such that $\cL(B)q''$ is unitarily conjugate into $\cL(A)$. Since $\cL(A)$ is a II$_1$ factor, it will follow that 
\begin{equation}\label{c2}
    u\cL(B)zu^*\subset \cL(A), \text{ for some unitary }u\in\mathscr U(\cM).
\end{equation}
Thus, take any non-zero projection $q'\in (\cL(B)'\cap\cM)z$. Since $\cL(B)q'\prec_{\cM}\cL(A)$, there exist projections $q\in\cL(B),r\in\cL(A)$, a non-zero partial isometry $w\in r\cM qq'$ and a $*$-homomorphism $\varphi: q\cL(B)qq'\to r\cL(A)r$ such that $\varphi(x)w=wx$, for any $x\in q\cL(B)qq'$. We can moreover assume that the support projection of $E_{\cL(\Sigma)}(ww^*)$ equals $r$. Let $P=\varphi(q\cL(B)qq')\subset r\cL(A)r$ and write $w^*w=qq_0$ for a projection $q_0\in q'(\cL(B)'\cap\cM)q'$. One can check that \eqref{c1} implies that $P\nprec_{\cL(A)} \cL(gAg^{-1}\cap A)$, for any $g\in G\setminus A$. By applying Lemma \ref{malnormalcontrol}, we derive that $ww^*\in \cL(A)$, and thus, $w(q\cL(B)qq_0)w^*\subset \cL(A).$ Let $z_0$ be the central support of $q$ in $\cL(B)$. Since $\cL(A)$ is a II$_1$ factor, it follows that there exists $\eta\in \mathscr U(\cM)$ such that $\eta \cL(B)z_0q_0\eta^*\subset \cL(A)$. We now take $q''=z_0q_0$ and therefore obtain that relation \eqref{c2} holds.

\noindent Thus, we take a unitary $u\in\mathscr U(\cM)$ such that $u\cL(B)zu^*\subset \cL(A)$ and let $e=uzu^*\in \cL(A)$. By \eqref{c1}, we have that $\mathscr{QN}_{e\cM e}^{(1)}(u\cL(B)zu^*)\subset e \cL(A)e$. By using the quasi-normalizer formulas Lemma \ref{QN1} and  Lemma \ref{QN2}, we deduce that $uz\cL(C)zu^*\subset e\cL(A)e$.
\hfill$\square$

\noindent
{\bf Claim 2.} There exists a non-zero projection $z'\in z\cL(C)z$ such that $uz'\cL(C)z'u^*=p\cL(A)p.$

\noindent
{\it Proof of Claim 2.} Denote $Q=u\cL(B)zu^*\subset e\cL(A)e$ and notice that $e\cL(A)e\prec_{\cM} Q$ since $\cL(A)$ is a II$_1$ factor. Thus, there exist projections $p\in e\cL(A)e, q\in Q,$ a non-zero partial isometry $v\in q\cM p$ and a $*$-homomorphism $\theta: p\cL(A)p\to qQq$ such that $\theta(x)v=vx$, for any $x\in p\cL(A)p$. Since $Q\subset e\cL(A)e$, we derive that $v\in \mathscr{QN}^{(1)}_{\cM}(\cL(A)=\cL(A)$. Moreover, we may assume that $v^*v=p$ and $p\in Q$ since $\cL(A)$ is a II$_1$ factor.

\noindent Next, if $x\in p\cL(A)p$, then $vx(pQp)\subset (qQq)v$. Hence, $vx\in W^*(\mathscr{QN}^{(1)}_{e\cL(A)e}(Q))$ (see the proof of \cite[Lemma 3.5]{Po03}). This shows that $p\cL(A)p\subset W^*(\mathscr{QN}^{(1)}_{e\cM e}(Q))$. Since $\mathscr{QN}^{(1)}_{e\cM e}(Q)\subset uz \cL (C)zu^*$, we derive that $p\cL(A)p\subset p(uz \cL (C)zu^*)p$. By letting $z'=u^*pu\in z\cL(C)z$, we have $p\cL(A)p\subset uz'\cL(C)z'u^*$. The claim is proven since the reversed inclusion follows from Claim 1. \hfill$\square$

\noindent We continue by proving that $[C:B]<\infty$. One can check that $(u^*vu)z\cL(C)z\subset \cL(B)zu^*vu$ and $u^*vu\cL(B)\subset \cL(B)u^*vu$. This shows that $u^*vu\in z\cL(C)z$ and hence, $\cL(C)\prec_{\cL(C)}\cL(B)$. By applying \cite[Lemma 2.5]{DHI16}, we get that $[C:B]<\infty$. This implies that $\mathscr{QN}^{(1)}_{\cM}(\cL(C))=\cL(C)$ and $\cL(C)\prec^s_{\cM} \cL(A)$.
Since $\cL(A)$ is a II$_1$ factor, we can use Claim 2 combined with \cite[Lemma 2.6]{CI17} and derive that there exists a unitary $w\in \cM$ such that $w\cL(A)w^*=\cL(C)$.   
\hfill$\blacksquare$

\subsection{Relative amenability}\label{section.amenable} 
\noindent A tracial von Neumann algebra $(\cM,\tau)$ is {\it amenable} if there exists a positive linear functional $\Phi:\mathbb B(L^2(\cM))\to\mathbb C$ such that $\Phi_{|\cM}=\tau$ and $\Phi$ is $\cM$-{\it central}, meaning $\Phi(xT)=\Phi(Tx),$ for all $x\in \cM$ and $T\in \mathbb B(L^2(\cM))$. The celebrated theorem of Connes asserts that a von Neumann algebra $M$ is amenable if and only if it is approximately finite dimensional \cite{Co76}. 

\noindent We recall now the relative version of this notion due to Ozawa and Popa \cite{OP07}. Let $(\cM,\tau)$ be a tracial von Neumann algebra. Let $p\in \cM$ be a projection and $P\subset p\cM p,\cQ\subset \cM$ be von Neumann subalgebras. Following \cite[Definition 2.2]{OP07}, we say that $\cP$ is {\it amenable relative to $Q$ inside $M$} if there exists a positive linear functional $\Phi:p\langle \cM,e_\cQ\rangle p\to\mathbb C$ such that $\Phi_{|p\cM p}=\tau$ and $\Phi$ is $\cP$-central. 
Note that $\cP$ is amenable relative to $\mathbb C$ inside $\cM$ if and only if $\cP$ is amenable. 
We also say that $\cP$ is {\it strongly non-amenable relative to $Q$ inside $M$} if $Pp'$ is non-amenable relative to $Q$ for any non-zero projection $p'\in P'\cap pMp.$

\noindent The following lemma is well known and we leave the proof to the reader.

\begin{lemma}\label{L:nonamenable}
Let $\Sigma<\Gamma$ be countable non-amenable groups. Then $\cL(\Sigma)q$ is non-amenable for any non-zero projection $q\in\cL(\Sigma)'\cap \cL(\Gamma)$.

\end{lemma}

\section{Bi-exact groups and an augmentation technique}

\noindent One of the technical ingredients needed in the proof of Theorem \ref{mainthm1} is Proposition \ref{P:bi-exact} which provides some intertwining results in von Neumann algebras of products of bi-exact groups. 

\begin{proposition}\label{P:bi-exact}
Let $\Gamma=\Gamma_1\times\dots\times\Gamma_n$ be a product of $n\ge 1$ non-amenable bi-exact icc groups and denote $\cM=\cL(\Gamma)$. Assume that $\cM=\bar\otimes_{0\leq j\leq k} \cP_j$ is a decomposition into $k+1$ II$_1$ factors such that $\cP_j$ is non-amenable for any $2\leq j\leq k.$ 
Let $\Sigma<\Lambda$ be some countable groups such that $\cM=\cL(\Lambda)$ and $\bar\otimes_{1\leq j\leq k} \cP_j\prec_\cM \cL(\Sigma)$. 

\noindent Then $\cL(\Sigma)\nprec_\cM \cL(\Gamma_S)$, for any subset $S\subset \{1,\dots,n\}$ that has at most $k-1$ elements. 

\end{proposition}

Although is not needed for our proof, we mention that each tensor factor $\cP_j$ has to be of a specific form. Indeed, the unique prime factorization result of Ozawa and Popa \cite{OP03,BO08} implies that there exists a partition $T_0\sqcup T_1\sqcup\dots\sqcup T_k=\{1,\dots,n\}$ such that for any $j\in\{0,1,\dots,k\}$, $\cP_j$ equals $\cL(\Gamma_{T_j})$, up to unitary conjugacy and amplification.

\noindent In the particular case when $\cL(\Sigma)\subset \cM$ is regular, the conclusion of Proposition \ref{P:bi-exact} follows immediately by applying Theorem \ref{Th:bo}. Indeed, assume that $\cL(\Sigma)\prec_\cM \cL(\Gamma_S)$ for a subset $S\subset \{1,\dots,n\}$ that has at most $k-1$ elements. Then by \cite[Lemma 2.4(2)]{DHI16} and \cite[Lemma 3.7]{Va08} we get that $\bar\otimes_{1\leq j\leq k} \cP_j\prec_\cM \cL(\Gamma_S)$, which will imply by repeated application of Theorem \ref{Th:bo} the contradiction that $P_1$ is not diffuse. The general case is more subtle, the idea is to exploit the group von Neumann algebra structure of $\cL(\Sigma)\subset \cL(\Lambda)$ and to make the analysis by considering a Bernoulli action of $\Lambda$.

\noindent We record now the following relative solidity result for von Neumann algebras arising from products of bi-exact groups. The result is a direct consequence of \cite[Theorem 15.1.5 and Lemma 15.3.3]{BO08}

\begin{theorem}\emph{\cite{BO08}}\label{Th:bo}
Let $\Gamma=\Gamma_1\times\dots\times\Gamma_n$ be a product of $n\ge 1$ non-amenable bi-exact groups and denote $\cM=\cL(\Gamma)$. Let $\cQ\subset q\cM q$ be a von Neumann subalgebra such that $\cQ'\cap q\cM q$ is non-amenable.

\noindent Then there exists $1\leq i\leq n$ such that $\cQ\prec_\cM \cL(\Gamma_{\hat i})$.
\end{theorem}

\noindent Before proceeding to the proof of Proposition \ref{P:bi-exact}, we need the following two useful lemmas.

\begin{lemma}\label{L:dual}
Let $\Lambda\curvearrowright \cB$ be a trace preserving action and denote $\cM=\cB\rtimes\Lambda$. Let $\cP\subset \cL(\Lambda)$ and $\cQ\subset \cM$ be some von Neumann subalgebras.  
Following \cite{PV09}, let $\Delta: \cM\to \cM\bar\otimes \cL(\Lambda)$ be the $*$-homomorphism given by $\Delta(b)=b\otimes 1$, for any $b\in \cB$ and $\Delta(v_g)=v_g\otimes v_g$, for any $g\in\Lambda$. 

\noindent Then the following are equivalent:

\begin{enumerate}
\item $\Delta(\cP)\prec_{\cM\bar\otimes \cM} \cM\bar\otimes \cQ$.
\item $\Delta(\cP)\prec_{\cM\bar\otimes \cM} \cQ\bar\otimes \cM$.

\end{enumerate}

\noindent Moreover, if $\Lambda$ is icc, $\Lambda\ca \cB$ is weakly mixing and  $\cP\subset \cL(\Lambda)$ is regular, then the above statements are also equivalent to the following:
\begin{enumerate}
\item[(3)] $\Delta(\cP)\prec_{\cM\bar\otimes \cM} \cQ\bar\otimes \cQ$.
\end{enumerate}

\end{lemma}

\noindent {\it Proof.}
We will first show that (1) and (2) are equivalent. By Kaplansky's density theorem, note that (1) does not hold if and only if there exists a sequence of unitaries $(u_n)_n\subset\mathscr U(\cP)$ such that $\|E_{\cM\bar\otimes \cQ}((1\otimes x)\Delta(u_n)(1\otimes y))\|_2\to 0$, for all $x,y\in \cM$. In this case, notice that the Fourier coefficients of $u_n=\sum_{g\in\Lambda}a_g^nv_g$ are scalars since $u_n\in \cL(\Lambda)$. Therefore, 
$$
\|E_{\cM\bar\otimes \cQ}((1\otimes x)\Delta(u_n)(1\otimes y))\|_2^2=\sum_{g\in\Lambda}|a_g^n|^2\|E_\cQ(xv_gy)\|_2^2=\|E_{\cQ\bar\otimes \cM}((x\otimes 1)\Delta(u_n)(y\otimes 1))\|_2^2\to 0,
$$
for all $x,y\in \cM$, which shows that (2) does not hold. Note that the previous formula actually shows that (1) and (2) are equivalent.

\noindent For proving the moreover part, we only have to show that (1) implies (3). Since $\Lambda$ is icc and $\Lambda\ca \cB$ is weakly mixing, it is a standard computation to check that $\Delta(\cL(\Lambda))'\cap (\cM\bar\otimes \cM)=\mathbb C1.$ Since $\cP\subset \cL(\Lambda)$ is regular, it follows by \cite[Lemma 2.4(2)]{DHI16} that $\Delta(\cP)\prec_{\cM\bar\otimes \cM}^s \cM\bar\otimes \cQ$ and $\Delta(\cP)\prec_{\cM\bar\otimes \cM}^s \cQ\bar\otimes \cM$. Finally, notice that $\cQ\bar\otimes \cM, \cM\bar\otimes \cQ\subset \cM\bar\otimes \cM$ form a commuting square, so we can apply \cite[Proposition 2.5]{Dr19a} and derive that $\Delta(\cP)\prec_\cM^s \cQ\bar\otimes \cQ$.
\hfill$\blacksquare$

\noindent {\bf Remark.} Note that the moreover part in the previous theorem holds in the case $\cB=\mathbb C1$ as well.

\begin{lemma}\label{L:reduce}
Let $\Lambda\overset{\rho}{\curvearrowright} \cB_0^\Lambda$ be a Bernoulli action and denote by $\cM=\cB_0^\Lambda\rtimes\Lambda$ the associated von Neumann algebra. Let $\Delta: \cM\to \cM\bar\otimes \cL(\Lambda)$ be the $\ast$-homomorphism given by $\Delta(b)=b\otimes 1$, for any $b\in \cB_0^\Lambda$ and $\Delta(v_g)=v_g\otimes v_g$, for any $g\in\Lambda$. Let $\cP,\cQ\subset \cL(\Lambda)$ be von Neumann subalgebras. 

\noindent  If $\Delta(\cP)\prec_{\cM\bar\otimes \cM} \cL(\Lambda)\bar\otimes \cQ$, then $\Delta(\cP)\prec_{\cL(\Lambda)\bar\otimes \cL(\Lambda)} \cL(\Lambda)\bar\otimes \cQ$.

\end{lemma}

\noindent {\it Proof.} Assuming the contrary, one can find a sequence of unitaries $(u_n)_n\subset\mathscr U(\cP)$ such that 
\begin{equation}\label{eq10}
\|E_{\cL(\Lambda)\bar\otimes \cQ}((1\otimes x)\Delta(u_n)(1\otimes y))\|_2\to 0, \text{ for all } x,y\in \cL(\Lambda).
\end{equation}

\noindent Observe that the Fourier coefficients of $u_n=\sum_{g\in\Lambda}a_g^nv_g$ are scalars as $u_n\in \cL(\Lambda)$. Therefore, to obtain a contradiction, it suffices to show that
\begin{equation}\label{eq11}
\|E_{\cL(\Lambda)\bar\otimes \cQ}((x_0\otimes x)\Delta(u_n)(y_0\otimes y))\|_2\to 0, \text{ for all }x_0,y_0\in \cB_0^\Lambda \text{ and }x,y\in \cM.
\end{equation}
Moreover, it is enough to consider $x=av_k, y=bv_h$, for some $k,h\in\Lambda$ and $a,b\in B_0^\Lambda$. If $a,b,x_0,y_0\in \mathbb C1$, then we are done by \eqref{eq10}. Hence, we can assume that $a\in  \cB_0^F, b\in  \cB_0^G,x_0\in \cB_0^H,y_0\in B_0^I$, where $F,G,H,I\subset\Lambda$ are some finite subsets and at least one of $a,b,x_0,y_0$ has trace zero. Without any loss of generality, assume that $\tau(a)=0$. Since $\|E_\cQ(a\rho_g(b)v_g)\|_2\leq \|E_{\cL(\Lambda)}(a\rho_g(b)v_g)\|_2=|\tau(a\rho_g(b))|$, for any $g\in\Lambda$, we have
\begin{equation*}
\begin{split}
\|E_{\cL(\Lambda)\bar\otimes \cQ}((x_0\otimes x)\Delta(u_n)(y_0\otimes y))\|_2^2&=
\sum_{g\in\Lambda}|a_g^n|^2|\tau(x_0\rho_g(y_0))|^2\|E_{\cQ}(a\rho_{kg}(b)v_{kgh})\|_2\\
&\leq \sum_{\{g\in\Lambda\,:\,F\cap kgG \neq\emptyset\}}|a_g^n|^2 |\tau(x_0\rho_g(y_0))|^2 |\tau(a\rho_{kg}(b))|^2.\\
\end{split}
\end{equation*}
Note that \eqref{eq10} implies that $a_g^n\to 0$, for any $g\in\Lambda$. As the last sum is a finite, this shows \eqref{eq11}.
\hfill$\blacksquare$

\noindent {\bf Proof of Proposition \ref{P:bi-exact}.} Assume by contradiction that there exists a subset $S\subset \{1,\dots,n\}$ that has at most $k-1$ elements satisfying $\cL(\Sigma)\prec_\cM\cL(\Gamma_S)$. We perform the following construction. Let $\Lambda\ca\cB$ be any Bernoulli action with abelian base and denote $\tilde\cM=\cB\rtimes\Lambda.$ Let $\Delta:\tilde \cM\to\tilde\cM\bar\otimes \cL(\Lambda)$ be the $*$-homomorphism given by $\Delta(b)=b\otimes 1$, for any $b\in \cB$ and $\Delta(v_g)=v_g\otimes v_g$, for any $g\in\Lambda$ as in \cite{PV09}.

\noindent Denote $\cP=\bar\otimes_{1\leq j\leq k} \cP_j$. The assumption implies that $\cP\prec_\mathcal M \cB\rtimes\Sigma.$ Since $\Lambda\ca\cB$ is free and mixing, we get that $(\cB\rtimes\Sigma)'\cap \mathcal M=\mathbb C 1$. By using \cite[Lemma 2.3]{Dr19b}, we get that $\Delta(\cP)\prec_{\tilde\cM\bar\otimes \tilde\cM} \Delta(\cB\rtimes\Sigma)z$, for any non-zero projection $z\in \Delta(\cB\rtimes\Sigma)'\cap (\tilde\cM\bar\otimes\tilde\cM)$. On the other hand, since $\Delta(\cB\rtimes\Sigma)\subset\tilde\cM\bar\otimes \cL(\Sigma)$, it follows by our assumption that $\Delta(\cB\rtimes\Sigma)\prec_{\tilde\cM\bar\otimes\tilde\cM} \tilde\cM\bar\otimes \cL(\Gamma_S)$. Therefore, by applying \cite[Lemma 2.4(2)]{Dr19b}, we get that $\Delta(\cP)\prec_{\tilde\cM\bar\otimes\tilde\cM} \tilde\cM\bar\otimes \cL(\Gamma_S)$. Using Lemma \ref{L:dual}, we get that $\Delta(\cP)\prec_{\tilde\cM\bar\otimes\tilde\cM} \cL(\Gamma_S)\bar\otimes \cL(\Gamma_S)$ and by Lemma \ref{L:reduce} we deduce that $\Delta(\cP)\prec_{\cM\bar\otimes \cM}\cM\bar\otimes \cL(\Gamma_S)$. By applying once again Lemma \ref{L:dual}, it follows that $\Delta(\cP)\prec_{ \cM\bar\otimes \cM} \cL(\Gamma_S)\bar\otimes \cL(\Gamma_S)$.

 Since $\Gamma$ is icc, it follows from a direct computation that $\Delta(\cM)'\cap (\cM\bar\otimes \cM)=\mathbb C1$. Since $\cM=\cP_0\bar\otimes\cP$, we get that $\mathcal{Z}(\Delta(\cP)'\cap (\cM\bar\otimes \cM))\subset \Delta(\cM)'\cap (\cM\bar\otimes \cM)$, which shows that $\Delta(\cP)'\cap (\cM\bar\otimes \cM)$ is a II$_1$ factor. We now apply \cite[Proposition 12]{OP03} and obtain a decomposition 
 $$\cM\bar\otimes \cM=\cL(\Gamma_S)\bar\otimes \cL(\Gamma_{\widehat S})\bar\otimes \cL(\Gamma_S)^t\bar\otimes \cL(\Gamma_{\widehat S})^{1/t},$$ 
 a positive number $t>0$ and a unitary $u\in\mathscr U(\cM\bar\otimes \cM)$ such that $u\Delta(\cP)u^*\subset \cL(\Gamma_S)\bar\otimes \cL(\Gamma_S)^t$. Next, since $\Delta(\cP_k)$ is non-amenable, we use Theorem \ref{Th:bo} and Lemma \ref{L:dual} and obtain an element $s_1\in S$ such that $\Delta(\bar\otimes_{1\leq i\leq k-1}\cP_i)\prec_{\cM\bar\otimes \cM} \cL(\Gamma_{S\setminus\{s_1\}})\bar\otimes \cL(\Gamma_{S\setminus\{s_1\}})$. By proceeding by induction, we get that $\Delta(\cP_1)\prec_{\cM\bar\otimes \cM}1\otimes 1$, showing that $P_1$ is not diffuse, contradiction.
\hfill$\blacksquare$

\section{A class of groups of Ioana-Popa-Vaes}\label{ipvclass}

\noindent Following \cite{IPV10}, we recall that the class $\mathcal {IPV}$ consists of all generalized wreath product groups $\G=A \wr_{I} G$ which satisfy the following properties:

\begin{enumerate}
    \item $A \cong \mathbb{Z}_2$ \text{ or } $\mathbb{Z}_3$; 
    \item $G$ is an icc  nonamenable  bi-exact group that contains an infinite property (T) normal subgroup;  
    \item The set $I = G/K$ on which $G$ acts is the set of left cosets with respect to an infinite amenable malnormal subgroup $K < G$.
    
\end{enumerate}
\noindent Concrete examples of groups in $\mathcal{IPV}$ can be obtained by considering various classes of groups intensively studied in geometric group theory. Below are two such families of examples: \begin{enumerate} \item [i)] $G$ is any torsion free, icc hyperbolic property (T) group (e.g. an uniform lattice in $Sp(n,1)$, $n\geq 2$) and $K \leqslant G$ is any infinite maximal amenable subgroup; 
\item [ii)] $G$ is any torsion free, icc, property (T) group that is hyperbolic relative to a family of amenable subgroups $\{H_1,H_2,...,H_n\}$ (see \cite[Theorem 1.1]{AMO} and \cite[Lemma 4.2(2)]{AMO}) and $K =H_i$, for some $i$. 
\end{enumerate}

\noindent Next, we record several properties for groups that belong to class $\mathcal{IPV}$ that will be useful in the next sections. 
\begin{theorem}\emph{\cite[Theorem 3.4.14]{De18}}\label{biexactwreathprod}
Any group in $\mathcal{IPV}$ is bi-exact. 
\end{theorem}

\noindent {\it Proof.}
Since the action $G \ca G/K$ is by translation for any $j = hK \in G/K$, its stabilizer satisfies ${\rm Stab}_{G}(j)=hKh^{-1}$. As $K$ is amenable so are its conjugates and hence ${\rm Stab}_G(j)$ is amenable. Finally, since $A$ is amenable then \cite[Theorem 3.4.14]{De18} implies that $\G$ is bi-exact. \hfill$\blacksquare$

\begin{theorem}\emph{\cite[Theorem 8.4]{IPV10}}
\label{ipv10} Let $\G \in \mathcal{IPV}$ and let $t > 0$. Assume that $\La$ is an arbitrary group such that there exists a $*$-isomorphism $\phi : \cL(\G) \to \cL(\La)^t$. Then $t=1$ and there exist $\delta\in {\rm Isom} (\G,\La)$,  $\omega\in {\rm Char}(\G)$ and $w \in \cU(\cL(\La))$ such that  $\phi={\rm ad} (w)\circ \Psi_{\omega, \delta},$ where $\Psi_{\omega, \delta}(u_g)=\omega(g) v_{\delta(g)},$ for any $g\in\Gamma$.

\end{theorem}

\section{$W^*$-superrigidity for product groups}

\noindent In the first part of the section we prove Theorem \ref{mainthm1} (see Theorem \ref{prodbiexactrig}) and therefore generalize the main results from \cite{CdSS15}. The technology that we use is slightly different from the one in \cite{CdSS15}, resembling more the methods developed in \cite{DHI16,Dr19a,Dr19b}.

\noindent In the second part we use the product rigidity in combination with other prior results \cite{IPV10} to show that any direct product of groups in class $\mathcal {IPV}$ is $W^*$-superrigid (see Corollary \ref{superrigidproducts}).

\noindent One of the crucial ingredients in the proof of Theorem \ref{prodbiexactrig} is an ultrapower technique \cite{Io11}, which we recall in the following form. This result is essentially contained in the proof of \cite[Theorem 3.1]{Io11} (see also \cite[Theorem 3.3]{CdSS15}), but the statement that we will use is a particular case of \cite[Theorem 4.1]{DHI16}.

\begin{theorem}\emph{\cite{Io11}}\label{Th:ultrapower}
Let $\Lambda$ be a countable icc group and denote by $\cM=\cL(\Lambda)$. Let $\Delta:\cM\to \cM\bar\otimes \cM$ be the $*$-homomorphism given by $\Delta (v_\lambda)=v_\lambda\otimes v_\lambda$, for all $\lambda\in\Lambda.$ Let $\cP,\cQ\subset \cM$ be von Neumann subalgebras such that $\Delta(\cP)\prec_{\cM\bar\otimes \cM}\cM\bar\otimes \cQ$.

\noindent Then there exists a decreasing sequence of subgroups $\Sigma_k<\Lambda$ such that $\cP\prec_\cM \cL(\Sigma_k)$, for every $k\ge 1$, and $\cQ'\cap \cM\prec_\cM \cL(\cup_{k\ge 1} C_\Lambda(\Sigma_k)).$
\end{theorem}


\noindent We are now ready to present the product rigidity result.
\begin{theorem} \label{prodbiexactrig}
	 Let $\G=\G_1\times\dots\times \G_n$ be a product of $n\ge 1$ icc, non-amenable, bi-exact groups and denote by $\cM=\El(\G)$. Let $\La$ be any countable group and $t>0$ such that $\emm^t=\El(\La)$. Then there exist a direct product decomposition $\La=\La_1\times\dots\times \La_n$, some scalars $t_1,\dots,t_n>0$ with $t_1\cdots t_n=t$, and a unitary $u\in \emm^t$ such that $u\El( \La_i)u^*=\El(\G_i)^{t_i}$, for any $1\leq i\leq n$.	
	 \end{theorem}
 
\noindent {\it Proof.} Without any loss of generality we can assume $t=1$, since the general case does not hide any technical difficulties.
Let $\Delta:\cM\to \cM\bar\otimes \cM$ be the $*$-embedding given by $\Delta(v_\lambda)=v_\lambda\otimes v_\lambda$, for any $\lambda\in\Lambda$ as in \cite{PV09}. First we prove the following

\begin{claim} For any $1\leq i\leq n$, we have $\Delta(\cL(\Gamma_{\hat i}))\prec_{\cM\bar\otimes \cM} \cM\bar\otimes \cL(\Gamma_{\hat j})$ for some $1\leq j\leq n.$
\end{claim}
\noindent {\it Proof of Claim 1.} Fix $1\leq i\leq n$. Since $\Delta(\cL(\Gamma_i))$ and $\Delta(\cL(\Gamma_{\hat i}))$ are commuting non-amenable subalgebras of $\cL(\Gamma)\bar\otimes \cL(\Gamma)$, it follows by Theorem \ref{Th:bo} that there exists $1\leq j\leq n$ such that $\Delta(\cL(\Gamma_{\hat i}))\prec_{\cM\bar\otimes \cM} \cM\bar\otimes \cL(\Gamma_{\hat j})$ or $\Delta(\cL(\Gamma_{\hat i}))\prec_{\cM\bar\otimes \cM}  \cL(\Gamma_{\hat j})\bar\otimes \cM$. The claim follows by using Lemma \ref{L:dual}; alternatively, one could use the flip automorphism $\sigma$ of $\cM\bar\otimes \cM$ since $\sigma\circ\Delta=\Delta$.
\hfill$\square$

\noindent Theorem \ref{Th:ultrapower} combined with Claim 1 imply that there exists a subgroup $\Sigma_i<\Lambda$ with non-amenable centralizer $C_{\Lambda}(\Sigma_i)$ such that $\cL(\Gamma_{\hat i})\prec_\cM \cL(\Sigma_i)$, for any $1\leq i\leq n.$ Moreover, we show next that the following holds:

\begin{claim}\label{backforthint1} For any $1\leq k\leq n$, we have $\cL(\Gamma_{\hat k})\prec_\cM^s \cL(\Sigma_k)$ and $\cL(\Sigma_k)\prec_\cM^s \cL(\Gamma_{\hat k})$.
\end{claim}

\noindent {\it Proof of Claim \ref{backforthint1}.} We will show the claim only for $k=1$ as the other cases are similar. First, we notice that since $\cL(\Gamma_{\hat 1})\subset \cM$ is regular then by  \cite[Lemma 2.4(2)]{DHI16}  we have that $\cL(\Gamma_{\hat 1})\prec_\cM^s \cL(\Sigma_1)$.  Next, we  show the second intertwining relation. Using \cite[Lemma 2.4(3)]{DHI16} there is a maximal projection $e_i\in \mathcal Z(\cL(\Sigma_1)'\cap M)$, possibly the zero projection, such that
\begin{equation}\label{eq1}
\cL(\Sigma_1)e_i\prec_\cM^s \cL(\Gamma_{\hat i}), \text{ for any  }1\leq i\leq n.
\end{equation}

\noindent Remark that $e:=e_1\vee\dots\vee e_n=1$. Indeed, otherwise the projection $f:=1-e\in \mathcal Z(\cL(\Sigma_1)'\cap \cL(\Lambda))$ is non-zero, and hence, $\cL(C_\Lambda(\Sigma_1
))f$ is non-amenable by Lemma \ref{L:nonamenable}. Thus by Theorem \ref{Th:bo} there exists $1\leq i\leq n$ such that $\cL(\Sigma_1)f\prec_\cM \cL(\Gamma_{\hat i})$. Since $f\leq 1-e_i$ this contradicts the maximality of $e_i$.

\noindent We continue by showing that for every $1\leq i\leq n$ we have either $e_i=0$ or $e_i=1$. 
Denote by $\Omega_0$ the set of all $\la\in\Lambda$ such that $\cL(\Gamma_{\hat 1})\prec_\cM \cL(\la\Sigma_1 \la^{-1}\cap\Sigma_1)$. First, we prove that 
\begin{equation}\label{commute}
v_\la e_iv_\la ^*=e_i, \text{ for all } \la\in\ \Omega_0 \text{ and } 1\leq i\leq n.
\end{equation}

\noindent If \eqref{commute} does not hold, then one can find $\la\in \Omega_0$ and $1\leq i\leq n$ such that $v_\la e_iv_\la^*\neq e_i.$ Hence, there is $j\neq i$ such that $v_\la e_iv_\la^*e_j\neq 0.$
By \eqref{eq1}, we get that $\cL(\la \Sigma_1\la^{-1}\cap \Sigma_1)v_\la e_iv_\la^*\prec_\cM^s \cL(\Gamma_{\hat i})$ and $\cL(\la\Sigma_1\la^{-1}\cap \Sigma_1)e_j\prec_\cM^s \cL(\Gamma_{\hat j})$. Using Lemma \ref{L:bigger}, we get non-zero projections $f_i,f_j\in \mathcal Z(\cL(\la\Sigma_1\la^{-1}\cap \Sigma_1)'\cap \cM)$ with $v_\la e_iv_\la^*\leq f_i$ and $e_j\leq f_j$ such that $\cL(\la\Sigma_1\la^{-1}\cap \Sigma_1)f_i\prec_M^s \cL(\Gamma_{\hat i})$ and $\cL(\la\Sigma_1\la^{-1}\cap \Sigma_1)f_j\prec_\cM^s \cL(\Gamma_{\hat j})$. Since $v_\la e_iv_\la^*e_j\neq 0$, we get that $f_0:=f_if_j\in \mathcal Z(\cL(\la\Sigma_1\la^{-1}\cap \Sigma_1)'\cap \cM)$ is a non-zero projection. By applying \cite[Lemma 2.8(2)]{DHI16}, we get that $\cL(\la\Sigma_1\la ^{-1}\cap \Sigma_1)f_0\prec_\cM^s \cL(\Gamma_{\widehat {\{i,j\}}})$. Finally, using Proposition \ref{P:bi-exact} and the fact that $\la\in\Omega_0$, we get a contradiction.
Hence, relation \eqref{commute} must hold.

\noindent If we let $\Omega$ be the subgroup generated by $\Omega_0$, we deduce that $v_\la e_iv_\la^*=e_i, \text{ for all } \la\in\ \Omega \text{ and } 1\leq i\leq n.$
By applying Lemma \ref{HPV}, we get that $\Omega<\Lambda$ has finite index. Since $\Lambda$ is icc, it follows that the set $\{\la \sigma \la^{-1}\,:\,\la\in\Omega\}$ is infinite for any $\sigma\in\Lambda\setminus\{1\}$. A standard computation reveals that \eqref{commute} implies $e_i\in\mathbb C1$. Since $e_i$ is a projection, it follows that $e_i=0$ or $e_i=1.$
Now, using $e_1\vee\dots\vee e_n=1$, one can find $i$ such that $\cL(\Sigma_1)\prec_\cM^s \cL(\Gamma_{\hat i})$. Since $\cL(\Gamma_{\hat 1})\prec_\cM^s \cL(\Sigma_1)$, then \cite[Lemma 3.7]{Va08} implies that $\cL(\Gamma_{\hat 1})\prec_\cM^s \cL(\Gamma_{\hat i})$. Since $\Gamma_i$ is an infinite group, it follows that $i=1$, thus ending the claim.
\hfill$\square$

\begin{claim}\label{backforthint2} There exists a subgroup $\Sigma_0<\Lambda$ such that $\cL(\Sigma_0)\prec_\cM^s \cL(\Gamma_1)$ and $\cL(\Gamma_1)\prec_\cM^s \cL(\Sigma_0) $.
\end{claim}

\noindent {\it Proof of Claim 3.} From Claim \ref{backforthint1} we have that $\cL(\Gamma_{\hat 2})\prec_\cM^s \cL(\Sigma_2)$ and $\cL(\Gamma_{\hat 3})\prec_\cM^s \cL(\Sigma_3)$. Using \cite[Lemma 2.7]{Va10} we find an element $\la_3\in\Lambda$ such that $\cL(\Gamma_{\widehat {\{2,3\}}})\prec_\cM^s \cL(\Sigma_2\cap \la_3\Sigma_3 \la_3^{-1})$. From Claim \ref{backforthint1} and \cite[Lemma 2.8(2)]{DHI16} we deduce that $\cL(\Sigma_2\cap \la_3\Sigma_3 \la_3^{-1})\prec_\cM^s \cL(\Gamma_{\widehat {\{2,3\}}})$.

\noindent Proceeding by induction for every $j\ge 2$, there exists $\la_j\in\Lambda$ such that $\cL(\Gamma_{\widehat {\{2,\dots, j\}}})\prec_\cM^s \cL(\Sigma_2\cap \la_3\Sigma_3 \la_3^{-1}\cap\dots \cap \la_j\Sigma_j \la_j^{-1})$ and $ \cL(\Sigma_2\cap \la_3\Sigma_3 \la_3^{-1}\cap\dots \cap \la_j\Sigma_j \la_j^{-1})\prec_\cM^s \cL(\Gamma_{\widehat {\{2,\dots, j\}}})$. Since $\widehat {\{2,\dots, n\}}=\{1\}$, Claim \ref{backforthint2} follows by taking $\Sigma_0=\Sigma_2\cap \la_3\Sigma_3 \la_3^{-1}\cap\dots \cap \la_n\Sigma_n \la_n^{-1}$. 
\hfill$\square$

\noindent Using the Claims \ref{backforthint1}-\ref{backforthint2} in combination with \cite[Theorem 6.1]{DHI16}, one can find a product group decomposition $\Lambda=\Lambda_1\times\Lambda'_1$, a tensor decomposition $\cM^t=\cL(\Gamma_1)^{t_1}\bar\otimes \cL(\Gamma_{\hat 1})^{1/t_1}$, for some scalar $t_1>0$, and a unitary $u\in \mathscr U(\cM)$ such that $u\cL(\Lambda_1)u^*=\cL(\Gamma_1)^{tt_1}$ and $u\cL(\Lambda'_1)u^*=\cL(\Gamma_{\hat 1})^{1/t_1}$. Since $\Gamma_{\hat 1}$ is a product of $n-1$ non-amenable bi-exact icc groups, we derive the conclusion by a standard induction argument. \hfill$\blacksquare$


\noindent Now combining our prior product rigidity results for bi-exact groups together with the superrigidity results from \cite{IPV10} we derive many examples of  $W^*$-superrigid groups of product type.

\begin{corollary}\label{superrigidproducts}
Let $\G_1,\G_2,\ldots,\G_k \in \mathcal{IPV}$ and denote by $\G = \G_1 \times \G_2 \times \cdots \times \G_k$. Let $t>0$ and assume that $\La$ is an arbitrary group such that there exists a $*$-isomorphism $\theta: \El(\G) \to \El(\La)^t$. Then $t=1$ and there exist $\delta \in {\rm Isom}(\G,\La)$, $ \omega \in {\rm Char}(\G)$ and  $u \in \mathcal{U}(\cL(\La))$ such that $\theta = {\rm ad}(u) \circ \Psi_{\omega,\delta}$. 
\end{corollary}

\noindent {\it Proof.} Using Theorems \ref{biexactwreathprod} and \ref{prodbiexactrig} there exist a $k$-folded product decomposition $\La = \La_1 \times ... \times \La_k$, scalars  $t_1,\ldots,t_k>0$ with $t_1 t_2 \cdots t_k = t$ and a unitary $w \in \mathcal{U}(\cL(\La))$ such that for every $1\leq i\leq k$ we have \[
w\theta(\cL(\G_i))^{t_i}w^* = \cL(\La_i). 
\]
Since $\G_i \in \mathcal{IPV}$ then Theorem \ref{ipv10} further implies that $t_i =1$ and there exist $  \omega_i \in {\rm Char}(\G_i)$, $\delta_i \in {\rm Isom}(\G_i,\La_i)$ and $w_i \in \mathcal U(\cL(\La_i))$ such that $w \theta(u_{\g_{i}})w^* = \omega_i(\g_{i}) w_i v_{\delta_{i}(\g_{i})}w_i^*$ for all $ \g_i \in \G_i$. Thus $t=1$ and letting $\omega = \prod_{i=1}^{k}{\omega_{i}} $,  $\delta = \prod_{i=1}^{k}{\delta_{i}}$ and $u=w^*\prod^k_{i=1}w_i$ we get the desired conclusion. 
\hfill$\blacksquare$

\section{A class of iterated amalgamated free products and HNN extension groups}\label{treesuperrigid}

\noindent In this section we present several properties of groups that belong to class $\mathscr D$ and  their associated von Neumann algebras. 


     \begin{lemma}\label{controlquasinorminhnn}
     Let $G= {\rm HNN
     }(K,\Sigma,\varphi)$ be an HNN-extension where $\Sigma<K$ are groups and $\varphi:\Sigma \rightarrow K$ is a monomorphism. Then the following hold: 
     
     \begin{enumerate}\item ${\rm QN}^{(1)}_G(\Sigma)=\Sigma$ if and only if ${\rm QN}^{(1)}_K(\Sigma)=\Sigma$, ${\rm QN}^{(1)}_K(\varphi(\Sigma))=\varphi(\Sigma)$ , $[\Sigma: \Sigma\cap g\varphi(\Sigma)g^{-1}]=\infty$ and  $[\varphi(\Sigma) : \varphi(\Sigma) \cap g \Sigma g^{-1}]=\infty$ for all $g\in K$. 
     \item Under the assumptions of 1.\   we have that $K\cap gKg^{-1}\leqslant \Sigma\cap \varphi(\Sigma)$ for every $g\in G\setminus K
     $. 
     \end{enumerate}
     \end{lemma}
     
     \noindent {\it Proof.}
     1. First we prove the forward implication. Since the first assertion follows trivially we will only justify the second one. Note that we can write $G=<K,t \, : \, \varphi(\sigma)=t^{-1}\sigma t, \text{ for any }\sigma\in \Sigma>$. 
          Assume there is $g\in K$ so that $[\Sigma :\Sigma \cap g\varphi(\Sigma)g^{-1}]<\infty$. Using the HNN relation this implies that $[\Sigma :\Sigma \cap gt\Sigma (gt)^{-1}]<\infty$ and hence $gt\in {\rm QN}_G^{(1)}(\Sigma)=\Sigma$. Thus $t\in g^{-1}\Sigma<K$, a contradiction.

     \noindent Now we prove the reverse implication. Fix $g\in {\rm QN}^{(1)}_G(\Sigma)$. Thus one can find finite index subgroups $\Sigma_1,\Sigma_2\leqslant \Sigma$ and a group isomorphism $\theta: \Sigma_1\rightarrow \Sigma_2$ such that \begin{equation}\label{conjfiniteindexsubgroups}\theta(h) =g hg^{-1}\text{ for all }h\in \Sigma_1.\end{equation}
     
     \noindent Using Britton's lemma we can write $g$ in reduced form, i.e.\ $g=g_0 t^{\varepsilon_1} g_1 t^{\varepsilon_2 }... g_{n-1}t^{\varepsilon_n} g_n$ where $g_i\in K$, $\varepsilon_i\in \{-1,1\}$ and the word $g$ does not contain any substring of the form $tht^{-1}$ for $h\in \Sigma$ or $t^{-1}k t$ for $k\in \varphi(\Sigma)$. Using this together with equation \eqref{conjfiniteindexsubgroups} we have that \begin{equation}\label{analysisword}1= \theta(h)gh^{-1}g^{-1}= \theta(h)g_0 t^{\varepsilon_1} g_1 t^{\varepsilon_2 }... g_{n-1}t^{\varepsilon_n} g_n h^{-1} g_n^{-1} t^{-\varepsilon_n} g^{-1}_{n-1} t^{-\varepsilon_{n-1}}... t^{-\varepsilon_1}g_0^{-1} \text{ for every }h\in \Sigma_1.\end{equation} Therefore using the Britton's normal form we have two cases to analyze: either I) $g_n hg_n^{-1}\in \Sigma$ for all $h\in \Sigma_1$ and $\varepsilon_n = 1$ or  II) $g_n hg_n^{-1}\in \varphi(\Sigma)$ for all $h\in \Sigma_1$ and $\varepsilon_n = -1$.  
     
     \noindent Assume $n\geq 2$. 
     If we are in case I) then we see that $g_n \Sigma_1g_n^{-1}<\Sigma$ and hence $g_n \Sigma_1g_n^{-1}\cap \Sigma= g_n \Sigma_1 g_n^{-1}$. Thus $[\Sigma : \Sigma \cap g_n^{-1} \Sigma g_n]=[g_n\Sigma g_n^{-1}: g_n\Sigma g_n^{-1} \cap  \Sigma ]\leq [g_n\Sigma g_n^{-1}: g_n\Sigma_1 g_n^{-1} \cap  \Sigma ] =[g_n\Sigma g_n^{-1}: g_n\Sigma g_n^{-1}  ] =[\Sigma:\Sigma_1]<\infty$. From the assumptions this implies that $g_n\in \Sigma$. In particular we have that $t^{\varepsilon_n} g_n \Sigma_1  g^{-1}_n t^{-\varepsilon_n}= \varphi(\Sigma_1')$ where $\Sigma_1'= g_n \Sigma_1    g_n^{-1}\leqslant \Sigma$ is a finite index subgroup. Thus equation \eqref{analysisword} again implies that either Ia) $g_{n-1} \varphi(\Sigma_1') g
     ^{-1}_{n-1}< \Sigma$ and $\varepsilon_{n-1}=1$ or Ib) $g_{n-1} \varphi(\Sigma_1') g^{-1}_{n-1}< \varphi(\Sigma)$ and $\varepsilon_{n-1}=-1$.
     
     \noindent Assume sub-case Ia). Thus  $[\varphi(\Sigma) : \varphi(\Sigma) \cap g_{n-1}^{-1} \Sigma g_{n-1}]=[g_{n-1}\varphi(\Sigma) g_{n-1}^{-1}: g_{n-1}\varphi(\Sigma) g_{n-1}^{-1} \cap  \varphi(\Sigma) ]\leq [g_{n-1}\varphi(\Sigma) g_{n-1}^{-1}: g_{n-1}\varphi(\Sigma_1') g_{n-1}^{-1} \cap  \varphi(\Sigma) ] =[g_{n-1}\varphi(\Sigma) g_{n-1}^{-1}: g_{n-1}\varphi(\Sigma_1') g_{n-1}^{-1}] =[\varphi(\Sigma):\varphi(\Sigma'_1)]<\infty$. However this contradicts the assumptions so this case cannot hold.
     
     \noindent Now assume sub-case Ib). Then we see that $[\varphi(\Sigma) : \varphi(\Sigma) \cap g_{n-1}^{-1} \varphi(\Sigma) g_{n-1}]=[g_{n-1}\varphi(\Sigma) g_{n-1}^{-1}: g_{n-1}\varphi(\Sigma) g_{n-1}^{-1} \cap  \Sigma ]\leq [g_{n-1}\varphi(\Sigma) g_{n-1}^{-1}: g_{n-1}\varphi(\Sigma_1') g_{n-1}^{-1} \cap  \Sigma ] =[g_{n-1}\varphi(\Sigma) g_{n-1}^{-1}: g_{n-1}\varphi(\Sigma_1') g_{n-1}^{-1}] =[\varphi(\Sigma):\varphi(\Sigma'_1)]<\infty$. Using the assumptions we infer that $g_{n-1}\in \varphi(\Sigma)$. However this together with the previous relations imply that $t^{\varepsilon_{n-1}} g_{n-1}t^{\varepsilon_n}= t^{-1} g_{n-1}t$ which contradicts that the word $g$ is reduced. So sub-case Ib) is impossible as well. 
     
     \noindent Altogether these show that case I is impossible. Proceeding in a similar manner one can show case II is impossible as well. In conclusion we must have $n\leq 1$. 
     
     \noindent Next assume  $n=1$. Also assume we are in case I. Proceeding as before we must have that $g_1\in \Sigma$ and using equation \eqref{analysisword} we see that $g_0 \Sigma_2 g_0^{-1}=g_0^{-1} \theta(\Sigma_1) g_0^{-1}= t g_1 \Sigma_1 g_1^{-1} t^{-1}=\varphi( g_1 \Sigma_1 g_1^{-1})$. Since $[\Sigma:\Sigma_1]<\infty$ we must have $[\Sigma:g_1\Sigma_1g_1^{-1}]<\infty$. Combining this with the previous relation this further entail that $[ g_0\Sigma g_0^{-1} :g_0\Sigma g_0^{-1}\cap  \varphi(\Sigma ) ]\leq [ g_0\Sigma g_0^{-1} :g_0\Sigma_2 g_0^{-1}\cap  \varphi(\Sigma ) ]=[ \Sigma  :\Sigma_2  ] <\infty$  thus contradicting the hypothesis assumptions. In a similar way case II also leads to a contradiction. Thus $n=0$ and hence equation \eqref{conjfiniteindexsubgroups} together with the hypothesis imply that $g=g_0\in \Sigma$, as desired.

     \noindent Part 2. follows by similar computations. We leave the details to the reader. \hfill$\blacksquare$

     \begin{theorem}\label{theorem.classD.qntrivial} If $G\in \mathscr D_i$, for some $i\ge 1$, and $a_i(G)=\Sigma$, then the following hold:
     \begin{enumerate}
         \item ${\rm QN}^{(1)}_{G}(\Sigma)=\Sigma$; 
         \item ${\rm vC}_G(\Sigma)=1$.
     \end{enumerate}
      \end{theorem}
     
     \noindent {\it Proof.} 1. This follows directly from the quasi-normalizers conditions from class $\mathscr{D}$ by applying \cite[Lemma 2.11]{CI17} and Lemma \ref{controlquasinorminhnn} inductively.

     \noindent 2. Let $g\in {\rm vC}_G(\Sigma)$. Thus there exists a finite index subgroup $\Sigma_0\leqslant \Sigma$ so that $g\in C_G(\Sigma_0)$; in particular $g \in {\rm QN}_G(\Sigma_0)$. However using the finite index condition one can easily check that ${\rm QN}_G(\Sigma_0)={\rm QN}_G(\Sigma)$ and combining with the first part we conclude that $g\in \Sigma$. In particular, this shows that ${\rm vC}_G(\Sigma) \subset \Sigma$ and since $\Sigma
     $ is icc we conclude that $vC_G(\Sigma)=1$, as desired.\hfill$\blacksquare$

     \begin{proposition}\label{P:amenable radical}
     Let $G\in \mathscr D_i$ with $i\ge 1$ and write  $a_i(G)=\Sigma$. Assume that that there exist $g_1,\dots, g_k\in G$ such that $\cap _{i=1}^k g_i\Sigma g_i^{-1}$ is finite.
     
    Then $\cL(G)$ does not admit a diffuse amenable regular von Neumann subalgebra. 
          In particular, $G$ has trivial amenable radical.
     \end{proposition}

{\it Proof.} Let $\cA\subset \cL(G)$ be an amenable regular von Neumann subalgebra. First, we show that $\cA \prec_{\cL(G)} {\cL(\Sigma)}$. From the definition, $G$ is either (1) an amalgamated free product $G=H_1\ast_\Sigma H_2$ with $H_1,H_2\in \mathscr D_{i-1}$ and $\Sigma$ amenable or (2) an HNN-extension $G= {\rm HNN}(H_1,\Sigma,\phi)$ with $H_1\in \mathscr D_{i-1}$ and $\Sigma$ amenable. By applying \cite[Theorem A and Theorem 4.1]{Va13}, we deduce that $\cA \prec_{\cL(G)} \cL(\Sigma)$ holds in both cases. 
Next, by applying \cite[Proposition 8]{HPV11} we derive that $\cA\prec_{\cL(G)} \cL(\cap _{i=1}^k g_i\Sigma g_i^{-1})$, implying that $\cA$ is not diffuse.

For the last part of the proof notice that since $G$ is icc, the first part of the proof implies that $G$ has trivial amenable radical. \hfill$\blacksquare$

     \begin{theorem}\label{Th:commute} Let $G\in \mathscr D$ and denote by $f(G)=\{G_1,G_2,...,G_n\}$ its factor set. For every $1\leq i\leq n$ denote by $s_i\geq 2$  the integer such that $G_i =\G_1^i \times \G_2^i\times ...\times \G^i_{s_i}$, where $\G^i_j\in \mathcal{IPV}$. Denote $\cM= \cL(G\times G)$ and let $p\in \cM$ be a projection. Assume that $\cA,\cB\subseteq p\cM p$ are commuting von Neumann subalgebras which contain property (T) diffuse subalgebras $\cA_0 \subseteq \cA$ and $\cB_0 \subseteq \cB$. Also assume that $\cA
     _0\nprec_\cM \cL(G\times A)$, $\cA_0 \nprec_\cM \cL(A\times G)$, $\cB
     _0\nprec_\cM \cL(G\times A)$ and $\cB_0 \nprec_\cM \cL(A\times G)$ for any amenable subgroup $A\leqslant G$.

\noindent     Then one can find $1\leq i\leq n$ and $1\leq j\leq s_i$ such that $\cA \prec_\cM \cL (G\times \Gamma^i_{\hat j})$ or $\cA \prec_\cM \cL (\Gamma^i_{\hat j} \times G)$.
     \end{theorem}
     
\noindent     Here, we denoted by $\Gamma^i_{\hat j}$ the product group $\times_{k
     \in \hat j}\Gamma_k^i$. 
     
     \noindent {\it Proof.} If $G\in\mathscr D_0$, then the result follows from Theorem \ref{Th:bo}. Hence, we
    assume that $G\in \mathscr D_m$ with $m\ge 1$ and denote $a_m(G)=\Sigma$. Firstly, we claim that there exist integers $1\leq k,l\leq n$, a projection  $0\neq z\in (\cA \vee \cB )'\cap \cM$ and unitary $u\in \cM$ such that \begin{equation}\label{firstreduction}
         u (\cA\vee \cB) zu^*\subseteq \cL(G_k\times G_l).
     \end{equation}
     From the definition, $G$ is either an amalgamated free product $G=H_1\ast_\Sigma H_2$ with $H_i\in \mathscr D_{m-1}$ and $\Sigma$ amenable or an HNN-extension $G= {\rm HNN}(H_1,\Sigma,\phi)$ with $H_1\in \mathscr D_{m-1}$ and $\Sigma$ amenable. Thus $\cM$ is canonically  either an amalgamated free product or an HNN-extension von Neumann algebra and since $\cA_0\subseteq \cM$ is a property (T) subalgebra then using either \cite[Theorem 5.1]{IPP05} or \cite[Theorem 3.4]{FV10} we have that $\cA_0 \prec_\cM \cL(G\times H_i):=\cP_i$ for some $H_i\in \mathscr D_{k-1}$.   Thus one can find projections $a\in \cA_0$, $p\in \cP_i$, a non-zero partial isometry $v\in p\cM a$, and an injective $\ast$-isomorphism  $\phi: a\cA_0 a\rar \phi(a \cA_0 a):=\cQ\subseteq p\cP_i p$  so that $\phi(x)v=vx$ for all $x\in a\cA_0 a$. Moreover, $v^*v\in a\cA_0 a'\cap a\cM a$, $vv^*\in \cQ'\cap p\cM p$ and we can assume that the support $s(E_{\cP_i}(vv^*))=q$.

     \noindent Next, observe that $\cQ\nprec_{\cP_i} \cL (G \times \Sigma)$. Indeed, otherwise composing this intertwining with $\phi$ we would obtain that $\cA_0 \prec_\cM \cL(G\times \Sigma)$. Since $\Sigma$ is amenable this would contradict the hypothesis assumptions. Therefore, by Lemma \ref{L:ipp}, we have that $vv^*\in \cQ'\cap p\cM p\subseteq \cL(G\times H_i)$. In particular, we have $v \cA_0 v^*= \cQ vv^*\subseteq \cL(G\times H_i)$ and moreover, if $u$ is a unitary extending $v^*v$, we get that $u v^*v (\cA_0 \vee \cA_0 '\cap p\cM p ) v^*v u^*\subseteq \cL(G\times H_i)$. As $\cL(G\times H_i)$ is a factor, after perturbing $u$ to a new unitary we further get that  $u  (\cA_0 \vee \cA_0 '\cap p\cM p)z u^*\subseteq \cL(G\times H_i)$, where $z$ is the central support of $v^*v
     $ in $\cA\vee \cA'\cap p\cM p $. Thus $u  (\cA_0 \vee \cB)z u^*\subseteq \cL(G\times H_i)$ and in particular $u\cB z u^*\subseteq \cL(G\times H_i)$. From the assumptions we also see that $u\cB zu^* \nprec_{\cP_i} \cL(G\times \Sigma)$ and therefore repeating the same argument as before on control of relative commutants we get that $u (\cB'\cap z\cM z) \vee \cB z u^*\subseteq \cL(G\times H_i)$; in particular, we conclude that $u (\cA \vee \cB)z u^*\subseteq \cL(G\times H_i)$.  Now, notice that $\cA z$ and $\cB z$ are still commuting von Neumann subalgebras containing property (T) diffuse subalgebras $\cA_0 z\subseteq \cA z$ and $\cB_0 z\subseteq \cB z$. Therefore, one can repeat the same argument finitely many times so that in the end there exist $H_k, H_l \in \mathscr D_0$ a unitary still denoted by $u \in \cM$ and a non-zero projection $z\in (\cA \vee \cB)'\cap \cM$   satisfying    $u(\cA \vee \cB)z u^* \subseteq \cL(H_k \times H_l) $. However, since the elements of $\mathscr D_0$ consists of factor subgroups of $G$, the claim \eqref{firstreduction} follows.

     
   \noindent  Finally, note that since the groups $\Gamma_j^i$'s are bi-exact and $B$ is non-amenable, we can apply Theorem \ref{Th:bo} and obtain the conclusion.
     \hfill$\blacksquare$

     \noindent We remark that the result above can also be obtained using bi-exactness methods from \cite{Oz03,BO08}. Moreover, the theorem still holds under the milder assumption that the algebras $\cA$ and $\cB$ have no amenable direct summand rather  containing property (T) diffuse subalgebras. For the interested reader we also note this result can be proved either using bi-exactness methods from \cite{Oz03,Oz04,BO08} or using Popa's deformation/ rigidity theory as in \cite{CH08,FV10}. We opted for this leaner version only for the brevity of the exposition as it follows relatively easily from existing results in the literature.
     
    {\bf Some examples of amalgamated free product groups in class $\mathscr D$.}  Let $K$ be a non-elementary, torsion free, property (T) group that is hyperbolic relative to a finitely generated, icc amenable subgroup $P<K$; using \cite[Theorem 1.1]{AMO}, for any given finitely generated group $P$, such a $K$ always exists. Morever, since amenable groups are biexact, then \cite{Oz04} (or \cite[Theorem 1.1]{Oy22}) implies that $K$ is biexact. Now fix $g\in K \setminus P$ a hyperbolic element and let $B:= E(g)< K$ be the (unique) maximal elementary subgroup defined by $E(g) = \{f\in K \,|\,  f^{-1}g^n f = g^{\pm n}, \text{ for some  }n\in \mathbb N\}$. Observe that both $P$ and $B$ are amenable and malnormal in $K$. Moreover, we have that $P\cap kBk^{-1}=1$, for all $k\in K$. The last two assertions follow from \cite[Lemma 3.1]{AMO} and \cite[Theorem 1.4]{Osi06b}. We consider the generalized wreath product $H= \mathbb Z_2 \wr_{K/B} K$ and notice that $H \in \mathcal {IPV}$. Let $G=\times_{i=1}^n H\in \mathscr D_0$ with $n\ge 2$. Therefore, the following groups belong to $\mathscr D_1$:
    
    \begin{itemize}
       \item $G*_\Sigma G$, where $\Sigma=\times_{i=1}^n P<G$ is the natural direct product embedding.
        
        \item $G*_\Sigma G$, where $\Sigma<G$ is the diagonal product embedding of $P$; more precisely, $\Sigma=\{(g,\dots ,g)\,|\,g\in P\}$.
        
    \end{itemize}

    Moreover, note that these groups actually belong to $\mathscr D_1^m$.
    For examples of groups in $\mathscr D_{i}$ with $i\ge 2$, one can easily iterate by following the procedure described in the definition of class $\mathscr D$.

  {\bf Some examples of HNN extensions in class $\mathscr D$.} Let $K$ be a non-elementary, torsion free, hyperbolic group that admits a normal infinite property (T) subgroup. Assume that $B,C,D< K$ are infinite cyclic subgroups that are malnormal and satisfy   $C\cap g B g^{-1}= D \cap gBg^{-1}=C \cap g D g^{-1} =1$, for all $g\in K$. Next we briefly indicate how to build such groups using Belegradek-Osin's Rips construction in geometric group theory \cite{BO06}. Consider the free group with three generators $\mathbb F_3=\langle a,b,c\rangle$. Then using \cite[Theorem 1.1]{BO06} there exist a torsion free, property (T) group $N$ and an action by automorphisms $\mathbb F_3\ca N$ such that the corresponding semidirect product $K=N\rtimes \mathbb F_3$ is hyperbolic relative to $\mathbb F_3$. Since $\mathbb F_3$ is itself hyperbolic it follows that $K$ is hyperbolic and also torsion free. Now consider the cyclic subgroups of $K$ given by the generators of $\mathbb F_3$, $B=\langle a\rangle$   $C=\langle b\rangle$  $D=\langle  c\rangle$. Since $K$ is hyperbolic relative to $\mathbb F_3= B\ast C\ast D$ then $F_3$ is malnormal in $K$ and therefore one can check easily that $B$, $C$, $D$ satisfy the required conditions. We also mention that one can build groups $K$ with the required properties that actually have property (T) in a similar manner, but using \cite[Theorem 1.1]{AMO} instead of \cite[Theorem 1.1]{BO06}.    
     
     Now, consider the generalized wreath product $\G= \mathbb Z_2 \wr_{K/B} K$ and notice that $\G \in \mathcal {IPV}$.

     \noindent Next, we consider the canonical subgroups $\Omega = \mathbb Z_2 \wr_{K/B} C$ and $\Upsilon = \mathbb Z_2 \wr_{K/B} D$ of $\G$ and we claim that $\Omega$ is isomorphic to $\Upsilon$. Towards this, we first notice that the actions by left translations on the base sets $C\ca K/B$ and $D\ca K/B$ have trivial stabilizers. Indeed, for every $g B \in K/B$ its stabilizer in $C$ is $g Bg^{-1}\cap C$ which by assumption is trivial. Moreover these actions have (countable) infinitely many orbits which are given specifically by the double cosets $K= \bigsqcup_{g\in F_1} C g B= \bigsqcup_{g\in F_2} Dg B$,  where $|F_1|=|F_2|=\aleph_0$. Indeed, just notice that if  $F_i$ would be finite then it would imply that $K$ is boundedly generated. However this would contradict for instance \cite{Mi04} or  \cite[Therem 1.9]{Os04}. Using these observations one can see the following sequence of isomorphisms hold: $\Omega= \mathbb Z_2 \wr_{K/B} C = \oplus_ {g\in F_1} ( \oplus_{h\in Cg B} \mathbb Z_2) \rtimes C \cong  \oplus_ {g\in F_1} ( \oplus_{h\in C} \mathbb Z_2) \rtimes C \cong \oplus_{\mathbb N} (\oplus_{\mathbb Z} \mathbb Z_2 )\rtimes \mathbb Z$, where the last semidirect product is associated with diagonal action of Bernoulli actions of  $\mathbb Z$ on $\oplus_{\mathbb N} (\oplus_{\mathbb Z} \mathbb Z_2 )$. A similar argument shows that $\Upsilon =\mathbb Z_2 \wr_{K/B} D\cong \oplus_{\mathbb N} (\oplus_{\mathbb Z} \mathbb Z_2 )\rtimes \mathbb Z$ and combining with the above we get the claim. Next, fix a group isomorphism $\psi: \Omega \rightarrow \Upsilon$.

     \noindent Now, let $n\geq 2$ be any integer and consider the $n$-folded product $H= \G\times ...\times \G$ together with the n-folded product subgroup $\Sigma= \Omega\times ...\times \Omega$. Also denote by $\varphi: \Omega\times ...\times \Omega \rightarrow \Upsilon\times ...\times \Upsilon$ the $n$-folded isomorphism induced by $\psi$. Now one can check that $\Sigma< H$ and $\varphi(\Sigma)<H$ satisfy all the conditions enumerated in Lemma \ref{controlquasinorminhnn} and consequently the one-sided quasinormalizer conditions in the definition of class $\mathscr D$.   So using this construction in conjunction with amalgamation and HNN-extensions we can build iteratively various examples of groups in the class $\mathscr D$ such as:  
      \begin{equation*}\begin{split}&{\rm HNN}(H,\Sigma, \varphi)\in \mathscr D_1,\qquad {\rm HNN}(H,\Sigma, \varphi)\ast_\Sigma (H\ast_\Sigma H) \in \mathscr D_2, \\
      &({\rm HNN}(H,\Sigma, \varphi)\ast_\Sigma (H\ast_\Sigma H))\ast_\Sigma ({\rm HNN}(H,\Sigma, \varphi)\ast_\Sigma (H\ast_\Sigma H))\in \mathscr D_3, {\rm etc}. \\
      \end{split}\end{equation*}

\section{A class of semidirect product groups with non-amenable core}\label{semidirectwsuperrigid}

\noindent Our class $\mathscr A$ introduces a new family of semidirect product groups $G$ whose von Neumann algebras $\cL(G)$ display excellent rigidity properties, as we will see in the next sections. We continue by recalling the definition of class $\mathscr A$.

\noindent Let $\G$ be a non trivial, icc, bi-exact, torsion free, property (T) group. Let $n\ge 2$ be a positive integer and let  $\G_1$, $\G_2$, ..., $\G_n$ be isomorphic copies of $\G$. For every $1\leq i\leq n$ consider the action $\G \curvearrowright^{\rho^i} \G_i$ by conjugation, i.e.\ $\rho^i_\g (\la)=\g\la \g^{-1}$ for all $\g\in \G, \la\in \G_i$. Next consider the action $\G\ca^{\rho} \G_1\ast \G_2\ast ...\ast \G_n$ on the  free product group $\G_1 \ast \G_2\ast...\ast \G_n$ given by the canonical free product automorphism $\rho_\g = \rho^1_\g \ast \rho^2_\g \ast ...\ast \rho^n_\g$ for all $\g\in \G$ and let $G= (\G_1\ast \G_2 \ast... \ast \G_n) \rtimes_{\rho} \G$ be the corresponding semidirect product. 

\noindent The class of these semidirect product groups is denoted throughout the paper by {\bf Class $\mathscr A$}.

\noindent {\bf Representation as amalgams.} The groups in the class $\mathscr A$ can be viewed alternatively as free product groups amalgamated over the acting group. Namely, one can canonically decompose $G=(\G_1 \ast \G_2 \ast....\ast \G_n)\rtimes_{ \rho} \G= (\G_1 \rtimes_{\rho^1} \G)\ast_\G (\G_2 \rtimes_{\rho^2} \G)\ast_\G ... \ast_\G (\G_n \rtimes_{\rho^n} \G )$. In addition, the semidirect product $\G_i \rtimes_{\rho^i} \G$ can be canonically identified with the semidirect product $(\G \times 1) \rtimes_\rho d(\G)$ where $d(\G)= \{ (\g,\g) \,:\, \g\in \G\}\leqslant \G \times \G$ is the diagonal group and the action is given by $\rho_{(\g,\g)}(\la,1)= (\g \la \g^{-1},1)$ for all $\g,\la\in \G$. In particular, this canonically shows that $\G_i\rtimes_{\rho^i} \G \cong \G\times \G$. Thus, using the aforementioned identifications we have \begin{equation}\label{amalgamdecomp}G = \left ((\G_1\times 1 )\rtimes_{\rho^1} d (\G)\right )\ast_ {d(\G)}\left ((\G_2\times 1 )\rtimes_{\rho^2} d (\G)\right )\ast_{d(\G)}... \ast_{d(\G)} \left( (\G_n\times 1 )\rtimes_{\rho^n} d(\G)\right ).\end{equation}    

\noindent This amalgam decomposition of $G$ along the retracts will be used extensively in the proofs of our main structural results. 

\noindent We end this section by recording a list of algebraic properties of groups in class $\mathscr A$ 
that are relevant to our von Neumann algebraic results.

\begin{proposition}\label{P:classA} Let $G= (\G_1\ast \G_2 \ast... \ast \G_n) \rtimes_{\rho} \G\in \mathscr A$. Then the following hold:
\begin{enumerate}
    \item [1)] 
    $G$ has trivial amenable radical, i.e. the only normal amenable subgroup of G is the trivial one.
    \item [2)]If $\G$ is residually finite then so is $G$.
    \item [3)]The class $\mathscr A$ has $2^{\aleph_0}$ elements.
\end{enumerate}

\end{proposition}

\noindent {\it Proof.} 
1) Denote by $K= \G_1\ast ...\ast \G_n$ and note that $G=K\rtimes_\rho \G$. Fix $\Sigma\lhd G$ an amenable normal subgroup. First we argue that $\Sigma \cap K =1$. Since $\Sigma$ is normal in $G$, then $\Sigma\cap K$ is also normal in $K$. As $\Sigma\cap K<K=\G_1\ast\G_2\ast ...\ast\G_n$ then by Kurosh's subgroup theorem we have that $\Sigma\cap K = F \ast (\ast^k_{j=1} N_j^{\g_j})$, where $F$ is a free group, $N_j\leqslant \G_{i_j}$ $1\leq i_j \leq n$ are subgroups, and $\g_j \in K$. Here, we denoted $N_j^{\gamma_j}=\gamma_jN_j\gamma_j^{-1}$. As $\Sigma\cap K$ is amenable and torsion free we must have that either (a) $F=1$ and $k=1$, or (b) $F=\mathbb Z$ and $k=0$. By assuming (a), it follows that $\Sigma\cap K= (\Sigma\cap K)^{\g_1^{-1}}=N_1< \G_{i_1}$. Now, pick $\g\in K\setminus \Gamma_{i_1}$.  Since $\G_{i_1}$ is malnormal in $K$ and $\Sigma\cap K$ is normal in $K$, we get  $(\Sigma\cap K) =(\Sigma\cap K)  \cap (\Sigma\cap K)^{\g}< \G_{i_1}\cap \G^{\g}_{i_1}=1$, as claimed. Assume now that (b) holds and let $a=a_1\dots a_s$ be a generator of $F$ written as a reduced word in $\Gamma_1*\dots*\Gamma_n$. If $s=1$, then the argument follows as in (a). If $s=2$, note that any element in $F$ can be written as a reduced word of even length. Since $F<K$ is normal, it follows that the reduced word of $a_1aa_1^{-1}\in F$ has length $3$, contradiction. In the case that $s\ge 3,$   we proceed as follows.
Since $gag^{-1}$ is also a generator of $F$ for any $g\in K$, we derive that we either reduce to the case $s\in\{1,2\}$ or we can assume that $a_1$ and $a_s$ do not belong to the same subgroup $\Gamma_j$ of $K$. In this second case, note that any element of $F$ can be written as a reduced word of for which its length is a multiple of $s$. However, the length of the reduced word of $a_1 a a_1^{-1}\in F$  
equals $s+1$, contradition. Therefore,  $\Sigma\cap K=1$.

\noindent Next, notice that since $\Sigma$ and $K$ are normal in $G$ it follows that the commutator $[\Sigma , K] <\Sigma \cap K=1$; in particular, $\Sigma < C_G( K) $. Next, we argue that the centralizer $C_G(K)=1$, which in particular gives the desired conclusion. Fix $\g=kl\in C_G(K)$, where $k\in K$ and $l\in \G$. This implies that for all $s\in K$ we have $s\gamma=\gamma s$ which implies $s kl =kl s$  and hence $k^{-1}sk=\rho_l(s)$. If we let $s\in \G_i$ we see the previous relation together with the malnormality of $\G_i$ in $K$ imply that $k\in \G_i$. Since this holds for all $i$ then $k\in \cap^n_{i=1} \G_i=1$ and so $k=1$. In conclusion, we must have that $s=\rho_l(s)$ for all $s\in K$ and since $\G$ is icc this further implies that $l=1$; hence $\g=1$, which finishes the proof.

\noindent 2) Notice that since $\G$ is residually finite then so is $\G\times \G$ and hence, $\G_i \rtimes _{\rho^i}\G$ is residually finite for all $1\leq i\leq n$. Then using the amalgam decomposition of $G$ along retracts \eqref{amalgamdecomp} together with \cite[Theorem 1]{BE73}, iteratively, we get that $G$ is residually finite as well.

\noindent 3) We will present a construction of a continuum of elements in $\mathscr A$ that relies heavily on several deep results in geometric group theory \cite{AMO,Os06}. We start by noticing that for every finitely generated, torsion free group $K$ there exists a group $H(K)$ containing $K$ as proper subgroup and satisfying the properties that $H(K)$ is torsion free, has property (T) and is hyperbolic relative to $K$. This essentially follows from the same arguments presented in the proof of \cite[Theorem 1.1]{AMO}. However, since in the aforementioned result the authors do not emphasize the torsion free aspect we repeat here a simplified version of their argument addressing this part. To this end, let $T$ be any torsion free, property (T), hyperbolic group (e.g.\ any uniform lattice in ${\rm Sp}(n,1)$, $n\geq 2$) and let $F$ be a finite set of generators of $K$. Now consider the free product $G= T\ast K$ and notice $G$ is hyperbolic relative to $\{K\}$. In addition, notice that $T$ is a suitable subgroup of $G$ in the sense of \cite[Defintion 2.2]{Os06}. Then using \cite[Theorem 2.4]{Os06} one can find an epimorphism $\phi: G \rightarrow H$ satisfying the following properties: a) the restriction $\phi_{|K}$ is injective; b) the group $H$ is hyperbolic relative to $\phi(K)$; c) $\phi(F)\subset \phi(T)$; d) every element of finite order in $H$ is the image under $\phi$ of an element of finite order in $G$. Clearly, a) and b) imply that $\phi(K)\cong K$ and $H$ is hyperbolic relative $\phi(K)$. Since $T$ and $K$ are torsion free then so is $G$ and by d) it follows that $H$ is torsion free as well. Finally, condition c) implies that $H= \phi(T\ast K)= \phi(T)$ and since $T$ has property (T) then $H$ has property (T) as well. Letting, $H(K):=H$ we get the desired statement.

Next, we claim that there exists a continuum family $\tilde{\mathcal K}$ of pairwise non-isomorphic non-elementary amenable groups, that are torsion free and have infinite center.  Indeed, from \cite[Theorem 6]{Ha54} there exists a continuum family $\mathcal K$ of groups $K_i$ that are $2$-generated,  torsion free, and  solvable (in particular amenable). Using this we define a new continuum family of groups $\tilde {\mathcal K}$ as follows.  First eliminate all possible elementary groups from $\mathcal K$ which are at most countably many so we are left again with a continuum family which we still denote by $\mathcal K$. Consider $\mathcal K_c \subseteq \mathcal K$ the subset of the groups in $\mathcal K$ with infinite center. If $|\mathcal K_c|=2^{\aleph_0}$ then let $\tilde {\mathcal K}:=\mathcal K_c$. If $|\mathcal K_c|\neq 2^{\aleph_0}$ then  $|\mathcal K\setminus \mathcal K_c|=2^{\aleph_0}$. Moreover, since all groups involved are torsion free then $\mathcal K\setminus \mathcal K_c$ consists only of groups with trivial center. Then in this scenario we define $\tilde{\mathcal K}:= \{\mathbb Z \times K\,:\, K\in \mathcal K\setminus \mathcal K_c\}$. This proves the claim. 

Now, we argue that the groups $H(S_i)$ where $S_i \in \tilde {\mathcal K}$ form a continuum family of pairwise non-isomorphic, icc, bi-exact property (T), torsion free groups. To conclude this we only need to show the non-isomorphism part as the rest follows from the prior paragraph. Assume $\theta:H(S_i)\rar H(S_j)$ be a group isomorphism. Fix an infinite order central element $a\in Z(S_i)$. Thus, $\theta(a)\in H(S_j)$ is an infinite order element as well. Assume $\theta(a)$ is a hyperbolic element of $H(S_j)=B$. Thus, by \cite[Theorem 2.1]{Os06} there exists an elementary group $E_B(\theta(a))$ such that $B$ is hyperbolic relative to $\{S_j\}\cup \{E_B(\theta(a))\}$. In particular, $E_B(\theta(a))$ is malnormal in $B$. As $\langle \theta(a)\rangle$ commutes with $\theta (S_i)$ it follows that $\theta(S_i)<E_B(\theta(a))$ which further entails that $\theta(S_i)$ and hence $S_i$ is elementary, a contradiction. In conclusion, $\theta(a)$ is parabolic and hence there exists  $h\in B$ such that $\langle \theta(a)\rangle^h\subseteq S_j$. Again since  $S_j<B$ is malnormal and $\theta(S_i)^h$ commutes with $\langle\theta(a)\rangle^h$ it follows that  $\theta(S_i)^h< S_j$. Using a similar argument for $\theta^{-1}$ one can find $k\in B$ such that
$\theta(S_i)^k< S_j$ and by malnormality again there is $s\in B$ such that $\theta(S_i)^s=S_j$; in particular, $S_i\cong S_j$ and hence $i=j$ which finishes the  argument. 

Finally, it is a basic exercise to see that if one starts with $K=H(S_i)$, $S_i\in \tilde{\mathcal K}$ in the semidirect product construction in the class $\mathscr A$ one gets non-isomorphic groups for different $i$'s. We leave the details to the reader. 
\hfill$\blacksquare$

\section{Height of elements in group von Neumann algebras and techniques for discretization of countable groups}
\noindent The notion of height of elements in crossed products and group von Neumann algebras was introduced and developed in \cite{Io10} and \cite{IPV10} and was highly instrumental in many of the recent classification results in von Neumann algebras \cite{Io10,IPV10,KV15,CI17,CU18,CDK19,CDHK20}. Following \cite[Section 3]{IPV10} for every $x\in \cL(\G)$ we denote by $h_{\G}(x)$ the largest Fourier coefficient of $x$, i.e.,
$h_{\G}(x)={\max}_{\g\in \G }\ |\tau(xu_\g^*)|$.  Moreover,  for every subset $\cG\subseteq \El(\G)$, we denote by
$h_{\G}(\cG)=\inf_{x\in \cG}h_{\G}(x)$, the height of $\cG$ with respect to $\G$.   Using the notion of height Ioana, Popa and Vaes proved in their seminal work, \cite[Theorem 3.1]{IPV10} that whenever $\G$, $\La$ are icc groups such that $\El(\G) = \El(\La)$ and $h_{\G}(\La)>0$, then $\G$ and $\La$ are isomorphic. Therefore, in order to reconstruct the underlying groups from their von Neumann algebras a first step is to develop an adequate analysis to control the lower bound of their height.  

\noindent There have been a few situations in the literature where it was possible to obtain lower bounds for the height.
At the heart of these results is the following common philosophy that was extensively exploited: given two group von Neumann decompositions of  $\cM =\cL(\G)=\cL(\La)$, to conclude that the height $h_\G(\La_0)>0$ for some subgroup $\La_0<\La$ sometimes it suffices to check that there are only a few subgroups $\G_i\in {\rm Sub}(\G)$ and $\La_i\in {\rm Sub}(\La)$ such that their von Neumann algebras can be identified  $\cL(\G_i)=\cL(\La_i)$ in $\cM$ or just merely  intertwined into each other. For example, this is the case of certain wreath products $\G=A^{(H)}\rtimes  H$ in \cite{IPV10} and left-right wreath products $\G= A^{(H)}\rtimes( H\times H)$ in \cite{BV12} where the von Neumann algebras of the core groups $A^{(H)}$ and of the acting groups $H$ and respectively $H\times H$ could be identified with the von Neumann algebras of certain subgroups in the mystery group $\Lambda$.  A similar statement was proved for semidirect products with no non-trivial stabilizers in \cite{CDK19,CDHK20}.

 \noindent Next, we highlight a rather different situation  where one can control the lower bound for height of unitary elements in the context of direct product groups. This is reminiscent to some of the techniques from \cite{CI17}. To properly state our result we first introduce the following definition:

\begin{definition}\label{D:special}   
\begin{enumerate}
    \item Let $n\geq 3$ be a positive integer and let $I\sqcup J =\{1,2...,n\}$ be a partition with $|I|\geq 2$. Let $\Sigma, \G_1, ... ,\G_n\in {\rm Sub} (\G)$ be a collection of subgroups and consider the (ordered) $n$-tuple of subgroups $\mathcal F=(\G_1, \G_2,..., \G_n)\in {\rm Sub}(\G)^n$. We say that $\Sigma$ is $I$-$J$-{\it fixable with respect to} $\mathcal F$ if the following property holds: for any finite subsets $F_i,K_i\subset \G$  where $1\leq i\leq n$ there exist  a finite set $G_j \subset F_j \G_j K_j$ when $j\in J$ and $l_i$ injective maps $\sigma^k_i: \Sigma\setminus$  $ \{1\} \rightarrow \G$ for $1\leq k\leq l_i$ and $i\in I$ such that whenever $g \in \Sigma \setminus \{1\}$ and $g_i\in F_i\G_i K_i \setminus \{1\}$ for $1\leq i\leq n$ are elements satisfying $g_j\in F_j\G_jK_j\setminus G_j$ for $j\in J$ and   $g g_1 g_2 ...  g_n =1$, then for every $i\in I$ we must have  $g_i = \sigma^{k}_i (g)$ for some $1\leq k\leq l_i$.
    
    \item If $J=\emptyset$ in (1), we simply say that $\Sigma$ is {\it fixable with respect to} $\mathcal F$.
    
    \item If we have a tuple of subgroups  $\mathcal G=(\G_1,\G_2,...,\G_m)\in {\rm Sub}(\G)^m$ and $\G_i$ is fixable with respect to $\hat {\mathcal G}_i=(\G_1,\G_2,...,\G_{i-1},\G_{i+1},..., \G_m)$ for all $1\leq i\leq m$, then we say that $\mathcal G$ is {\it fixable}.
\end{enumerate}

\end{definition}
\noindent While this definition seems somewhat technical there are in fact many natural examples of groups $\Sigma<\Gamma$ such that $\Sigma$ is fixable with respect to certain families of subgroups of $\Gamma$. This includes, for instance, the collections of the so-called ``diagonal subgroups''. More precisely, we have the following result for which the proof we leave it to the reader.

\begin{proposition}\label{thinexamples}
\begin{enumerate}
    \item  Let $n\ge2$ be a positive integer and let $\Sigma, \G_1, \G_2, ... ,\G_n$ be some groups. Assume that $\pi_i : \Sigma \rightarrow \G_i$ is a monomorphism for all $1\leq i\leq n$ and consider the diagonal subgroup $\delta(\Sigma)= \{(\pi_1(g), \pi_2(g),..., \pi_n(g))\in \G_1\times \G_2 \times \cdots \times \G_n \,:\, \g\in \Sigma\} \leqslant \G_1\times \G_2\times \cdots \times\G_n=\G$. Then the $(n+1)$-tuple  $\mathcal F=( \delta(\Sigma),\G_1, \G_2,..., \G_n) \in {\rm Sub}(\G)^{n+1}$ is fixable.  
    \item  Let $\G =A\rtimes_\rho G$ be a semidirect product and let $H\leqslant G$ be a subgroup. Assume that there exists a map $c: H \rightarrow A\setminus \{1\}$ such that $c_{gh} = c_g \rho_g(c_h)$ for all $g,h\in H$. If we denote by $\delta(H)=\{ c_h h \,:\,h\in H\}$, then $\delta(H)$ is fixable with respect to $\{A,G\}$.

    \item  Let $\G= A \rtimes_\rho G$ be a semidirect product group. Then $G$ is $\{1,3\}$-$\{2\}$-fixable with respect to $\mathcal F=(A,G,A)\in {\rm Sub}(\G)^3$. 
\end{enumerate} \end{proposition}

\noindent With {these} preparations at hand we are now ready to derive the first main result of the section. Specifically, we show that in the presence of groups that are $I$-$J$-fixable, it is possible to control the lower bound for the heights of elements that satisfy various relations in the von Neumann algebra setting.

\begin{lemma}\label{inclexcl} Let $\Sigma \leqslant\Gamma$ be an inclusion of groups and let $F,K\subset \G$ be subsets. Then for every $x\in L(\Gamma)$ we have 
\begin{equation}\label{inftynormcontrol}
    \|P_{F \Sigma K}(x)\|_\infty \leq (2^{|F||K|}-1)\|x\|_\infty.
\end{equation}
    \end{lemma}
\begin{proof} We only need to show \eqref{inftynormcontrol} when $F$ and $K$ are finite, the other cases being tautological. Towards this, observe that for every $s,t\in \Gamma$ and $x\in L(\Gamma)$ we have that $P_{s\Sigma t}(x)= u_s E_{L(\Sigma)}(u_{s^{-1}} x u_{t^{-1}})u_t$. In particular, we have $\|P_{s \Sigma t}(x)\|_\infty \leq \|x\|_\infty$. This already proves our statement when $F$ and $K$ are singletons. The general case follows from this combined with the inclusion-exclusion principle for orthogonal projections. To see this, consider the sets $\Sigma_{s,t}= s\Sigma t$ for $s\in F$, $t\in K$ and enumerate them as $\{S_i\}_i$, for $1\leq i\leq k$ where $k = |F||K|$. Thus, one can check that 
\begin{equation}\label{incl-excl}P_{F\Sigma K}= P_{\cup^k_{i=1} S_{i}}= \sum^k_{i=1} (-1)^{i+1}  \sum_{1\leq j_1<\cdots<j_i\leq k} P_{ S_{j_1}\cap \cdots \cap S_{j_i}}\end{equation}
Next, we notice that if $\Sigma'<\Sigma$ is a subgroup and $s,t\in\Gamma$, then $\Sigma'\cap s\Sigma t$ is either trivial or of the form $a\Sigma''$ for some subgroup $\Sigma''<\Sigma'$ and $a\in\Gamma$. Indeed, if $\Sigma'\cap s\Sigma t$ is not trivial, then there exist $\sigma'\in\Sigma'$ and $\sigma\in\Sigma$ such that $\sigma'=s\sigma t$. Hence, $\Sigma'\cap s\Sigma t=\Sigma'\cap \sigma' t^{-1}\Sigma t=\sigma'(\Sigma'\cap t^{-1}\Sigma t)$.

Now, the previous paragraph implies that every $S_{j_1}\cap \cdots \cap S_{j_i}$ is either empty or of the form $g\Sigma' h$ for some subgroup $\Sigma'\leqslant \Gamma$ and $g,h\in \Gamma$. 
Hence, the first part of the proof implies that $\|P_{S_{j_1}\cap \cdots \cap S_{j_i}}(x)\|_\infty \leq \|x\|_\infty$, for all $1\leq i\leq k$ and $x\in L(\Gamma)$. This together with \eqref{incl-excl} and the triangle inequality imply that $\|P_{F\Sigma K}(x)\|_\infty \leq \sum^k_{i=1}\sum_{1\leq j_1<\cdots<j_i\leq k} \| P_{S_{j_1}\cap \cdots \cap S_{j_i}}(x)\|_\infty \leq (2^k -1)\|x\|_\infty$, as desired. \end{proof}

\begin{theorem}\label{thingroupheight} Let $\Sigma, \G_1, ... ,\G_n\leqslant \G$ with 
$n\geq 2$ and let $I\sqcup J=\{1,...,n\}$ be a partition.  Assume that $\Sigma$ is $I$-$J$-fixable with respect to $\mathcal F=(\G_1,\ldots, \G_n)$. Also let $\cM_1, ..., \cM_n\subseteq \cL(\G)=\cM$ be von Neumann subalgebras such that $\cM_i\prec^s_\cM \cL(\G_i)$ for all $i\geq 1$. Also for every $j\in J$ assume $(x^j_k)_k\subseteq (\cM_j)_1$ is a sequence so that $x^j_k \rightarrow 0$ in the WOT topology, as $k\rightarrow \infty$.  Let $\mathscr G\subseteq \mathcal U(\cL(\Sigma))$ such that for every $x\in \mathscr G$ there exist some elements $x_i\in \mathscr U( \cM_i)$ for $i\in I$ such that for all $k\in \mathbb N$ we have $x a^1_ka^2_k...a^n_k=1$, where $a^i_k =x_i$  if $i\in I$ and $a^j_k= x^j_k$ if $j\in J$. Then the height $h_{\Sigma}(\mathscr G)>0$.
\end{theorem}

\noindent {\it Proof.} By hypothesis we have $\cM_t\prec^s_\cM \cL(\G_t)$ for all $t\geq 1$. Fix $0<\varepsilon<\frac{1}{3}$. Using \cite{Va10} recursively, for every $1\leq t\leq n$ one can find finite subsets $F_t, K_t \subset \G$ so that for all $y\in (\cM_t)_1$ we have \begin{equation}\label{ineqproj1}
    \|y-P_{F_t \G_t K_t}(y)\|_2\leq \frac{\varepsilon}{n \prod^{t-1}_{m=1} 2^{|F_m||K_m|}}. 
\end{equation}
Here, and throughout the rest of the proof we make the convention that $\prod^{0}_{m=1} 2^{|F_m||K_m|}=1$.

\noindent Fix $x\in \mathscr G$. By hypothesis we have that $x a^1_ka^2_k...a^n_k=1$, where $a^i_k =x_i$  if $i\in I$ and $a^j_k= x^j_k$ if $j\in J$. Also for simplicity of the writing denote by $S_i= F_i \G_i K_i$ for $1\leq i\leq n$. From Lemma \ref{inclexcl} we have that $\|P_{S_i}(x)\|_\infty \leq  2^{|F_i||K_i|}\|x\|_\infty$ for all $x\in \mathcal M$. Using this in combination with the triangle inequality and inequalities \eqref{ineqproj1} we see that 
\begin{equation}\label{norm2}\begin{split}
 \|1- x\prod^n_{t=1} P_{S_t}(a^t_k) \|_2&=  \|xa^1_k a^2_k...a^n_k- x\prod^n_{t=1} P_{S_t}(a^t_k) \|_2 \\& 
 \leq \sum^n_{t=1} \|(\prod^{t-1}_{m=1} P_{S_m}(x^m_k))(a^t_k - P_{S_t} (a^t_k)) (\prod^n_{m=t+1} a^m_k)\|_2\\
 &\leq \sum^n_{t=1} \prod^{t-1}_{m=1} \|P_{S_m}(x_k^m)\|_\infty \|a^t_m - P_{S_t} (a^t_k)\|_2\\& 
 \leq \sum^n_{t=1} \prod^{t-1}_{m=1} 2^{|F_m||K_m|}  \frac{\varepsilon}{n \prod^{t-1}_{m=1} 2^{|F_m||K_m|}}
 \\&
 = \sum^n_{t=1}\frac{\varepsilon}{n}=\varepsilon.
 \end{split}
\end{equation}
Combining the previous inequality with $|1- \tau (x\prod^n_{t=1} P_{S_t}(a^t_k))|\leq \|1- x\prod^n_{t=1} P_{S_t}(a^t_k) \|_2$  and the triangle inequality we further see that for all $k$ we have  \begin{equation}\label{ineqproj4}1-\varepsilon \leq |\tau (x\prod^n_{t=1} P_{S_t}(a^t_k))|.\end{equation}

\noindent For every $j\in J$ pick a finite subset $G_j\subset F_j\G_jK_j=S_j$ satisfying the condition in Definition \ref{D:special}. Let $J= \{j_1,j_2,...,j_r\}$ for some $r$. Using that $x^j_k \rightarrow 0$ in WOT as $k\rightarrow \infty$ for every $j\in J$, we can choose an $k$ such that for every $1\leq s\leq r$ there is a finite set $R_{j_s} \subset F_{j_s}\G_{j_s}K_{j_s} \setminus G_{j_s}$ satisfying 

\begin{equation}\label{ineqproj1more}
    \|x^{j_s}_k-P_{R_{j_s}}(x^{j_s}_k)\|_2\leq \frac{2\varepsilon}{n (\prod_{i\in I} 2^{|F_i||K_i|})( \prod^{s-1}_{v=1} |R_{j_v}|)( \prod^r_{v=s+1} |2^{|F_{j_v}||K_{j_v}|})}. 
\end{equation}
Again, in this formulas we  convene that $\prod^{0}_{v=1} |R_{j_v}| =\prod^r_{v=r+1} 2^{|F_{j_v}||K_{j_v}|}=1$.

\noindent Exploiting \eqref{ineqproj1more} and \eqref{ineqproj4}, in the same way as in \eqref{norm2} we get that $1-3\varepsilon \leq |\tau (x\prod^n_{t=1} P_{W_t}(a^t_k))|$ where we have denoted by $W_i= S_i$ if $i\in I$ and $W_j= R_j$ if $j\in J$. This further implies that  

 
\begin{equation}\label{ineqheight4}\begin{split}
    1-3\varepsilon &\leq |\tau (x\prod^n_{t=1} P_{W_t}(a^t_k))|\\&= |\sum_{\substack{g\in \Sigma\\  g_t\in W_t,1\leq t\leq n\\ gg_1g_2...g_n=1}} \tau(x u_{g^{-1}}) \prod^n_{t=1} \tau(a^t_k u_{g^{-1}_t})| \\&
    \leq \sum_{\substack{g\in \Sigma\\  g_t\in W_t,1\leq t\leq n\\ gg_1g_2...g_n=1}} |\tau(x u_{g^{-1}})| \prod^n_{t=1} |\tau(a^t_k u_{g^{-1}_t})|.
    \end{split}
\end{equation}
    
\noindent Next, since $\Sigma$ is $I$-$J$-fixable with respect to $\mathcal F=(\G_1,..., \G_n)$ and $g g_1 g_2...g_n=1$ then for every $t\in I$ there are injections $\sigma^s_t: \Sigma\setminus \{1\} \rightarrow G$ for $1\leq s\leq l_t$ such that   $g_t = \sigma^{s_t}_t (g)$ for some $1\leq s_t\leq l_t$. Choose $\sigma_t^s(1)\in G$ such that $\sigma_t^s:\Sigma\to G$ is still injective.
Let $I=\{i_1,i_2,\dots,i_p\}$ for some $p\ge 2$. This together with Cauchy-Schwarz inequality show the last term in \eqref{ineqheight4} is smaller than  
    
    \begin{equation*}\begin{split} &\sum_{\substack{g\in \Sigma\\  g_t\in W_t,1\leq t\leq n\\ gg_1g_2...g_n=1}} |\tau(xu_{g^{-1}})|\prod^n_{t=1} |\tau(a^t_k u_{g_t^{-1}})|
    \\& \leq  h_{\Sigma}(x)\sum_{\substack{g\in \Sigma\\  g_t\in W_t,1\leq t\leq n\\ gg_1g_2...g_n=1}}  \prod_{t\in I} |\tau(x_t u_{g_t^{-1}})|\prod_{t\in J} |\tau(x^t_k u_{g_t^{-1}})|
    \\&\leq  h_{\Sigma}(x)\sum_{\substack{g\in \Sigma \\ g_t\in W_t, t\in J}}  \prod_{t\in I} |\tau(x_t u_{\sigma^{s_t}_t(g)^{-1}})|\prod_{t\in J} |\tau(x^t_k u_{g_t^{-1}})|
    \\&\leq  h_{\Sigma}(x)\sum_{\substack{g\in \Sigma \\ g_t\in W_t, t\in J}}  \prod_{t\in I} |\tau(x_t u_{\sigma^{s_t}_t(g)^{-1}})|
    \\& \leq  h_{\Sigma}(x)  (\prod_{t\in J}|W_t|) 
    \sum_{g\in \Sigma}  \prod_{t\in I} (\sum_{ 1\leq s\leq l_t} |\tau(x_t u_{\sigma^{s}_t(g)^{-1}})|) 
    \\& \leq h_{\Sigma}(x)  (\prod_{t\in J}|W_t|) 
     (\prod_{t\in I\setminus \{i_1,i_2\}} l_t)  (\sum_{\substack{g\in \Sigma \\ 1\leq s\leq l_{i_1} \\ 1\leq r\leq l_{i_2}}} |\tau(x_{i_1} u_{\sigma^{s}_{i_1}(g)^{-1}})| |\tau(x_{i_2} u_{\sigma^{r}_{i_2}(g)^{-1}})|)
    \\&
    \leq h_{\Sigma}(x)  (\prod_{t\in J}|W_t|) 
     (\prod_{t\in I\setminus \{i_1,i_2\}} l_t)  (\sum_{\substack{g\in \Sigma \\ 1\leq s\leq l_{i_1} \\ 1\leq r\leq l_{i_2}}} |\tau(x_{i_1} u_{\sigma^{s}_{i_1}(g)^{-1}})|^2)^{1/2}  (\sum_{\substack{g\in \Sigma \\ 1\leq s\leq l_{i_1} \\ 1\leq r\leq l_{i_2}}}|\tau(x_{i_2} u_{\sigma^{r}_{i_2}(g)^{-1}})|^2)^{1/2}
    \\&
    \leq h_{\Sigma}(x)  (\prod_{t\in J}|W_t|) 
     (\prod_{t\in I\setminus \{i_1,i_2\}} l_t)  ( l_{i_1} l_{i_2})^{1/2} \|x_{i_1}\|_2  ( l_{i_1} l_{i_2})^{1/2} \|x_{i_2}\|_2\\&=h_{\Sigma}(x)  (\prod_{t\in J}|W_t|) 
     (\prod_{t\in I} l_t). 
     \end{split}
\end{equation*}
Altogether, these imply that $h_\Sigma(x)\geq \frac{1-3\varepsilon}{(\prod_{t\in J}|W_t|) 
     (\prod_{t\in I} l_t)}$ for all $x\in \mathcal F$, as desired. 
\hfill$\blacksquare$

\noindent The next corollary will be particularly useful in the proofs of the main results. 
\begin{corollary}\label{identificationdiagonalgroups} Let $\G$ be an icc nonamenable bi-exact group and denote by $\cM= \cL(\G\times \G)$. Let  $\La$ be an arbitrary group together with a subgroup $\Omega\leqslant \La$ such that  $\cM= \mathcal L(\La)$ and $\mathcal L(d(\G))=\cL(\Omega)$. 

\noindent Then one can find a unitary $w\in \cM$ such that $\mathbb T w(\G \times \G)w^*=\mathbb T \La$.    

\end{corollary} 

\noindent {\it Proof.} Denote by $\G_1 =\G\times 1$ and $\G_2 = 1\times \G$. Next, we claim that the conditions of \cite[Theorem 5.1]{CI17} are satisfied for $\Sigma=d(\G)$. Towards this, let $\rho_i: \Gamma_1\times \Gamma_2 \rightarrow \Gamma_i$ be the canonical group projection $\rho_i(g_1,g_2)=g_i$, for $i=1,2$. Now notice that  the restrictions $\rho_i : \Sigma\rightarrow \Gamma_i$ are injective. Fix $h_i \in \Gamma_i\setminus \{1\}$. From definitions one can see that $\{ \rho_i (g,g) h_i \rho_i (g,g)^{-1}\,:\, (g,g)\in \Sigma \}=\{ g h_i g^{-1} \,:\, g\in \Gamma_i\} $. As $\G_i$ is icc, it follows that the previous set is infinite which yields our claim.

Thus, using the conclusion of \cite[Theorem 5.1]{CI17} one can find a unitary $u\in \cM$ and a product decomposition $\La=\La_1\times \La_2$ such that $u\cL(\G_1)u^*=\cL(\La_1)$,  $u\cL(\G_2)u^*=\cL(\La_2)$ and $\mathbb T u\Sigma u^*=\mathbb T \Omega $. Moreover, the
second paragraph of the proof of \cite[Theorem 5.1]{CI17} shows that the restriction of the projection $\Lambda\to\Lambda_i$ to $\Omega$ is a monomorphism for any $i\in\{1,2\}$. Note also that we can identify $\Omega$ by $\{ (\pi_1(\g), \pi_2(\g))\,:\, \g\in \Omega\}$. 
Thus, denoting by $G_i:= u\G_i u^*$ and $H:= u\Sigma u^*$ we have \begin{equation}\label{equalityalg}
\cL(G_1)=\cL(\La_1),\qquad \cL(G_2)=\cL(\La_2), \text{ and} \quad \mathbb T H =\mathbb T \Omega.
\end{equation} 

\noindent Notice that by part 1) in Proposition \ref{thinexamples} the triple $(\La_1,\La_2,\Upsilon )\in {\rm Sub}(\La)^3$ is fixable. Also for every $g\times 1 \in G_1$ we have that $u_{g\times 1} u_{g^{-1}\times g^{-1}} u_{1\times g}=1$ where $u_{g\times 1} \in \cL(G_1)$, $u_{g^{-1}\times g^{-1}}\in \cL (H)$ and $u_{1\times g}\in \cL(G_2)$. Thus using relations \eqref{equalityalg} and Theorem \ref{thingroupheight} we have that $h_{\La_1}(G_1)>0$. Therefore by \cite[Theorem 3.1]{IPV10} there exists a unitary $u_1\in \cL(\La_1)$ such that $\mathbb T u_1 G_1 u_1^*= \mathbb T \La_1$. By a similar argument there is a unitary $u_2\in \cL(\La_2)$ such that $\mathbb T u_2 G_2 u_2^*= \mathbb T \La_2$. Altogether, these relations imply that $\mathbb T (u_1\otimes u_2 ) u (\G\times\G) u^*(u_1^*\otimes u_2^*) = \mathbb T (u_1\otimes u_2 ) (G_1\times G_2 )(u^*_1\otimes u^*_2 )= \mathbb T \La$ and letting $w= (u_1\otimes u_2)u$ we get the desired conclusion. \hfill$\blacksquare$

\noindent We end this section with two technical results regarding discretization of underlying groups in the von Neumann algebra regime. These will be used in an essential way to derive the main results of the next sections. 

\noindent The first result asserts that the discretization of two subgroups with infinite and ``sufficiently malnormal'' intersection can be bumped up to the group they generate.
\begin{theorem}\label{discretizationgeneratinggroups} Let $\G_1,\G_2\leqslant \G$ be  groups. Assume that  the subgroup $\Sigma=\G_1\cap \G_2\leqslant \G$ is icc and satisfies ${\rm QN}_\G^{(1)}(\Sigma)=\Sigma$. Let $\La$ be an arbitrary group such that $\cN=\cL (\G)=\cL(\La)$. Assume there exist $w_1,w_2\in \mathcal U(\cN)$ and subgroups $\La_1,\La_2\leqslant \La$ such that $\mathbb T w_1 \G_1 w_1^*= \mathbb T \La_1$ and $\mathbb T w_2 \G_2 w_2^*= \mathbb T \La_2$. 

\noindent Then one can find a unitary $w\in\mathcal U(\cN)$ such that  $\mathbb T w (\G_1\vee \G_2) w^*= \mathbb T (\La_1\vee \La_2)$.
\end{theorem}

\noindent {\it Proof.}
From assumptions there are group isomorphisms $\delta_i : \G_i\rightarrow \La_i$ and characters $\eta_i: \G_i \rightarrow \mathbb T$ so that 
\begin{equation}\label{eqsisom1}
w_i u_{\g_i} w_i^*= \eta_i (\g_i) v_{\delta_i(\g_i)} \text{ for all } \g_i\in \G_i \text{ and } 1\leq i\leq 2.
\end{equation}
These relations show that for all $\g\in \Sigma=\G_1\cap \G_2$ we have  $\eta_1(\g) w_1^* v_{\delta_1(\g)}w_1=u_\g= \eta_2(\g) w_2^* v_{\delta_2(\g)}w_2$. Thus, if we let $d_\g = \eta_1 (\g)^{-1}\eta_2(\g)$, we see that \begin{equation}v_{\delta_1(\g)}= d_\g w_1w_2^* v_{\delta_2(\g)} w_2w_1^*\text{ for all }\g\in \Sigma.\end{equation}
In particular, this relation entails that $\cL(\delta_1(\Sigma))
\prec_{\cN} \cL(\delta_2(\Sigma))$. By \cite[Lemma 2.6]{CI17} one can find $\la\in \La$ such that $[\delta_1(\Sigma): \la \delta_2(\Sigma)\la^{-1}\cap \delta_1(\Sigma) ]<\infty$. Therefore, replacing $w_2$ by $v_\la w_2$ and $\delta_2$ by ${\rm ad} (\la)\circ \delta_2$, we can assume that $[\delta_1(\Sigma):  \delta_2(\Sigma)\cap \delta_1(\Sigma) ]<\infty$ and relations \eqref{eqsisom1} still hold.

\noindent Also, \eqref{eqsisom1} show that  $\eta_1(\delta_1^{-1}(\la))^{-1} w_1 u_{\delta_1^{-1}(\la)} w_1^*= v_\la = \eta_2(\delta_2^{-1}(\la))^{-1} w_2 u_{\delta_2^{-1}(\la)} w_2^*$  for every $\la \in \delta_1(\Sigma)\cap \delta_2(\Sigma)$. Letting $e_\la= \eta_1(\delta_1^{-1}(\la))^{-1}\eta_2(\delta_2^{-1}(\la))\in \mathbb T$ and $w=w_2^*w_1$ this further shows that  
\begin{equation}\label{eqsisom3}e_\la w u_{\delta_1^{-1}(\la)} =  u_{\delta_2^{-1}(\la)} w \text{ for all } \la \in \delta_1(\Sigma)\cap \delta_2(\Sigma).
\end{equation}
\noindent Since $\delta_1^{-1}(\delta_1 (\Sigma) \cap \delta_2(\Sigma))\leqslant \Sigma$ has finite index there are $h_1,...,h_n\in \Sigma$ such that $\Sigma = \bigcup^n_{i=1} \delta_1^{-1}(\delta_1 (\Sigma) \cap \delta_2(\Sigma)) h_i$. Using this in combination with \eqref{eqsisom3} and  $\delta_2^{-1}(\delta_1 (\Sigma) \cap \delta_2(\Sigma))\leqslant \Sigma$ we get that \begin{equation*}\begin{split}
   w \cL(\Sigma) &\subseteq \sum^n_{i=1} w\cL(\delta_1^{-1}(\delta_1 (\Sigma) \cap \delta_2(\Sigma)) u_{h_i} = \sum^n_{i=1} \cL(\delta_2^{-1}(\delta_1 (\Sigma) \cap \delta_2(\Sigma)) w u_{h_i}\subseteq \sum^n_{i=1} \cL(\Sigma) w u_{h_i}.
   \end{split}
\end{equation*}
\noindent In particular, this shows that $w\in \mathscr {QN}^{(1)}_{\cL(\G)}(\cL(\Sigma))'' =\cL({\rm QN}^{(1)}_\G(\Sigma))=\cL(\Sigma)$. Consider the Fourier decomposition $w=\sum_\g a_\g u_\g$ in $\cL(\Sigma)$. Then using this in combination with relation \eqref{eqsisom3}  it follows that  $a_\g= e_\la a_{\delta_2^{-1}(\la) \g \delta_1^{-1}(\la^{-1})} $ for all $\la \in \delta_1 (\Sigma)\cap \delta_2 (\Sigma), \g\in \Sigma$. As $|e_\la|=1$ this implies that $a_\g\neq 0$ only if the set $O_\g=\{\delta_2^{-1} (\la)\g \delta_1^{-1}(\la)^{-1} \,:\, \la\in \delta_1(\Sigma)\cap \delta_2(\Sigma) \}$ is finite. This implies that there exists a finite index subgroup  $\Sigma_\g\leqslant \delta_1(\Sigma)\cap \delta_2(\Sigma) $ such that  $\delta_2^{-1} (\la)\g \delta_1^{-1}(\la^{-1})=\g$ and hence  \begin{equation}\label{partinnisom1}\delta_2^{-1} (\la)=\g \delta_1^{-1}(\la)\g ^{-1}\text{ for all }\la\in \Sigma_\g.\end{equation}

\noindent Next, let $\g,\mu \in \Sigma$ such that $O_\g$ and $O_\mu$ are finite.  By \eqref{partinnisom1}, for every $\la \in \Sigma_\g \cap \Sigma_\mu$ we get   $\mu \delta_1^{-1}(\la) \mu ^{-1}= \delta_2^{-1} (\la)=\g \delta_1^{-1}(\la) \g ^{-1}$ and hence $\mu^{-1}\g \in C_\G( \delta_1^{-1}(\Sigma_\la \cap \Sigma_\mu) )$; in particular,  the unitary  $u_{\mu^{-1}\g}\in \cL(\delta_1^{-1}(\Sigma_\la\cap \Sigma_\mu))'\cap \cL(\Sigma)$. Now, notice that since $[\Sigma : \delta_1^{-1}(\delta_1(\Sigma)\cap \delta_2(\Sigma))]<\infty, [\delta_1(\Sigma)\cap \delta_2(\Sigma) :\Sigma_\la], [\delta_1(\Sigma)\cap \delta_2(\Sigma):\Sigma_\mu]<\infty$, we deduce that $[\Sigma: \delta_1^{-1}(\Sigma_\la\cap\Sigma_\mu)]<\infty$. Therefore,  $\cL(\delta_1^{-1}(\Sigma_\la\cap \Sigma_\mu))'\cap \cL(\Sigma)\subseteq \cL(vC_\Sigma(\Sigma))$. Since $\Sigma$ is icc, we have that $vC_\Sigma(\Sigma)=1$, and thus, we conclude that $\cL(\delta_1^{-1}(\Sigma_\la\cap \Sigma_\mu))'\cap \cL(\Sigma)=\mathbb C1$. In particular, this implies that $u_{\mu^{-1}\g} =1$ and hence $\g=\mu$. Altogether, this shows that $w= c u_{\sigma}$ for some $\sigma\in \Sigma$  and $c\in \mathbb T$. Thus relations \eqref{eqsisom1} become $w_2 u_{\sigma \g_1 \sigma^{-1}} w_2^*= \eta_1(\g_1) v_{\delta_1(\g_1)}$ for all $\g_1\in \G_1$  and $w_2 u_{\g_2} w_2^*= \eta_2(\g_2) v_{\delta_2(\g_2)}$ for all $\g_2\in \G_2$. These relations clearly imply that  $\mathbb T w_2 ( \G_1 \vee \G_2) w_2^*= \mathbb T( \La_1\vee \La_2 )$, as desired.\hfill$\blacksquare$

\noindent The second result asserts that the elements that conjugate discretized subgroups can be discretized.

\begin{theorem}\label{discretizationhnngroups} Let  $\Theta, \Omega \leqslant \G$ be  groups and let $\Sigma\leqslant \Theta$ be a subgroup so that  $vC_\G(\Sigma)=1$. Let $\phi: \Sigma \rightarrow \Omega$ be a group homomorphism   and $t\in \G$ such that $\phi(\sigma ) = t \sigma t^{-1}$ for all $\sigma\in \Sigma$.   Let $\La$ be an arbitrary group such that $\cN=\cL (\G)=\cL(\La)$. Also assume there exist subgroups $\Phi, \Upsilon \leqslant \La$ and unitaries $x,y\in \mathcal U(\cN)$ such that $\mathbb T x \Theta x^*= \mathbb T \Phi$ and $\mathbb T y \Omega y^*= \mathbb T \Upsilon$. 

\noindent Then one can find $\la \in \La$ such that $y u_t x^*\in \mathbb T v_\la$.
    
\end{theorem}

\noindent {\it Proof.} From the assumptions there exist group monomorphisms $\delta: \Theta \rightarrow \Phi$, $\omega: \Omega \rightarrow \Upsilon$  and characters $\eta:\Theta \rightarrow \mathbb T$, $\nu:\Omega \rightarrow \mathbb T$ such that \begin{equation} \begin{split}\eta(\g) xu_\g x^* = v_{\delta(\g)} \text{ for all }\g\in \Theta\\
\nu(\g) y u_\g y^* = v_{\omega(\g)} \text{ for all }\g\in \Omega
\end{split}\end{equation}

\noindent Using this in combination with the hypothesis we see that for every $\sigma \in \Sigma$ we have that $u_t= u_{\phi(\sigma) ^{-1}} u_t u_\sigma =  d_\sigma  y^*v_{\omega (\phi(\sigma))^{-1}} yu_t x^* v_{\delta(\sigma)} x$ and hence  \begin{equation}\label{intelement} yu_tx^*=   d_\sigma  v_{\omega (\phi(\sigma))^{-1}} yu_t x^* v_{\delta(\sigma)},
\end{equation}     
where we have denoted by $d_\sigma = \nu (\phi(\sigma)) \eta(\sigma^{-1})$.

\noindent Consider the Fourier decomposition $yu_tx^*= \sum_\la c_\la v_\la$ with respect to $\cN =\cL(\La)$. Using this in   \eqref{intelement} we see that $c_\la= d_\sigma c_{\omega (\phi(\sigma))\la \delta(\sigma)^{-1}}$ for every $\sigma \in \Sigma$, $\la \in \Lambda$. As  $|d_\sigma|=1$ it follows that $|c_\la|= |c_{\omega (\phi(\sigma))\la \delta(\sigma)^{-1}}|$ for every $\sigma \in \Sigma$, $\la \in \Lambda$. This implies that $c_\la\neq 0$ only if the set $O_\la=\{ \omega (\phi(\sigma))\la \delta(\sigma)^{-1} \,:\, \sigma\in \Sigma \}$ is finite. This implies that there exists a finite index subgroup  $\Sigma_\la\leqslant \Sigma$ such that  $\omega (\phi(\sigma))\la \delta(\sigma)^{-1}=\la$ and hence  \begin{equation}\label{partinnisom}\omega (\phi(\sigma))=\la \delta(\sigma) \la ^{-1}\text{ for all }\sigma\in \Sigma_\la.\end{equation} 

\noindent Next, let $\la,\mu \in \La$ such that $O_\la$ and $O_\mu$ are finite.  By \eqref{partinnisom}, for every $\sigma \in \Sigma_\la \cap \Sigma_\mu$ we get   $\mu \delta(\sigma) \mu ^{-1}= \delta (\phi(\sigma))=\la \delta(\sigma) \la ^{-1}$ and hence $\mu^{-1}\la \in C_\La( \delta(\Sigma_\la \cap \Sigma_\mu) )$; in particular,  the unitary  $u_{\mu^{-1}\la}\in \cL(\Sigma_\la\cap \Sigma_\mu)'\cap \cN =\cL(\Sigma_\la\cap \Sigma_\mu)'\cap \cL(\G)$. Now, notice that since $[\Sigma :\Sigma_\la], [\Sigma:\Sigma_\mu]<\infty$ then $[\Sigma: \Sigma_\la\cap\Sigma_\mu]<\infty$. Therefore,  $\cL(\Sigma_\la\cap \Sigma_\mu)'\cap \cL(\G)\subseteq \cL(vC_\G(\Sigma))$ and since $vC_\G(\Sigma)=1$  we conclude that $\cL(\Sigma_\la\cap \Sigma_\mu)'\cap \cL(\G)=\mathbb C1$. In particular, this implies that $u_{\mu^{-1}\la} =1$ and hence $\la=\mu$. Altogether, this shows that there exists $\la_0\in \La$ such that $yu_tx^* = c_{\la_0} v_{\la_0}$. As $u_t$ is unitary we have that $|c_{\la_0}|=1$, as desired.\hfill$\blacksquare$

\begin{corollary}\label{C:hnn} Let $\Sigma\leqslant  \Theta\leqslant \G$ be  groups such that $vC_\G(\Sigma)=1$. Let $\phi: \Sigma \rightarrow \Theta$ be a group homomorphism   and $t\in \G$ such that $\phi(\sigma ) = t \sigma t^{-1}$ for all $\sigma\in \Sigma$.   Let $\La$ be an arbitrary group such that $\cN=\cL (\G)=\cL(\La)$. Assume there exist a subgroup $\Phi\leqslant \La$ and a unitary $w\in \mathcal U(\cN)$ such that $\mathbb T w \Theta w^*= \mathbb T \Phi$. Then one can find a subgroup $\Phi < \Xi \leqslant \La$ such that  $\mathbb T w\langle \Theta, t\rangle w^*= \mathbb T \Xi$.
    
\end{corollary}
\noindent {\it Proof.}
Applying Theorem \ref{discretizationhnngroups} result for $\Omega=\Theta$, $\Upsilon=\Phi$ and $x=y=w$, there exist $c\in \mathbb T$ and $\la \in \La$ such that $w u_t w^*=cv_\la$. Thus, the result follows by letting $\Xi= \langle \Phi, \la \rangle\leqslant \La$. 
\hfill$\blacksquare$

\begin{corollary} Let  $\Sigma \leqslant \G$ be groups so that $vC_\G (\Sigma)=1$ and ${\rm QN}^{(1)}_\G(\Sigma)=\G$.  Let $\La$ be an arbitrary group such that $\cN=\cL (\G)=\cL(\La)$. Assume there exist a subgroup $\Upsilon \leqslant \La$ and a unitary $u\in \mathcal U(\cN)$ so that $\mathbb T u \Sigma u^*= \mathbb T \Upsilon$. Then we have $\mathbb T  u \G u^*= \mathbb T \La$.
    
\end{corollary}

\noindent {\it Proof.} Fix an arbitrary $g\in \G={\rm QN}^{(1)}_\G(\Sigma)$. Since $[\Sigma: g\Sigma g^{-1}\cap\Sigma]<\infty$, we derive that $vC_\G(g\Sigma g^{-1}\cap\Sigma)=vC_\G(\Sigma)=1$. By applying Corollary \ref{C:hnn} for the map $\phi:g\Sigma g^{-1}\cap\Sigma\to \Sigma$ defined by $\phi(\sigma)=g^{-1}\sigma g$ for all $\sigma\in g\Sigma g^{-1}\cap\Sigma$, we deduce that there is a subgroup $\Upsilon_g<\Lambda$ such that $\mathbb T u \langle \Sigma,g \rangle u^*= \mathbb T \Upsilon_g$. Since this is true for any $g\in\Gamma$, we deduce that $\mathbb T u \Gamma u^*= \mathbb T (\vee_{g\in\Gamma}  \Upsilon_g)=\mathbb T \Lambda$.
\hfill$\blacksquare$

\section{Identification of ``peripheral structure''}
\noindent 

In this section we present several technical results that will be used for reconstructing the factor subgroups of $G$ from the von Neumann algebra $\cL(G)$, whenever $G$  belongs to $\mathscr D$ or $\mathscr A$. Specifically, this means that whenever $\cL(G)=\cL(H)$, for an arbitrary group $H$, then for every factor subgroup $G_0<G$ one can find a subgroup $H_0<H$ and a unitary $u\in \cL(G)$ such that $\cL(G_0) = u\cL(H_0)u^*$. In other words, we can identify a collection of subgroups of $H$ that up to von Neumann equivalence play the same role as the factor subgroups of $G$. Such a result is called \emph{identification of peripheral structure}. Our first result introduces sufficient conditions for identification of peripheral structure, up to corners.  

\begin{theorem}\label{identificationcornergroups1} Let $G$ be an icc non-amenable group and denote $\cM=\cL (G).$ Let $K\leqslant H\leqslant G$ be an inclusion of icc non-amenable groups satisfying the following properties:\begin{enumerate}
    \item ${\rm QN}_G(K)={\rm QN}^{(1)}_G(H)=H$, 
    and $[H:KC_H(K)]<\infty$;
    \item For every $p\in \mathscr P(\cL(H))$ and every von Neumann algebra with no amenable direct summand $\cA \subseteq p\cL(H)p$ such that  $\cA'\cap p\cM p$ is diffuse, we have that $\cA'\cap p\cM p \subseteq p\cL(H)p$;
    \item For every $p_i\in \mathscr P(\cL(H))$ and von Neumann algebras $\cA_i \subseteq p_i\cL(H)p_i$ with $i=1,2$ such that $\cA_1$ has no amenable direct summand and $\cA_1'\cap p_1 \cL(H)p_1$ is diffuse, if $\cA_1 \prec_\cM \cA_2$ then $\cA_1\prec_{\cL(H)}\cA_2$; 
    \item For every $p\in \mathscr P(\cL(H))$, whenever $\cD,\cE,\cF \subseteq p\cL(H)p$ are mutually commuting von Neumann subalgebras so that $\cD$ is isomorphic to a corner of $\cL(K)$ and $\cE$ has no amenable direct summand then $\cF$ is purely atomic.
 \end{enumerate}

\noindent Let $\La$ be an arbitrary group such that $\cM=\cL(\La)$ and assume there is a subgroup $\Omega<\La$ with non-amenable centralizer $C_\La(\Omega)$ such that $\cL(K)\prec_\cM  \cL(\Omega)$. 

\noindent Then one can find a subgroup $\Omega C_\La(\Omega)\leqslant {\rm QN}_\La(\Omega) \leqslant \Sigma<\La$ with  $[\Sigma : {\rm QN}_\La (\Omega)]<\infty$ and ${\rm QN}^{(1)}_\La(\Sigma)=\Sigma$, a non-zero projection $c\in \mathcal Z(\cL(\Sigma))$ and $w_0\in \mathscr U(\cM)$ with $w_0c w^*_0=n\in\cL(H)$  such that $w_0 \cL(\Sigma)c w^*_0= n\cL (H)n$.
\end{theorem}

\noindent {\it Proof.} Since $\cL(K)\prec_\cM  \cL(\Omega)$ one can find projections $a\in \cL(K)$, $f\in \cL(\Omega)$, a non-zero partial isometry $v\in f\cM a$ and a $\ast$-isomorphism onto its image $\phi:a \cL(K)a \rightarrow \cB:=\phi(a\cL(K)a) \subseteq f\cL(\Omega)f $ such that \begin{equation}\label{inteq1}
    \phi(x)v=vx \text{ for all } x\in a\cL(K)a.
\end{equation}
Notice that $vv^*\in \cB'\cap f\cM f$ and $v^*v\in a\cL(K)a'\cap a\cM a$. The equation \eqref{inteq1}  implies that $\cB vv^* = v \cL(K)v^* = u_1 \cL(K) v^*v u^*_1$, where $u_1\in \cM$ is a unitary extending $v$. Taking relative commutants we get $vv^*(\cB' \cap f\cM f)vv^*= u_1 v^*v( (a\cL(K)a)'\cap a\cM a )v^*v u^*_1$. 
Condition (1) implies that ${\rm vC}_G(K)\leqslant H$, and hence, $\cL(K)'\cap \cM\subset \cL({\rm vC}_G(K)) \subset \cL(H)$. It follows that $vv^*(\cB \vee \cB'\cap f\cM f)vv^*=\cB vv^*\vee vv^*(\cB' \cap f\cM f)vv^*\subseteq u_1 \cL(H) u^*_1$. Therefore, since $\cL(H)$ is a factor one can find a new unitary $u_2 \in \mathscr U(\cM)$ such that \begin{equation}(\cB \vee \cB'\cap f\cM f) z_2 \subseteq u_2 \cL(H) u^*_2,\end{equation} where $z_2$ is the central support of $vv^*$ in $\cB \vee \cB'\cap f\cM f$. In particular, we have $z_2\in \mathcal Z(\cB'\cap f\cM f )$ and $vv^*\leq z_2 \leq f$.

\noindent Observe that $\cL(C_\La(\Omega))z_2\subseteq ((f\cL(\Omega )f)'\cap f \cM f )z_2\subseteq (\cB'\cap f \cM f)z_2\subseteq u_2 \cL(H)u_2^*$. As $z_2 \in( \cL(C_\La(\Omega))f)'\cap f\cM f$ then using hypothesis (2) we further have that   $z_2 (\cL(C_\La(\Omega))f\vee (\cL(C_\La(\Omega))f)'\cap f\cM f) z_2\subseteq u_2\cL(H)u_2^*$. Again since $\cL(H)$ is a factor there is $u\in \mathscr U(\cM)$ so that 
\begin{equation}\label{inteq2} (\cL(C_\La(\Omega))f\vee (\cL(C_\La(\Omega))f)'\cap f\cM f) z\subseteq u\cL(H)u^*,\end{equation} where $z$ is the central support of $z_2$ in $\cL(C_\La(\Omega))f\vee (\cL(C_\La(\Omega))f)'\cap f\cM f$. In particular, we have $vv^*\leq z_2 \leq z \leq f$. Now, since  $f\cL(\Omega )f\subseteq  (\cL(C_\La(\Omega))f)'\cap f\cM f$, then by \eqref{inteq2} we get $(f\cL(\Omega )f\vee \cL(C_\La(\Omega))f)z\subseteq u\cL(H)u^*$ and hence 
\begin{equation}\label{inteq3} u^*(\cL(C_\La(\Omega))f\vee f \cL(\Omega )f) z u\subseteq \cL(H).\end{equation}
Since $vv^*\leq z\in (f\cL(\Omega)f)'\cap f\cM f$ and $\cB$ is a factor then the map $\phi': a\cL(K)a \rightarrow u^* \cB z  u \subseteq f \cL(\Omega)fz$ given by $\phi'(x)=u^* \phi(x)z u$ still defines a $\ast$-isomorphism that satisfies $\phi'(x) w=wx$, for any $x\in a\cL(K)a$, where $w= u^*zv$ is a non-zero partial isometry. Hence, $\cL(K) \prec_{\cM}  u^*f\cL(\Omega)fzu$. By the hypothesis (3), it follows that $\cL(K)\prec_{\cL(H)} u^*f\cL(\Omega)fzu$.

\noindent To this end using \cite[Proposition 2.4]{CKP14} and its proof we can find non-zero $p \in \mathscr P(\cL(K))$, $r=u^*ezu \in u^*f\cL(\Omega)fzu$ with $e\in \mathscr P(f\cL(\Omega)f)$, a von Neumann subalgebra $\cC \subseteq  u^*e\cL(\Omega)ezu $ , and a $\ast$-isomorphism $\theta : p\cL(K)p \rightarrow \cC$ such that the following properties are satisfied:
\begin{enumerate}
\item [a)] the inclusion $\cC \vee (\cC' \cap  u^*e\cL(\Omega)ezu ) \subseteq  u^*e\cL(\Omega)ezu $  has finite index in the sense of Pimsner-Popa \cite{PP86};
\item [b)]there is a non-zero partial isometry $y \in \cL(H)$ such that $\theta(x) y =y x$ for all $x\in p\cL(K)p$, where $y^*y\in p\cL(K)p'\cap p\cM p\subset p\cL(H)p$ and $yy^*\in \cC '\cap r\cM r$. 
\end{enumerate}

Note that $\cC$, $\cC' \cap   u^* e\cL(\Omega)e z u$ and $u^*\cL(C_\La(\Omega))e z u$ are commuting subalgebras of $u^*ez \cL(H) ze u$. 
Since $C_\La(\Omega)$ is non-amenable, it follows that
$u^*\cL(C_\La(\Omega))ezu$ has no amenable direct summand. Hence, since $\cC$ is isomorphic to a corner of $\cL(K)$, it follows from hypothesis (4)  that $\cC' \cap  u^* e\cL(\Omega)ezu$ is purely atomic. Thus, one can find a non-zero projection $q \in \mathcal Z(\cC' \cap   u^* e\cL(\Omega)e zu )$ such that $(\cC' \cap   u^* e\cL(\Omega)e zu )q=\mathbb C q$. Hence, after compressing the containment in a) by $q$ and replacing $\cC$ by $\cC q$, $y$ by $qy$  and  $\theta(x)$ by $\theta(x) q$ in b) we can assume in addition that  $\cC \subseteq  u^*e\cL(\Omega)ezu $  is a finite index inclusion of non-amenable II$_1$ factors. By \cite[Proposition 1.3]{PP86}, it follows that $\cC \subseteq u^* e\cL(\Omega)e zu$ admits a finite Pimsner-Popa basis, which implies that there exist $x_1,\dots,x_m\in u^* e\cL(\Omega)e zu$ such that $u^* e\cL(\Omega)e zu=\sum_{i=1}^m x_i \cC$. 
Note also that $u^* e\cL(\Omega)e zu \subset r\cM r$ since $r=u^*ezu$. Hence,
\begin{equation}\label{jj}
 \mathcal {QN}_{ r\cM r}(\cC ) ''=\mathcal {QN}_{ r\cM r}( u^* e\cL(\Omega)ezu )''.   
\end{equation}
Also, the intertwining relation in b) shows that $\cC yy^*= y p\cL(K)p y^*= l p\cL(K)p y^*y l^*$, where $yy^*\in \cC'\cap r \cM r$ and $l\in \cL(H)$ is a unitary extending $y$ (i.e. $y=ly^*y$). Therefore, using the quasi-normalizer formulas for group von Neumann algebras  and for compressions, Lemma \ref{QN2}  and Lemma \ref{QN1} (and Remark \ref{QN.remark}), respectively,  we see that
\begin{equation}\label{eqquasicorner}\begin{split}
  ly^*y \cL(H) y^*yl^* & = ly^*y \cL({\rm QN}_G (K)) y^*yl^*   \overset{L. \ref{QN2}}{=}l y^*y \mathcal{QN}_{\cM}(\cL(K))''y^*yl^*\\
  & \overset{L. \ref{QN1}}{=}\mathcal {QN}_{l y^*y \cM y^*yl^*}(l p\cL(K)p y^*y l^*)''= \mathcal {QN}_{yy^*\cM yy^*}(\cC yy^*)''
  \\
  &\overset{R. \ref{QN.remark}}{=}yy^*\mathcal {QN}_{ r\cM r}(\cC ) ''yy^*\overset{\eqref{jj}}{=}yy^*\mathcal {QN}_{ r\cM r}( u^* e\cL(\Omega)ezu )'' yy^*\\
  &\overset{L. \ref{QN1}}{=}yy^*  u^*z e\mathcal {QN}_{ \cL(\La)}(\cL(\Omega) )''ez u yy^*\overset{L. \ref{QN2}}{=}yy^*  u^* z e\cL ({\rm QN}_{\La}(\Omega) )e zu  yy^*.
\end{split}\end{equation}

\noindent In the above equalities, we also used that $y\cM y^*=yy^*\cM yy^*$ and $u^* e\cL(\Omega)ezu = r u^* \cL(\Omega) ur$.

\noindent Next, denote by $\Upsilon = {\rm QN}_\La(\Omega)$ and by $\Upsilon<\Sigma = \langle {\rm QN}^{(1)}_\La(\Upsilon)\rangle <\La$. As ${\rm QN}^{(1)}_G(H)= H$, then formula \eqref{eqquasicorner} together with the corresponding formulas for one-sided quasinormalizers, Lemma \ref{QN1} and Lemma \ref{QN2}, show that \begin{equation}\label{eqquasicorner2}yy^* u^* ze\cL (\Upsilon)e zu  yy^* =yy^*  u^* z e\cL (\Sigma)ez u yy^*= ly^*y \cL(H) y^*yl^*.\end{equation}
 In particular, by \cite[Lemma 2.2]{CI17} we have $[\Sigma:\Upsilon]<\infty$, and hence, ${\rm QN}^{(1)}_\La(\Sigma )={\rm QN}^{(1)}_\La(\Upsilon )=\Sigma$. By applying Lemma \ref{QN2}, we obtain
 \begin{equation}\label{to.apply.CI17}
     \mathcal{QN}_{\cM}^{(1)}(\cL(\Sigma))=\cL(\Sigma).
 \end{equation}
\noindent Notice that relation \eqref{eqquasicorner2}
also show that $yy^*= u^* d  u$ for some projection $d\in ze\cL (\Sigma)ez$, and thus,  \eqref{eqquasicorner2} entails that $ u^*  d\cL (\Sigma) d u = ly^*y \cL(H) y^*yl^*$. By letting $w_0:=ul \in \mathscr U(\cM)$ and $t:=y^*y\in \cL(H)$,  we conclude that $w^*_0 d \cL(\Sigma)d w_0= t \cL(H) t$. Since $l^*yy^*l= w_0^* d w_0$ and $yy^*=ly^*$, we have $t=w^*_0 d w_0$. Moreover, if we replace $w^*_0 \Sigma w_0$ by $\Sigma$ and use $w^*_0 d w_0= t$, we have that $t \cL(\Sigma)t = t \cL(H)t$. 

Next, as $\cL(H)$ is a factor one can find a unitary $w_1 \in \cM$  such that  $\cL(\Sigma) c \subseteq w_1\cL(H)w_1^*$, where $c$ denotes the central support of $t \in \cL(\Sigma )$. In particular, it follows that there exists a projection $h \in \cL(H)$ such that $t= w_1 hw_1^*$. Moreover, since $\cL(H)$ is a factor and $t\in \cL(H)$, there is a unitary $w_2 \in \cL(H)$ so that $t = w_2 hw^*_2$. Altogether, these relations show that $wt=tw$, where $w :=w_1 w_2^*$. Also, note that $\cL(\Sigma)c\subseteq w \cL(H)w^*$. Multiplying on both sides by $t$, we get 
\begin{equation}\label{eq.j}
t \cL(H)t=t\cL(\Sigma)t\subseteq t w\cL(H)w^* t.    
\end{equation}
Multiplying on the left by $tw^*$ and using that $tw=wt$, we further obtain that $tw^*t \cL(H) t\subseteq t \cL(H)t w^* t$. In particular, using the hypothesis (1) together with Lemmas \ref{QN1} and \ref{QN2}, we get $w^*t =tw^*t\in \mathcal{QN}^{(1)}_{t\cM t}(t\cL(H)t)=t \cL(H)t$, and hence, $wt \in t \cL(H)t$. Using relation \eqref{eq.j}, we deduce that $t \cL(\Sigma)t= t\cL(H)t=w t \cL(H)t w^*$. By using relation \eqref{to.apply.CI17} and $\cL(\Sigma) c\subseteq cw \cL(H)w^*c$, we apply the moreover part of \cite[Lemma 2.6]{CI17} and derive that $\cL(\Sigma) c= c w \cL(H) w^*c$. This shows the desired conclusion. \hfill$\blacksquare$



\noindent Next, we present a technical result that will be needed to derive the main superrigidity results for groups $G$ in the classes $\mathscr D$ and $\mathscr A$. The result is essentially showing that if $G_0$ is a factor subgroup of $G$, then we can upgrade from recognizing some subalgebras in corners of $\cL(G_0)$, to recover the corners of $\cL(G_0)$. 
The proof of this result generalizes the proof of \cite[Theorem 3.2]{CI17} and for reader's convenience we include all the details.

\begin{theorem}\label{unitpart2} Assume that $G$ is a group in one of the following classes: \begin{enumerate}
\item [(i)]$G \in \mathscr D$ and let $f(G)=\{G_1,...,G_n\}$ be its canonical factors.
\item [(ii)] $G= (\G_1\ast \G_2\ast ...\ast \G_n) \rtimes_\rho \G \in \mathscr A$ and denote by $G_i = \G_i\rtimes_{\rho^i} \G$ for all $1\leq i\leq n$.
\end{enumerate} 

\noindent Let $\La$ be a group such that $\cN=\cL(G)=\cL(\La)$. In addition, assume that for every $i\in \{1,...,n\}$ there is a subgroup  $\La_i \leqslant \La$ satisfying the following relations: \begin{enumerate}
    \item $\La_i$ contains two commuting non-amenable subgroups;
    \item ${\rm QN}^{(1)}_\La(\La_i)=\La_i$;
    \item There is a subset $i\in J_i \subseteq \{1,...,n\}$, projections $0\neq z^i_k \in \mathscr Z(\cL(\La_i))$ with $k\in J_i$ which satisfy $\sum_{k\in J_i}z^i_k= 1$; 
    \item There is $u_i\in \mathscr U(\cN )$  such that \begin{equation}\begin{split} &u_i \cL( \La_i) z^i_i u_i^*=  p_i \cL(G_i) p_i , \text{ where }p_i= u_i z^i_i u_i^*\in \mathscr P(\cL(G_i))\text{, and } \\ &u_i \cL( \La_i) z^i_j u_i^*\subset  \cL(G_j) \text{ for all }j\in J_i \setminus \{i\}.
    \end{split} 
\end{equation} 
\end{enumerate}
Then one can find  a partition  $T_1\sqcup \cdots\sqcup T_l=\{1,...,n\}$ and for every $1\leq i\leq l$ there is a subgroup  $\Omega_i \leqslant \La$  such that the following relations hold: \begin{enumerate}
    \item ${\rm QN}^{(1)}_\La(\Omega_i)=\Omega_i$;
    \item There are projections $0\neq \tilde z^i_k \in \mathscr Z(\cL(\Omega_i))$ with $k\in T_i$ which satisfy $\sum_{k\in T_i}\tilde z^i_k= 1$; 
    \item There exist $\tilde u_i\in \mathscr U(\cN )$  such that $\tilde u_i \cL( \Omega_i)\tilde z^i_k \tilde u_i^*=  p^i_k \cL(G_k) p^i_k $, where $p_k^i=\tilde u_i\tilde z^i_k \tilde u_i^*\in \mathscr P(\cL(G_k))$, for every $k\in T_i$.
 
\end{enumerate}
\end{theorem}

\noindent {\it Proof.} From the hypothesis,  for every $1\leq i\leq n$  there is a subgroup $\La_i \leqslant \La$ such that ${\rm QN}^{(1)}_\La(\La_i)=\La_i$. Moreover,  there exist $i\in J_i \subseteq \{1,...,n\}$ and  $0\neq z^i_k \in \mathscr P(\mathscr Z(\cL(\La_i)))$, $k\in J_i$ which satisfy \begin{equation}\label{unitpart3}\sum_{k\in J_i}z^i_k= 1.\end{equation} and there is $u_i\in \mathscr U(\cN )$  such that  \begin{equation}\label{cont2}\begin{split} &u_i \cL( \La_i) z^i_i u_i^*=  p_i \cL(G_i) p_i , \text{ where }p_i= u_i z^i_i u_i^*\in \mathscr P(\cL(G_i))\text{, and } \\ &u_i \cL( \La_i) z^i_k u_i^*\subset  \cL(G_k) \text{ for all }k\in J_i \setminus \{i\}.
\end{split} 
\end{equation}  
Fix $k\in J_1$. Since $u_1z^1_ku_1^*, p_k\in \cL(G_k)$ are non-zero projections and $\cL(G_k)$ is a factor one can find a projection $\cL(G_k)\ni q\leq u_1z_k^1u_1^*$ and a unitary  $w\in \cL (G_k)$ such that $wqw^*\leq p_k$. This together with relations \eqref{cont2} imply that $q u_1 \cL( \La_1) z^1_k u_1^* q \subseteq  q\cL(G_k) q \subseteq w^* p_k\cL(G_k)p_k w   = w^* u_k \cL(\La_k) z_k^k u_k^* w$. In particular, we have  $\cL(\La_1)\prec_\cN \cL(\La_k)$ and hence by \cite[Lemma 2.6]{CI17} there is a $h_k\in \La$ such that $[\La_1: h_k \La_k h^{-1}_k\cap \La_1]<\infty$.  Therefore replacing $\La_k$ by $h_k \La_k h^{-1}_k$, relations \eqref{cont2} still hold and in addition we can assume that \begin{equation}\label{finiteindexgr}[\La_1: \La_k\cap\La_1]<\infty \text{ for all } k\in J_1.\end{equation}

\noindent Next, we claim that for all $i\in J_1,j\in J_k$ satisfying $i\neq j$ we have  $z_i^1z_j^k=0$. To see this, assume by contradiction there are $J_1\ni i\neq j\in J_k$ so that $z^1_i z_k^j\neq 0$. Relations \eqref{cont2} give that $u_1 \cL( \La_k\cap \La_1) z^1_i u_1^*\subseteq  \cL(G_i)$ and $u_k \cL( \La_k\cap \La_1) z^k_j u_k^*\subseteq   \cL(G_j)$. Thus, by \cite[Lemma 2.4]{CI17} one can find $g\in G$ such that \begin{equation}\label{intertwiningintersection}\cL(\La_k\cap \La_1)\prec \cL(G_i \cap g G_k g^{-1}).\end{equation} 
From here we treat each case separately. First, assume (i).  Since $i\neq j$  then $G_i \cap g G_j g^{-1}$ is amenable which forces $\La_k\cap \La_1$ amenable. Using \eqref{finiteindexgr} we would get that $\La_1$ is amenable which is a contradiction. Next, assume (ii). Since $i\neq j$  then $G_i \cap g G_j g^{-1}$ is either trivial or it is equal to $d(\G)$. First possibility is obviously impossible so assume that $G_i \cap g G_k g^{-1}=d(\G)$. From condition (1) and \eqref{finiteindexgr} it follows that $\La_k \cap \La_1$ contains two non-amenable commuting subgroups. This however together with \eqref{intertwiningintersection} contradicts the fact that $d(\G)$ is bi-exact. This concludes our claim. 

\noindent Now, fix an arbitrary $i\in J_1$. By  \eqref{unitpart3} there exists $l\in J_k$ such that $z_i^1z^k_l\neq 0$. Then the above claim implies that $l=i$ and also $z^1_i\leq z_i^k$ by using once again \eqref{unitpart3}. In particular, we have $J_1\subseteq J_k$. Arguing by symmetry we conclude that \begin{equation}\label{equality1} J_1=J_k\text{ and  }z^1_i= z_i^k \text{ for all } i\in J_1.
\end{equation}
Next, we notice that relations \eqref{cont2} imply $u_1 \cL(\La_k\cap \La_1)z^1_i u_1^* \subseteq \cL(G_i)$ and $u_k \cL(\La_k\cap \La_1)z^k_i u_k^* \subseteq \cL(G_i)$ for all $i\in J_1$. We continue by arguing that in both cases we have \begin{equation}\label{nointcore}u_1 \cL(\La_k\cap \La_1)z^1_i u_1^* , u_k \cL(\La_k\cap \La_1)z^k_i u_k^* \nprec_{\cL(G_i)} \cL(G_i\cap G_j)\text{ for all }j\neq i.\end{equation} In case (i) this is immediate because $G_i\cap G_j$ is amenable while $\La_k\cap \La_1$ is not.  In case (ii) this follows because $G_i \cap G_j= d(\G)$ is a bi-exact group while, as before, $\La_k\cap \La_1$ contains two commuting non-amenable subgroups.

\noindent Next, since $(u_1z_i^k u_k^*) u_k \cL(\La_k\cap \La_1)z^k_i u_k^*= u_1 \cL(\La_k\cap \La_1)z^1_i u_1^* (u_1z_i^k u_k^*)$ then \eqref{nointcore} and \cite[Theorem 1.2.1]{IPP05} further imply that $u_1z_i^k u_k^*\in \cL(G_i)$ for all $i\in J_1$.  Thus conjugating the algebras in \eqref{cont2} by $u_1z_i^k u_k^*$ we can assume without any loss of generality that $u_k=u_1$. In particular, we have $u_1 \cL(\La_1)z^1_1 u_1^*= p_1 \cL(G_1)p_1$ and $u_1 \cL(\La_k)z^k_1 u_1^*\subset  \cL(G_1)$.  Since by \eqref{equality1} we have $z_1^i=z^k_1$ and $u_1z_1^1u_1^*=p_1$ these relations imply that  $u_1 \cL(\La_k)z^k_1 u_1^*\subset  p_1\cL(G_1)p_1=u_1 \cL(\La_1)z^k_1 u_1^*$. Thus $\cL(\La_k)z^k_1 \subset  \cL(\La_1)z^k_1 $ and moreover $\cL(\La_k\cap \La_1)z^k_1\subseteq \cL(\La_k)z^k_1 \subseteq  \cL(\La_1)z^k_1 $. Thus by \eqref{finiteindexgr} the inclusion $\cL(\La_k\cap \La_1)z^k_1\subseteq \cL(\La_k)z^k_1$ admits a finite Pimsner-Popa basis and hence $[\La_k: \La_k\cap \La_1]<\infty$ by \cite[Lemma 2.6]{CI17}. However, this combined with  \eqref{finiteindexgr} and  part (2) imply that $\La_k=\La_1$. Altogether, these and \eqref{cont2} show that $u_1 \cL(\La_1) z_1^1 u_1^*=p_1 \cL(G_1)p_1$ and $u_1 \cL(\La_1) z^1_k u_1^*=p_k \cL(G_k)p_k$. Since the above arguments work for all $k\in J_1$, letting $T_1:=J_1$ and $p^1_k :=p_k $ for $k\in T_1$ we get the statement for the first element of the partition. Also we let $\Omega_1:=\La_1$.

\noindent If $T_1 =\{1,...,n\}$ the proof is completed. If not, pick the smallest $s\in \{1,...,n\}\setminus T_1$ and repeat the same arguments as before starting with set $J_s$, etc. We leave the details to the reader.
\hfill$\blacksquare$



\section{Superrigidity results for groups in class $\mathscr D$}

\noindent In this section we derive the main superrigidity results for von Neumann algebras $\cL(G)$ associated with groups in class $G\in \mathscr D$. These add numerous examples of $W^*$-superrigid amalgams to the ones found in \cite{CI17} and also give the first examples of $W^*$-superrigid HNN extensions. Next we highlight a  point of contrast between our results and the prior results \cite{CI17}. In the one hand, the factor subgroups $f_i(G)$ of our groups $G$ consist of direct products of groups in the class $\mathcal{IPV}$ which is slightly more restrictive when compared to the factor subgroups of the amalgamated free products considered in \cite{CI17}. One the other hand, the amalgamated subgroups $a_i(G)$ of our $G$ are more general, exceeding the specific case of diagonal subgroups used in \cite{CI17}. For instance, our examples cover amalgamated subgroups which are of direct product type themselves.

\noindent The following is the main result of the section and it should be compared with \cite[Theorem A]{CI17}.  
\begin{theorem}\label{tt1} Let $G\in \mathscr D$ and let $f(G)=\{G_1,G_2,...,G_n\}$ be its factor set. Assume that $\La$ is an arbitrary group such that $\cM=\El(G) =\El(\La)$.

Then one can find a unitary $u_i\in \mathscr U(\cM)$ and a subgroup $\La_i<\La$ such that $u_i \cL(\La_i )u_i^*= \cL(G_i)$ for any $1\leq i\leq n$.   

\end{theorem}

\noindent {\it Proof.} First, if $G\in \mathscr{D_0}$, then the conclusion follows from Corollary \ref{superrigidproducts}. Hence, we assume $G\in \mathscr{D}_m$ for some $m\ge 1$ and denote $a_m(G)=\Sigma$.
Next, we prove the following claim:

{\bf Claim. }\label{eqpartcorners6} For every $1\leq i\leq n$ there exists a subgroup $\La_i<\La$ satisfying the following relations: \begin{enumerate}
    \item $\La_i$ contains two commuting non-amenable subgroups;
    \item ${\rm QN}^{(1)}_\La(\La_i)=\La_i$;
    \item There is a subset $i\in J_i \subseteq \{1,...,n\}$ and projections $0\neq z^i_k \in \mathscr Z(\cL(\La_i))$ with $k\in J_i$ so that $\sum_{k\in J_i}z^i_k= 1$; 
    \item There exist $u_i\in \mathcal U(\cN )$  such that  
\begin{enumerate} \item $ u_i \cL( \La_i) z^i_i u_i^*=  p_i \cL(G_i) p_i $, where $p_i= u_i z^i_i u_i^*\in \mathscr P(\cL(G_i))$, and \item $u_i \cL( \La_i) z^i_k u_i^*\subset  \cL(G_i ) $ for all $k\in J_i \setminus \{i\}$.
\end{enumerate}   
\end{enumerate}

\noindent\textit{Proof of the claim}. Let $G_i=\Gamma^i_1\times\dots\times\Gamma^i_{s_i}$ be a product of groups that belong to $\mathcal{IPV}$, for any $1\leq i\leq n$. We denote by $\Gamma^i_{\hat j}$ the product group $\times_{k
     \in \hat j}\Gamma_k^i$, for any $1\leq j\leq s_i.$  
Consider $\Delta : \cM \rightarrow \cM \bar\otimes \cM$ the comultiplication along $\La$, given by $\Delta (v_\la)=v_\la \otimes v_\la$ for all $\la \in \La$. Fix $1\leq i \leq n$ and consider the inclusion $\Delta(\cL (G_i))\subseteq \cM \bar\otimes \cM = \cL(G  \times G)$. Fix $1\leq j\leq s_i$ and denote by $\cA:=\Delta (\cL(\G_{\hat j}^i))$ and $\cB = \Delta(\cL(\G_j^i))$ and notice that $\cA$ and $\cB$ are commuting von Neumann algebras with no amenable direct summand. Thus by \cite[Proposition 7.2(4)]{IPV10} it follows that $\cA \nprec_{\cM\bar\otimes \cM} \cL(G\times A)$, $\cA \nprec_{\cM\bar\otimes \cM} \cL(A\times G)$, $\cB \nprec_{\cM\bar\otimes \cM} \cL(G\times A)$ and $\cB \nprec_{\cM\bar\otimes \cM} \cL(A\times G)$, for every amenable subgroup $A<G$. Hence, Theorem \ref{Th:commute} further implies that one can find $1\leq k\leq n$ and $1\leq l\leq s_k$ such that 
\begin{equation}
    \cA \prec_{\cM\bar\otimes \cM} \cL(G \times \Gamma^k_{\hat l})\text{ or }\cL( \Gamma^k_{\hat l}\times G).
\end{equation}

\noindent By symmetry it suffices to treat only the first case.  Using Theorem \ref{Th:ultrapower} one can find a non-amenable subgroup $\Omega < \La$ with non-amenable centralizer $C_\La(\Omega)$ so that \begin{equation}
 \cL(\G_{\hat j}^i)\prec_\cM \cL(\Omega).   
\end{equation}

\noindent To this end, we notice that if we let $K=\G^{\hat j}_i$ we see that all the conditions (1)-(4) in the Theorem \ref{identificationcornergroups1} are satisfied. Therefore, the conclusion of Theorem \ref{identificationcornergroups1} implies that there exist a subgroup ${\rm QN}^{1}_\Lambda(\Omega)<_f\La_i <\La$ with ${\rm QN}^{(1)}_\La(\La_i)=\La_i$ and $\La_i$ contains two commuting non-amenable subgroups,  a central projection $z^i_i\in \cL(\La_i)$ and a unitary $u^i_i\in \cM$ with $t^i_i=u^i_i z_i(u^i_i)^*\in \cL(G_i)$ so that \begin{equation}\label{asd1}
u^i_i \cL(\La_i) z^i_i (u^i_i)^* = t^i_i \cL(G_i)t^i_i.
\end{equation}  
\noindent Now, for every $1\leq k\leq n$ let  $y^i_k \in \mathscr Z(\cL(\La_i))$ be the maximal projection for which there  exists a unitary $v_k^i\in \cM$ such that  
\begin{equation}\label{asd2}
v^i_k\cL(\La_i) y^i_k (v^i_k)^* \subseteq  t_k^ i\cL(G_k) t_k^ i,  \end{equation}
where $t_k^i=v_k^iy_k^i (v_k ^i)^*$. It is an easy exercise to see that since $\cL(G_k)$ is a factor such projections always exist. We also notice that $y^i_k$ are mutually orthogonal. Indeed, otherwise by \cite[Lemma 2.4]{CI17} we would get that $\cL(\La_i)\prec \cL(\Sigma)$ which is not possible as $\La_i$ is non-amenable while $\Sigma$ is amenable.

\noindent Next, we show that $\sum_k y^i_k=1$. Towards this let $z=1- \sum_k y_k^i\in \mathscr Z(\cL(\La_i))$ and assume by contradiction that $z\neq 0$. Now since $\Omega$ and $C_\La(\Omega)$ are commuting non-amenable subgroups and $G$ is bi-exact relative to $f(G)$ (in the sense of  \cite[Definition 15.1.2]{BO08}), then by \cite[Theorem 15.1.5]{BO08} there is $1\leq l\leq n $ such that $\cL(\Omega )z\prec_\cM \cL(G_l)$. Thus, one can find some projections $r \in \cL(\Omega) z$, $q\in \cL(G_l)$, a non-zero partial isometry $v\in q\cM r$ and a $\ast$-isomorphism onto its image $\Phi: r\cL(\Omega)r\rightarrow \cD:=\Phi(r\cL(\Omega)r)\subseteq q\cL(G_l)q$ such that $\Phi(x)v=vx$ for all $x\in r\cL(\Omega)r$. Since $\Sigma$ is an amenable group and $\cD$ has no amenable summand, then $\cD\nprec \cL(\Sigma)$. Therefore, by using Lemma \ref{L:ipp} we have $\cD'\cap q\cM q \subseteq q\cL(G_l)q$. In particular, $vv^*\in q\cL(G_l)q$ and hence the intertwining relation implies that $vr\cL (\Omega)  rv^*\subseteq \cL(G_l)$. Thus since $\cL(G_l)$ is a factor one can find a unitary $w\in \cM $ such that $w \cL(\Omega )r_o w^*\subseteq \cL(G_l)$ where $r_o$ is the central support of $r$ in $\cL(\Omega)$. From the conclusion of Theorem \ref{identificationcornergroups1} we have that ${\rm QN}^{(1)}_\La({\rm QN}_\La^1(\Omega))=\La_i$ and therefore repeating  the same arguments as before (two times) one can find a new unitary $w_1\in \cM$ and a projection $a\in \mathscr Z(\cL(\La_i))$ 
with $a\geq r_o\geq r$ such that $ w_1\cL(\La_i)a w_1^*\subseteq \cL(G_l)$.  Notice by construction we have that $0\neq b:=az\in \mathscr Z (\cL(\La_i))$. In particular, $by_l^i=0$ and since
$ w_1\cL(\La_i)b w_1^*\subseteq \cL(G_l)$, $ v^i_l\cL(\La_i)y_l^i (v^i_l)^*\subseteq \cL(G_l)$ and $\cL(G_l)$ is a factor one can perturb $w_1$ to a new unitary such that that there exists $t\in \mathscr U(\cM)$ satisfying $t \cL(\La_i)(y_l^i+ b)t^*\subseteq \cL(G_l)$. This obviously contradicts the maximality of $y^i_l$, so $z=1$. 

\noindent We continue by showing that $z_i^i=y_i^i$. From construction we have that   $z_i^i\leq y_i^i$ and assume by contradiction  $c:=y_i^i-z_i^i\neq 0$. Notice that $ v^i_i\cL(\La_i)c (v^i_i)^*\subseteq \cL(G_i)$ and $ u^i_i\cL(\La_i)z_i^i (u^i_i)^*= t^i_i\cL(G_i)t_i^i$.  Since $\cL(G_i)$ is a factor we can perturb $v_i^i$ to a new unitary so that there is a projection $e\in \cL(G_i)$ satisfying $e t_i^i=0$ and $ v^i_i\cL(\La_i)c (v^i_i)^*\subseteq e\cL(G_i)e$. Thus, the element $f= v_i^ic+ u_i^iz_i^i$ satisfies $f^*f= c+z_i^i$, $ff^*= e+ t_i^i$ and $f \cL(\La_i) (c+z_i^i) f^*\subseteq (e+t_i^i) \cL(G_i)(e+t_i^i)$. Also, by using Lemma \ref{QN2} we have that $\mathscr{QN}^{(1)}_{f \cM f }( f \cL(\La_i) (c+z_i^i) f )= f \cL(\La_i) (c+z_i^i) f$ and then obviously $\mathscr{QN}^{(1)}_{ (e+t_i^i) \cL(G_i )(e+t_i^i)}( f \cL(\La_i) (c+z_i^i) f )= f \cL(\La_i) (c+z_i^i) f$. Therefore, since $(e+t_i^i) \cL(G_i )(e+t_i^i)$ is a factor, applying the moreover part in \cite[Lemma 2.6]{CI17}  we conclude that  $f \cL(\La_i) (c+z_i^i) f^*= (e+t_i^i) \cL(G_i)(e+t_i^i)$. However this is impossible as the center of the algebra on the left-hand side is two-dimensional while the center the right-hand side one is one-dimensional. 
In conclusion, $z_i^i=y_i^i$.

Next, we denote $z^i_k:= y^i_k$, $u_k^i:=v_k^i$ and $J_i=\{1\leq k\leq n \,:\, z_k^i\neq 0\}$. Since the common part $\cL(\Sigma)$ is a II$_1$ factor by perturbing the $(u_k^i)$'s to new unitaries one can assume that $(t^i_k)_k\subset \cL(\Sigma)$ are mutually orthogonal projections satisfying $\sum_k t^i_k=1$. These relations imply that $u_i= \sum_j u^i_k z_k^i\in \cM $ is a unitary and moreover the equations \eqref{asd1} and \eqref{asd2} entail that $u_i \cL(\La_i) z^i_k u^*_i\subseteq \cL(G_k)$ for all $k\neq i$ and $u_i\cL(\La_i) z^i_i u^*_i= t^i_i \cL(G_i)t^i_i$. This concludes the proof of the claim.
$\hfill\square$

\noindent To this end, we note that the Claim together with Theorem \ref{unitpart2} imply that one can find  a partition  $J_1\sqcup \cdots\sqcup J_l=\{1,...,n\}$ and for every $1\leq i\leq l$ there is a subgroup  $\Omega_i \leqslant \La$  such that the following relations hold: \begin{enumerate}
    \item ${\rm QN}^{(1)}_\La(\Omega_i)=\Omega_i$;
    \item There are projections $0\neq z^i_k \in \mathscr Z(\cL(\Omega_i))$ with $k\in J_i$ which satisfy $\sum_{k\in J_i}z^i_k= 1$; 
    \item There exist $u_i\in \mathscr U(\cN )$  such that $ u_i \cL( \Omega_i) z^i_k u_i^*=  p^i_k \cL(G_k) p^i_k $, where $p^i_k= u_i z^i_k u_i^*\in \mathscr P(\cL(G_k))$.
 
\end{enumerate}
 Next, we claim that $l=n$ and each set $J_i$ consists of a singleton. Indeed, assume by contradiction there is i such that $|J_i|\geq 2$. Also by replacing $\Lambda$ by $u_i \Lambda u_i^*$ we can assume that $\cL(\Omega_i)z_k^i= z_k^i \cL(G_k) z_k^i$ for all $k\in J_i$. Note that $\bigcap_{k\in J_i}G_k=\Sigma$. Now using the same argument from the proof of \cite[Proposition 4.1]{CI17} together with Theorem \ref{malnormalcontrol} one obtains a contradiction with $\Sigma$ being icc.  
 Thus, the $J_i$'s are singletons and therefore for every $1\leq i\leq n$ one can find a unitary $u_i \in \cM$ and an icc subgroup $\Omega_i<\La$ with ${\rm QN}^{(1)}_\La(\Omega_i)=\Omega_i$ such that $u_i \cL(G_i)u_i^*=\cL(\Omega_i)$.
\hfill$\blacksquare$

\begin{theorem}\label{superrigidtreegroups}
Let $G \in \mathscr D$.  Assume that $\La$ is an arbitrary group and let $\theta: \El(G) \to \El(\La)$ be a $*$-isomorphism. Then  there exist $\delta \in {\rm Isom}(G,\La)$, $ \omega \in {\rm Char}(G)$ and  $u \in \mathcal{U}(\cL(\La))$ such that $\theta = {\rm ad}(u) \circ \Psi_{\omega,\delta}$. 
\end{theorem}

\noindent {\it Proof.} Let $f(G)=\{G_1,...,G_n\}$ be the factor set of $G$. From the hypothesis we have $\theta(\cL(\G))=\cL(\La)$ and using the previous theorem one can find for each $1\leq i\leq n$ a unitary $u_i\in \cM$ and subgroup $\La_i<\La$ such that $u_i\cL(G_i)u_i^*=\cL(\Lambda_i)$. Using Theorem \ref{ipv10} after perturbing the $u_i$'s to new unitaries we have that $\mathbb T u_i \theta (G_i)u_i^*= \mathbb T \La_i$. Then using Theorem \ref{discretizationgeneratinggroups} and Corollary  \ref{C:hnn} iteratively after finitely many steps one can find a unitary $u\in\cM$  such that $\mathbb T u \theta (G)u^*= \mathbb T \La$, which gives the desired conclusion.   
\hfill$\blacksquare$

\noindent This result also implies that the groups in class $\mathscr D$ are completely recognizable from the $C^*$-setting as well. This adds a new class of non-amenable $C^*$-superrigid groups to the only other previously known \cite{CI17,CD-AD21}.

\noindent Next we record two immediate applications of the prior result.  For the definition of $\cup _{i\geq 1}\mathcal D^m_i$ we encourage the reader to revisit Section \ref{treesuperrigid}.

\begin{corollary}  Let $G\in \cup_{i\geq 1}\mathscr D_i^m$. Assume that $\La$ is an arbitrary group and let $\theta: C^*_r(G) \to C^*_r(\La)$ be a $*$-isomorphism. Then  there exist $\delta \in {\rm Isom}(G,\La)$, $ \omega \in {\rm Char}(G)$ and  $u \in \mathcal{U}(\cL(\La))$ such that $\theta = {\rm ad}(u) \circ \Psi_{\omega,\delta}$.  
\end{corollary}

\noindent {\it Proof.} Note that $G$ has trivial amenable radical by Proposition \ref{P:amenable radical}. Then it follows from \cite{BKKO14} that $C^*_r(G)$ has unique trace and thus $\theta$ lifts to a $\ast$-isomorphism of the corresponding von Neumann algebras  $\theta: \cL(G) \to \cL(\La)$. The statement follows then from the previous theorem. \hfill$\blacksquare$

\begin{corollary}  Let $G\in \cup_{i\geq 1}\mathscr D_i^m$. Then for any $\theta\in {\rm Aut}( C^*_r(G))$ there exist $\delta \in {\rm Aut}(G)$, $ \omega \in {\rm Char}(G)$ and  $u \in \mathscr{U}(\cL(G))$ such that $\theta = {\rm ad}(u) \circ \Psi_{\omega,\delta}$.  
\end{corollary}

\section{Superrigidity results for groups in class $\mathscr A$}
\noindent In this section we show that the semidirect product groups in class $\mathscr A$ are both $W^*$ and $C^*$- superrigid. These add new examples, to the prior ones found in \cite{IPV10,BV12,Be14}. The only known examples of $C^*$-superrigid groups were found in \cite{CI17}. Our results provide a second such class and the first
of semidirect product type.  After this work was completed, the authors found a second class of $C^*$-superrigid semidirect products, this time from the realm of generalized wreath product groups with nonamenable core, \cite{CD-AD21}.

\begin{theorem}\label{eqfactors}\label{tt2} Let $G=(\G_1 \ast \G_2 \ast ...\ast \G_n) \rtimes_{ \rho}\G \in \mathscr A$. Assume that $\La$ is an arbitrary group such that $\cM=\El(G) =\El(\La)$. 

Then one can find a unitary $u_i\in \mathscr U(\cM)$ and a subgroup $\La_i<\La$ such that $u_i \cL(\La_i )u_i^*= \cL(\G_i \rtimes_{\rho^i} \G)$ for any $1\leq i\leq n$.

\end{theorem}

\noindent {\it Proof.} The proof will be obtained in several steps. Some of them follow directly from the prior results in \cite{CI17} while for the others we include detailed proofs. We encourage the reader to consult beforehand \cite[Theorem A and Proposition 4.1]{CI17} as some parts of the proofs rely heavily on these results.  We start by proving the following:  

{\bf Claim.}\label{eqpartcorners1} For every $1\leq i\leq n$ there exists a property (T) subgroup $\La_i<\La$ satisfying the following relations: \begin{enumerate}
    \item ${\rm QN}^{(1)}_\La(\La_i)=\La_i$;
    \item There is a subset $i\in J_i \subseteq \{1,...,n\}$ and projections $0\neq z^i_k \in \mathscr Z(\cL(\La_i))$ with $k\in J_i$ so that $\sum_{k\in J_i}z^i_k= 1$; 
    \item There exist $u_i\in \mathcal U(\cM )$  such that  
\begin{enumerate} \item $ u_i \cL( \La_i) z^i_i u_i^*=  p_i \cL(\G_i\rtimes \G) p_i $, where $p_i= u_i z^i_i u_i^*\in \mathscr P(\cL(\G_i\rtimes \G))$, and \item $u_i \cL( \La_i) z^i_k u_i^*\subset  \cL(\G_k \rtimes \G) $ for all $k\in J_i \setminus \{i\}$.
\end{enumerate}   
\end{enumerate}

\noindent{\it Proof of the claim}. Denote by $G^1_i\times G^2_i=G_i = \G_i\rtimes_{\rho^i} \G$ where $G_i^j\cong \G$. By using \eqref{amalgamdecomp}, we view $G$ as an amalgam $G= G_1\ast_\Sigma G_2\ast_\Sigma... \ast_\Sigma G_n$ where $\Sigma=d(\G)$.   Consider $\Delta : \cM \rightarrow \cM \bar\otimes \cM$ the commultiplication along $\La$, given by $\Delta (v_\la)=v_\la \otimes v_\la$ for all $\la \in \La$. Fix $1\leq i \leq n$ and consider the inclusion $\Delta(\cL (G_i))\subseteq \cM \bar\otimes \cM = \cL((\ast^j_\Sigma G_j)  \times (\ast^l_\Sigma G_l))$. Using   \cite[Theorem 5.1]{IPP05} and the fact that $G_i$ has property (T), there exist $1\leq j,l \leq n$, a projection $0\neq z\in \Delta (\cL(G_i))'\cap \cM\bar\otimes \cM$ and $u\in \mathscr U(\cM\bar\otimes \cM)$ such that \begin{equation}
    u \Delta(\cL(G_i))z u^*\subseteq \cL(G_j \times G_l).
\end{equation}   
Since $G^1_j, G^2_j, G^1_l, G^2_l$ are bi-exact we get by Theorem \ref{Th:bo} that  there exist $1\leq k,t\leq 2$ such that 
$\Delta(\cL(G^k_i))z \prec \cL(G_j \times G^t_l)$ or $ \Delta(\cL(G^k_i))z \prec \cL(G^t_j \times G_l)$. Due to symmetry, it suffices to treat only one of these possibilities; thus, assume  $\Delta(\cL(G^k_i))z \prec_{\cM\bar\otimes \cM} \cL(G_j \times G^t_l)$. Using Theorem \ref{Th:ultrapower} one can find a non-amenable subgroup $\Omega < \La$ with non-amenable centralizer $C_\La(\Omega)$ so that \begin{equation}
 \cL(G^k_i)\prec_\cM \cL(\Omega).   
\end{equation}

\noindent Next, we notice that if we let $G_i^k= K$ then all conditions (1)-(4) in the Theorem \ref{identificationcornergroups1} are satisfied. Therefore, the conclusion of Theorem \ref{identificationcornergroups1} implies that there exist a subgroup $\Omega C_\La(\Omega )<\La_i <\La$ with ${\rm QN}^{(1)}_\La(\La_i)=\La_i$,  a central projection $z^i_i\in \cL(\La_i)$ and a unitary $u^i_i\in \cM$ with $t^i_i=u^i_i z_i(u^i_i)^*\in \cL(G_i)$ so that \begin{equation}\label{eqcorner3}u^i_i \cL(\La_i) z^i_i (u^i_i)^* = t^i_i \cL(G_i)t^i_i.\end{equation}  

\noindent Since $\cL(G_i)$ has property (T), then \eqref{eqcorner3}, \cite[Lemma 2.13]{CI17} and \cite{CJ85} show that $\La_i$ is a property (T) group as well. Thus, using \cite[Theorem 5.1]{IPP05} again for every $j\neq i$ one can find projections $z^i_j \in \cL(\La_i)'\cap  \cM$ with $\sum_{j\neq i} z^i_j =1-z^i_i$,  unitaries $u^i_j\in \cM$ and projections $t_j^i\in \cL(G_j)$ such that \begin{equation}\label{contcorner1}
    u^i_j\cL(\La_i)z^i_j (u^i_j)^*\subseteq t_j^i\cL(G_j)t_j^i.
\end{equation}  
Also notice that since the common part $\cL(\Sigma)$ is a II$_1$ factor by perturbing the $(u_j^i)$'s to new unitaries one can assume that $(t^i_j)_j\subset \cL(\Sigma)$ are mutually orthogonal projections satisfying $\sum_j t^i_j=1$. These relations imply that $u_i= \sum_j u^i_j z_j^i\in \cM $ is a unitary and moreover the equations \eqref{eqcorner3} and \eqref{contcorner1} entail that $u_i \cL(\La_i) z^i_j u^*_i\subseteq \cL(G_j)$ for all $j\neq i$ and $u_i\cL(\La_i) z^i_i u^*_i= t^i_i \cL(G_i)t^i_i$. By letting $J_i=\{1\leq j\leq n \,:\, z_j^i\neq 0\}$, this concludes the proof of the claim. $\hfill\square$

\noindent To this end, we note that the Claim together with Theorem \ref{unitpart2} imply that one can find  a partition  $J_1\sqcup \cdots\sqcup J_l=\{1,...,n\}$ and for every $1\leq i\leq l$ there is a property (T) subgroup  $\La_i \leqslant \La$  such that the following relations hold: \begin{enumerate}
    \item ${\rm QN}^{(1)}_\La(\La_i)=\La_i$;
    \item There are projections $0\neq z^i_k \in \mathscr Z(\cL(\La_i))$ with $k\in J_i$ which satisfy $\sum_{k\in J_i}z^i_k= 1$; 
    \item There exist $u_i\in \mathscr U(\cM )$  such that $ u_i \cL( \La_i) z^i_k u_i^*=  p^i_k \cL(\G_k\rtimes_{\rho^k}\G) p^i_k $, for any $k\in J_i$, where $p_k^i= u_i z^i_k u_i^*\in \mathscr P(\cL(G_k))$.
 
\end{enumerate}
 Next, we claim that $l=n$ and each set $J_i$ consists of a singleton. Indeed, if we assume that for some $i$ we have $|J_i|\geq 2$ then applying verbatim the arguments from  the proofs of \cite[Proposition 4.1 and Theorem A]{CI17} one obtains a contradiction. We leave the details to the reader. In particular, our claim entails that for every $1\leq i\leq n$ there is a unitary $u_i \in \cM$ so that $u_i \cL(\La_i)u_i^*= \cL(\G_i\rtimes \G )$. 
 \hfill$\blacksquare$

\noindent Now, we are ready to derive the main results of the section.

\begin{theorem} Let $G=(\G_1 \ast \G_2 \ast ...\ast \G_n) \rtimes_{ \rho}\G\in \mathscr A$. Assume that $\La$ is an arbitrary group and let $\theta: \El(G) \to \El(\La)$ be a $*$-isomorphism. Then  there exist $\delta \in {\rm Isom}(G,\La)$, $ \omega \in {\rm Char}(\G)$ and  $u \in \mathscr{U}(\cL(\La))$ such that $\theta = {\rm ad}(u) \circ \Psi_{\omega,\delta}$.  
\end{theorem}

\noindent {\it Proof.} From the hypothesis we have that $\theta(\cL(G))= \cL(\La)$. Thus by Theorem \ref{eqfactors}  there exists a unitary $u\in \mathscr U(\cM)$ and for each $1\leq i\leq n$ there is a subgroup $\La_i<\La$ so that $\theta (\cL(\G_i\rtimes_{\rho^i} \G ))=u \cL(\La_i)u^*$. Therefore $\La$ admits an amalgam decomposition $\La = \La_1 \ast_\Omega \La_2\ast_\Omega .... \ast_\Omega \La_n$ and viewing $\G_i\rtimes_{\rho^i} \G$ as $\G_i \times \G_i$ and the acting group as the diagonal group $d(\G)$ we have that $\theta (\cL(\G_i\times \G_i)=u \cL(\La_i)u^*$  and  $\theta(d(\G))= u\cL(\Omega) u^*$ for all $1\leq i\leq n $. By Corollary \ref{identificationdiagonalgroups} there is a unitary $u_i\in \cL (\La_i)$ such that $\mathbb T \theta (\G_i \times \G_i) = \mathbb Tu_i \La_i u_i^*$. Therefore using Lemma \ref{quasinormalizerdiagonal} and Theorem \ref{discretizationgeneratinggroups} recursively one can find a unitary $u \in\mathscr U(\cM)$ such that $\mathbb T \theta(G)=\mathbb T u \La u^*$. This gives the desired conclusion.\hfill$\blacksquare$

\noindent This result also implies that the groups in class $\mathscr A$ are completely reconstructible from the $C^*$-setting as well. This adds a new class of nonamenable $C^*$-superrigid groups to the only other previously known, \cite{CI17}. 

\begin{corollary}  Let $G\in \mathscr A$. Assume that $\La$ is an arbitrary group and let $\theta: C^*_r(G) \to C^*_r(\La)$ be a $*$-isomorphism. Then  there exist $\delta \in {\rm Isom}(G,\La)$, $ \omega \in {\rm Char}(G)$ and  $u \in \mathcal{U}(\cL(\La))$ such that $\theta = {\rm ad}(u) \circ \Psi_{\omega,\delta}$.  
\end{corollary}

\noindent {\it Proof.} Note that $G$ has trivial amenable radical by Proposition \ref{P:classA}. Then it follows from \cite{BKKO14} that $C^*_r(G)$ has unique trace and thus $\theta$ lifts to an $\ast$-isomorphism of the corresponding von Neumann algebras  $\theta: \cL(G) \to \cL(\La)$. The statement follows then from the previous theorem. \hfill$\blacksquare$

\begin{corollary}  Let $G\in \mathscr A$. Then for any $\theta\in {\rm Aut}( C^*_r(G))$ there exist $\delta \in {\rm Aut}(G)$, $ \omega \in {\rm Char}(G)$ and  $u \in \mathscr{U}(\cL(G))$ such that $\theta = {\rm ad}(u) \circ \Psi_{\omega,\delta}$.  
\end{corollary}

\section*{Appendix}

\noindent In this appendix we provide an alternative proof of the direct product rigidity Theorem \ref{prodbiexactrig} for groups in the class $\mathcal{IPV}$ which by-passes the usage of prior techniques for bi-exact groups. 
This approach builds upon the methods developed in \cite{OP08,CPS11,CdSS15}.

\noindent A key ingredient for our proof is a structural result which classifies all weak compact embeddings into tensor products by wreath product von Neumann algebras in the same spirit with some results in \cite{CPS11}. In fact this result does not appear anywhere in the literature and deserves some attention on its own. This is one of the main reasons we decided to include this appendix in the paper. To properly introduce the result we first recall briefly the definition of a weakly compact action introduced in \cite{OP07}. 
\begin{definition}\label{OP}
Let $\cA \subseteq \cM$ be an inclusion of tracial von Neumann algebras and consider a subgroup of normalizers $\mathscr H  \leqslant \mathscr{N}_\cM(\cA)$. Then the  conjugation action $\mathscr H \ca \cA$ is called {\it weakly compact} if we can find a net $ \eta_n \in L^2(\cA \bar\otimes \bar{\cA})$ of positive unit vectors satisfying the following conditions: 
\begin{enumerate}
    \item $\lim_n{\|(a \otimes \overline{a}) \eta_n - \eta_n\|_2}=0$,  for all $a \in \mathscr{U}(\cA)$, 
    \item $\lim_n{\|[u \otimes \overline{u},\eta_n]\|_2} =0$,  for all $u \in \mathscr H$,
    \item $\langle (x \otimes 1) \eta_n , \eta_n \rangle = \langle (1 \otimes \overline x) \eta_n, \eta_n \rangle = \tau(x)$, for all $n$ and $x \in \cA$.
\end{enumerate}
\end{definition}
 
\noindent  With this definition at hand we are now ready to state and prove the result.

\begin{theorem}\label{convweakcomp}
Let $K_0<G_0$ and $A$ be some countable groups such that $K_0$ and $A$ are amenable.
Let $\cM$ be a finite von Neumann algebra and denote by $\cN:=\cM \bar\otimes \cL(A \wr_{K_{0}} G_{0})$. Let $\cB \subset p\cN p$ be a diffuse von Neumann subalgebra and let $\mathscr H \subset \mathcal{N}_{p\cN p}(\cB)$ a subgroup of normalizers such that the natural action by conjugation $\mathscr H \ca \cB$ is weakly compact. 

If the von Neumann algebra $\mathscr{H} ''$ is strongly non-amenable relative to $\cM \otimes 1$, then the deformation $1 \otimes \alpha_t \to id$ uniformly on the unit ball $(\cB)_1$. 

Here the path $\alpha_t$ is the wreath product core-length deformation on $\cL(A \wr_{G_0/K_{0}} G_0)$ introduced in \cite{Io06} (see also \cite{IPV10}). 
\end{theorem}

\noindent {\it Proof.}
Our proof is similar to \cite[Theorem B]{OP08} and \cite[Theorem 6.2]{CPS11}. For the proof we can assume without loss of generality that $p \in \cM \bar\otimes \cL(G_0)$ and therefore $(1 \otimes \alpha_t)(p)=p$, for all $t \in \mathbb{R}$. Let $z_0\in \mathcal Z(\mathscr H'\cap p\cN p)$ be the maximal projection such that $\alpha_t\to id$ uniformly on the unit ball of $\cB z_0$. 
Assume by contradiction that $z_0\neq 1$ and take an arbitrary non-zero projection  $z\in \mathcal Z(\mathscr H'\cap p\cN p)$ with $z\leq 1-z_0$. This implies that $1 \otimes \alpha_t$ does not converge uniformly on $\mathscr{U}(\cB z)$. Using the transversality property of $\alpha_t$ from \cite[Lemma 2.1]{Po08} there exists $c > 0$ and sequences $t_k \searrow 0$ and $(u_kz)_k \subset \mathscr{U}(\cB z)$ such that \[
\|(1 \otimes \alpha_{t_{k}})(u_kz) - E_{\cN}((1 \otimes \alpha_{t_{k}})(u_kz))\|_2 \geq c\|z\|_2 \text{, } \text{ for all }k \in \mathbb{N} \text{.}
\]
Using Pythagoras's theorem we get that \[
\|E_{\cN}((1 \otimes \alpha_{t_{k}})(u_kz))\|_2 \leq \sqrt{1 -c^2}\|z\|_2 \text{, } \text{ for all }k \in \mathbb{N} \text{.}
\]

Now, pick $0 < \delta < \frac{1-\sqrt{1- c^2}}{6}\|z\|_2$. Choose and fix $k \in \mathbb{N}$ such that $\alpha = \alpha_{t_k}$ satisfies the following relations
\begin{enumerate}
    \item [a.] $\|z - \alpha(z)\|_2 < \delta$
\end{enumerate}
Let $v = u_k$ and let $(\eta_n)_n$ be a net of vectors as in Definition \ref{OP} which corresponds to the weakly compact action $\mathscr H\car \cB$  and consider the following notations 
\begin{enumerate}
    \item [b.] $\tilde{\eta}_{j,n}=(\alpha_{t_j} \otimes {\rm id})(\eta_n) \in L^2(\tilde{\cN}) \otimes L^2({\cN})$ 
    \item [c.] $b_{j,n} = (e_{\cN} \otimes 1)(\tilde{\eta}_{j,n}) \in L^2(\cN) \otimes L^2(\cN)$
    \item [d.] $b_{j,n}^{\perp} = \tilde{\eta}_{j,n} - b_{j,n} \in (L^2(\tilde{\cN}) \ominus L^2(\cN)) \otimes L^2(\cN)$.
\end{enumerate}
For ease of notation, denote $\tilde{\eta}_n=\tilde{\eta}_{k,n}$, $b_n=b_{k,n}$ and $b_n^\perp= b_{k,n}^{\perp}$. Notice that $(p \otimes 1) b_n^{\perp} (p \otimes 1) = b_n^{\perp}$ and 
\[
\|(xp \otimes 1)\tilde{\eta}_n\|_2 = \tau(E_{\cN}(\alpha^{-1}(px^*xp))) = \|xp\|_2^2.
\]
Also, as in the proof of \cite[Theorem 4.9]{OP08} we get \[
\|[u \otimes \bar{u}, b_n^{\perp}]\|_2 \leq {\|(\alpha \otimes 1)([u \otimes \bar{u},\eta_n])\|_2}+2\|u-\alpha(u)\|_2, \text{ for all }u\in\mathscr U.
\]
Next, we claim that \[\lim_n\|(z \otimes 1)b_n^{\perp}\|_2 \geq \delta.\] 
Assume this is not the case, since $e_\cN z = ze_\cN$ and $zv=vz$ we get that 
\[
\lim_n{\|(z \otimes 1) \tilde{\eta}_n - (e_\cN \alpha(v) z \otimes \overline{v})b_n\|_2} \leq \lim_{n}{\|(z \otimes 1) \tilde{\eta}_n - (e_\cN \alpha(v) z \otimes \overline{v})\tilde{\eta}_n\|_2} + \lim_{n}{\|(z \otimes 1) b_{n}^{\perp}\|_2}
\]
\[
\leq \lim_{n}{\|(z \otimes 1) \tilde{\eta}_n - (e_{\cN}z \alpha(v) \otimes \overline{v}) \tilde{\eta}_n\|_2} + \|[\alpha(v),z]\|_2 + \delta \leq \lim_{n}{\|\alpha \otimes 1 (\eta_n - (v \otimes \overline v) \eta_n)\|_2} + 4 \delta = 4 \delta \text{.}
\]
Therefore,  this further implies that
\[
\|E_{\cN}(\alpha(vz))\|_2 \geq \|E_{\cN}(\alpha(v))z\|_2 - \|z - \alpha(z)\|_2 \geq \lim_{n}{\|E_\cN(\alpha(v))z \otimes \overline{v} \tilde{\eta}_n\|_2} - \delta 
\]
\[
\geq \lim_{n}{\|(e_\cN \alpha(v)z \otimes \overline{v}) b_n\|_2} - \delta \geq \lim_{n}{\|(z \otimes 1) \tilde{\eta}_n\|_2} - 5 \delta \geq \|z\|_2 - 5 \delta \geq \sqrt{1-c^2}\|z\|_2
\]
which contradicts 
$\|E_{\cN}((1 \otimes \alpha_{t_{k}})(u_kz))\|_2 \leq \sqrt{1 -c^2}\|z\|_2 \text{, } \text{ for all }k \in \mathbb{N} \text{.}$ This concludes the proof of the claim.

\noindent Pick $n$ large enough such that $b = b_n^{\perp} \in (L^2(\tilde{\cN}) \ominus L^2(\cN))  \otimes L^2(\cN)$. For any $x \in \cN$ we have that \[
\|(x \otimes 1) b^{\perp}\|_2^2 = \|(x \otimes 1) (e_\cM^{\perp} \otimes 1) \tilde{\eta}_n\|_2^2 = \|(e_\cM^{\perp} \otimes 1) (x \otimes 1) \tilde{b}_n\|_2^2 \leq \|(x \otimes 1) \tilde{\eta}_n\|_2^2 = \|x\|_2^2\text{.}
\]

\noindent Next, we employ an argument similar with the proof of \cite[Theorem B]{OP08}. Denote by $\cK = L^2(\tilde{\cN}) \ominus L^2(\cN)$ and notice that it is an $\cN$-bimodule with the natural left and right action by $\cN$. Also, consider the von Neumann algebra $\cP= B(\cK) \cap \rho(\cN^{op})'$, where $\rho(\cN^{op})$ is the right action on $\cK$. Let $\eta_{n,k} = \|(z\otimes 1)b_{k,n}^{\perp}\|^{-1}(z\otimes 1) b_{k,n}^{\perp}$ and consider the functional $\phi_k : \cP \to \mathbb{C}$ given by $\phi_{k}(x) = \lim_{n}{\langle (x \otimes 1) \eta_{k,n}, \eta_{k,n} \rangle}$. Now, $\phi_k$ is a well-defined state on $\cP$ satisfying $\phi_k(zx) = \phi_k(xz) = \phi_k(x)$, for all $x \in \cP$.

\noindent Now, we prove the following:

{\bf Claim.}
For every $y \in \mathscr{H}''$ we have that 
\[
\lim_{k}{|\phi_k(xy)-\phi_{k}(yx)|} = 0,
\]
uniformly for $x \in (\cP)_1$.

{\it Proof of the Claim.}
Fix $u \in \mathscr{H}$. Then, for every $x \in \cP$ we have 
\begin{equation*}
\begin{split}
&|\phi_k(uxu^*) - \phi_k(x)| = \lim_{n}{|\langle (uxu^* \otimes 1) \eta_{k,n} , \eta_{k,n} \rangle - \langle (x \otimes 1) \eta_{k,n} , \eta_{k,n}\rangle|}\\
&= \lim_{n}{\frac{1}{\|b_{k,n}\|^2}|\langle (x \otimes 1)(u \otimes \overline{u}) b_{k,n}^{\perp}(u^* \otimes \overline{u^*}) , (u \otimes \overline{u}) b_{k,n}^{\perp}(u^* \otimes \overline{u^*}) \rangle - \langle (x \otimes 1) b_{k,n}^{\perp} , b_{k,n}^{\perp} \rangle|}
\\
&\leq 2\|x\|_{\infty} \lim_{n}{\frac{\|[u \otimes \overline{u} , b_{k,n}^{\perp}]\|_{2}}{\|b_{k,n}^{\perp}\|}} \leq \frac{4}{\delta} \|x\|_{\infty}  \|u - \alpha_{k}(u)\|_2 \text{.}
\end{split}
\end{equation*}
Thus, for ever $y \in {\rm span} \mathscr{H}$ we have that \[
\lim_{k}{|\phi_k(yx) - \phi_k(xy)|} = 0,
\]
uniformly on $x \in (\cP)_1$.

\noindent Using (5) one can check that 
\[
\lim_{k}{|\phi_k(xy)|} \leq \lim_{k}{\frac{1}{\|b_{k,n}^{\perp}\|^2}|\langle (xy \otimes 1) b_{k,n}^{\perp} , b_{k,n}^{\perp}\rangle|} \leq \frac{1}{\delta} \|x\|_{\infty} \|y\|_2 \text{.}
\]
The same inequality can be proven for $\phi_k(yx)$ and using Kaplansky's density theorem we get the claim.
\hfill$\square$

We notice that by the calculation done in \cite[Lemma 4.2]{CPS11} we have that $\mathcal{K}=L^{2}(\tilde{\cN}) \ominus L^{2}(\cN)) \simeq \oplus_{s}{L^{2}(\langle \cN , e_{\cM\bar\otimes\cK_{s}} \rangle )}$ where $\cK_{s} = \cL(A^{I -{\Delta_{s}}} \wr {\rm stab}_{G_{0}}(\tilde{\eta}_{s}))$ where $\Delta_{s}$ is the support of $\tilde{\eta}_{s}$. Therefore, using Connes fusion we have $\cK \simeq \oplus_{s}\left[{L^2(\cN) \otimes_{\cM \otimes \cK_{s}} L^2(\cN)}\right]$. Since, $\cK_{s}$ is amenable we get that $L^{2}(\cN) \otimes_{\cM \bar\otimes \cK_{s}} L^2(\cN)$ is weakly contained in $L^2(\cN) \otimes_{\cM} L^2(\cN)$ for every $s \in S$. Therefore, $\cK$ is weakly contained in $L^{2}(\cN) \otimes_\cM L^2(\cN)=:\mathcal{T}$.  Using \cite[Lemma A.3]{Is16} one can find a ucp map \[
\Phi: \cQ := \mathcal{B}(\mathcal{T}) \cap \rho(\cN^{op})' \to \mathcal{B}(\mathcal{K}) \cap \rho(\cN^{op})'=\mathcal{P}
\]
such that $\Phi(\lambda_{\mathcal{T}}(n))=\lambda_{\mathcal{K}}(n)$ for all $n \in \cN$ and the sub script denotes the actions of $\cN$ on $\mathcal{T}$ and $\mathcal{K}$, respectively. Now, consider the state $\psi_{k}: \cQ \to \mathbb{C}$ given by $\psi_k = \phi_k \circ \Phi$. Since the left action is in the multiplicative domain of $\Phi$ using the Claim, for every $u \in \mathcal{H}''$ we have that \[
\lim_{k}{|\psi_{k}((uz)^*xuz-x)|}=\lim_{k}{|\phi_k(\Phi((uz)^*xuz)-\Phi(x))|}=\lim_{k}{|\phi_{k}(u\Phi(x)u^*-\Phi(x))|} = 0 
\]
uniformly for $x \in (\cQ)_1$. 
Now, using a standard averaging argument in conjunction with Hahn-Banach separation theorem and the functional calculus one can find $\beta_k \in L^1(\cQ)_{+}$ such that $0 \leq \cE(\beta_k) \leq 1$ and for all $u \in\mathscr{H}''$ we have \[
\lim_{k}{\|\beta_k - (uz)^*\beta_k uz\|_1} =0 \text{.}
\]
Here $\cE: L^1(\cQ) \to L^1(\cN)$ is the canonical map such that $\tau(\cE(s)x)=Tr(sx)$, for all $x \in \cN, s \in \cQ$. Using an appropriate normalization we can assume that $\beta_k = z \beta_k z$ and $\|\beta_k\|_1 = 1$. Letting $z_k = \beta_k^{1/2}$ and using the generalized Power-Stormer inequality we further get \[
\lim_{k}{\|z_k - (uz) z_k (uz)^*\|_2} = 0,
\]
for all $u \in \mathscr{H}''$. Now, fix $F \subset \mathscr{H}''$ an arbitrary finite subset. Using the identification $L^2(\cQ) = \mathcal{T} \otimes_{\cM} \mathcal{T}$, assuming $z_k \in \mathcal{T} \otimes_{\cM} \mathcal{T}$ and using the above equality, we get that \[
\begin{split}
|F| &=\lim_k \|\sum_{u \in F}{z_k}\|_2 \leq \lim_k\sum_{u \in F}{\|z_k - (uz)z_k(uz)^*\|_2} + \lim_k\|\sum_{u\in F}{(uz)z_k(uz)^*}\|_2\\ &\leq \lim_k\|\sum_{u\in F}{uz \otimes \overline{uz}}\|_{\mathcal{T} \otimes_{\cM} \mathcal{T}}.
\end{split}
\]
Since this holds for all $z \in \mathcal{Z}(\mathscr{H}' \cap \cN)$ and all $F \subset \mathscr{H}''$ finite it follows from \cite[Corollary 2.4]{Si11} that $\mathcal{T}$ is a left amenable $\cN-\cN$ bimodule over $\mathscr{H}''$. Since $\mathcal{T}$ can be identified to $L^2(\cN) \otimes_{\cM \otimes 1} L^{2}(\cN)$, then $\mathscr{H}''$ is amenable relative to $\cM \otimes 1$ inside $\cN$, which contradicts our assumption. 
\hfill$\blacksquare$

\begin{theorem}
Let $G_1,...,G_m \in \mathcal{IPV}$ and let $G= G_1  \times ... \times G_m$. Assume that $H$ is an arbitrary group and let $\theta: \cL(G) \to \cL(H)$ be a $*$-isomorphism, then there exist $u \in \mathscr{U}(\cL(H))$ and $H_1,...,H_m \leqslant H$ such that $H = H_1  \times ... \times H_m$ and $t_1,...,t_m > 0$ such that $t_1 ... t_m = 1$ and $
u \theta(\cL(G_i))^{t_{i}} u^* = \cL(H_i)$, for any $1\leq i\leq m$.
\end{theorem}

\noindent {\it Proof.}
For the reader's convenience we recycle the notations used in \cite[Theorems 4.3 and 4.16]{CdSS15}. In fact we follow the proofs of these theorems only adding in the new aspects of the technique. Thus, we suggest the reader review these proofs beforehand as we only include a proof of Claim 4.8, this being the only piece needed.

{\bf Claim.}
$\Sigma \cap \Omega$ is finite.

Let $O_i'=O_i \cap \Sigma$ and notice that $\Sigma \cap \Omega = \bigcup_{i=1}^{\infty}{O_i'}$. For every $k$ consider $R_k = \langle O_i' , i \in \{1,..,k\} \rangle$ and notice that it forms an ascending sequence of normal subgroups of $\Sigma$ such that $\bigcup_{k} R_k= \Sigma \cap \Omega$. Moreover, $[\Sigma : \Sigma_k] < \infty$ where $\Sigma_k = C_{\Sigma}(R_k)$. Since, $R_k \cap \Sigma_k$ is abelian and $[\Sigma: \Sigma_k] < \infty$ it follows that $R_k$ is virtually abelian. In particular, $\Sigma \cap \Omega $ is an amenable group. In the first part of the proof of \cite[Theorem A]{CdSS15} we have obtained that $\mathscr Q \subset qL(\Sigma)q$ is a finite index inclusion of II$_1$ factors. Letting $z=z(q) \in \mathscr{Z}(\cL(\Sigma))$ be the central support of $q$, we have that for $s > 0$, $\mathscr Q^s \subseteq (q\cL(\Sigma)q)^s = \cL(\Sigma)z$ is a finite inclusion of II$_1$ factors. Then perform the basic construction for $\mathscr Q^s \subseteq \cL(\Sigma)z \subseteq \langle L(\Sigma)z,e_{\mathscr Q^{s}} \rangle = \mathscr Q^{\mu}$ where $\mu = s[qL(\Sigma)q:\mathscr Q]^2$. First, we argue that each $R_k$ is finite. Since $C_k = R_k \cap \Sigma_k \leqslant R_k$ has finite index, it suffices to show that $C_k$ is finite. From construction, we have that \[
\cL(C_k) \subseteq \mathscr Z(\cL(\Sigma_k)) \subseteq \cL(\Sigma_k) ' \cap \cL(\Sigma) \text{.}
\]
By passing to a finite index subgroup we can assume that $\Sigma_k \leqslant \Sigma$ is normal and $[\Sigma:\Sigma_k] = r < \infty$. Let $\gamma_1,\gamma_2,...,\gamma_r$ be a complete set of representatives for $\Sigma_k \leqslant \Sigma$. One can check that the map $E: \cL(\Sigma_k)' \cap \cL(\Sigma) \to \mathscr Z(\cL(\Sigma))$ given by $E_{\mathscr Z(\cL(\Sigma))}(x)=\frac{1}{r}\sum\limits_{i=1}^{r}{u_{\gamma_{i}}xu_{\gamma_{i}^{-1}}}$ is a conditional expectation satisfying:
\[
\|E_{\mathscr Z(\cL(\Sigma))}(x)\|_2^2 \geq r^{-1}\|x\|_2^2
\]
and hence $[\cL(\Sigma_k)' \cap \cL(\Sigma):\mathscr Z(\cL(\Sigma))] \leq r$ and thus $(\cL(\Sigma_k) \cap \cL(\Sigma))z$ is finite dimensional. Hence, there exists a $z_0 \in \mathscr{P}(\cL(\Sigma_k)' \cap \cL(\Sigma))$ such that $(\cL(\Sigma_k)' \cap \cL(\Sigma))z_0 = \mathbb{C} z_0$. This of course implies that $\cL(C_k)z_0 = \mathbb{C}z_0$. By \cite[Corollary 2.7]{CdSS15} one gets that $C_k$ is finite.

\noindent We now show that $\Sigma \cap \Omega$ is finite. Assume by contradiction that $\Sigma \cap \Omega$ is infinite. Now, $\cL(\Sigma \cap \Omega)z \subset \cL(\Sigma)z$ is a diffuse subalgebra where $z \in \mathscr{Z}(\cL(\Sigma)) \subseteq \cL(\Sigma \cap \Omega)$ and $\mathscr{Z}(\cL(\Sigma)) = \mathbb{C}z$. Let $\mathcal{H} = \{u_{\sigma}z : \sigma \in \Sigma\} \subset \mathscr{N}_{\cL(\Sigma)z}(\cL(\Sigma \cap \Omega)z)$ and notice that the action $\mathcal{H} \curvearrowright \cL(\Sigma \cap \Omega)z$ is weakly compact. Indeed, consider the self-adjoint element $\xi_k = |R_k|^{-\frac{1}{2}} \sum\limits_{a \in R_{k}}{u_{a}z \otimes \overline{u_{a}z}} \in \cL(R_k)z \bar\otimes \overline{\cL(R_{k})z}$. Note that $(u_{\gamma}z \otimes \overline{u_{\gamma}z})\xi_k = \xi_k(u_{\gamma}z \otimes \overline{u_{\gamma}z})$, $\forall \gamma \in \Sigma$. Now, for all $a \in R_k$, $\l \geq k$ we have $(u_a z \otimes \overline{u_a z})\xi_{l} = \xi_{l}$ and hence $\lim_{n}{\|(u_{a}z \otimes \overline{u_a z})\xi_n - \xi_n\|_{2,z} }= 0$, for all $ a \in \Sigma \cap \Omega$. Here, $\|\cdot \|_{2,z}$ is the $2$-norm induced by the trace $\tau_{z}(y)=\frac{\tau(yz)}{\tau(z)}$ on $L(\Sigma)z$ where $\tau$ is the canonical trace on $\cL(\Sigma)$. 

\noindent To this end notice we also have that 

\begin{equation}\begin{split}
&\langle (xz \otimes z)\xi_k, \xi_k \rangle = |R_k|^{-1} \sum\limits_{a,b \in R_{k}}{\langle xu_{a}z,u_b z \rangle \langle \overline{u_{a}z}, \overline{u_b z}} \rangle
\\
&
= |R_{k}|^{-1} \sum\limits_{a,b \in R_{k}}{\tau(xu_{a}zu_{b^{-1}}) \tau(zu_{a^{-1}b})} = |R_{k}|^{-1} \sum\limits_{a,s \in R_{k}}{\tau(xu_{a}zu_{s^{-1}a^{-1}})\tau(zu_{s})} 
\\
&= |R_k|^{-1} \sum\limits_{a \in R_{k}}{\tau(xu_a z(\sum\limits_{s \in R_k}{\tau(zu_s)u_{s^{-1}}})u_{a^{-1}})} \\
&= |R_{k}|^{-1} \sum\limits_{a \in R_{k}}{\tau(x u_{a}z E_{\cL(R_{k})}(z)u_{a^{-1}}}) = \tau(xzE_{\cL(R_{k})}(z) )\text{.} 
\end{split}
\end{equation}
Since $\bigcup_{k}{R_{k}} = \Sigma \cap \Omega$ we have that \[
\lim_{k}{\langle (xz \otimes z) \xi_n, \xi_n \rangle_{z}} = \lim_{k}{\tau_{z}(xzE_{\cL(R_k)}(z)) = \tau_z(x) \text{.}} 
\]
Similarly, $\lim_{k}{\langle (z \otimes \overline{xz}) \xi_k , \xi_k \rangle} = \tau_z(x)$, for all $ x \in \cL(\Sigma)$. 
Since $\cQ^{\mu} = \cL(G_I)^t$ where $t = \tau(p)\mu$ and $I=\hat 1$, from above we have that 
\[
\cL(\Omega \cap \Sigma)z \subset \cL(\Sigma)z \subset \cL(G_I)^t \text{.}
\]
Note that the last inclusion is an irreducible inclusion of finite index II$_1$ factors. Next, we show this leads to a  contradiction. When $|I|=1$ this already follows from \cite[Theorem 6.1]{CSU13}, so assume that $|I| \geq 2$. 
Write $\cL(G_I)^t = e (\cL(G_I) \otimes \mathcal M_n(\mathbb{C}))e$ for some projection $e \in \mathscr{P}(\cL(G_I) \otimes \mathcal M_n(\mathbb{C}))$. Fix, $i \in I$. First, we observe that $\cL(\Sigma)z$ is strongly non-amenable relative to $\cL(G_{I - \{i\}}) \otimes \mathcal M_n(\mathbb{C})$. Assume otherwise, since $\cL(\Sigma)z \subset \cL(G_{I})^t$ has finite index then $\cL(G_{I})^t$ is amenable relative to $\cL(G_{I-\{i\}}) \otimes M_n(\mathbb{C})$. Thus, by \cite[Proposition 2.4(3)]{OP07} we would have that $\cL(G_{I})^t$ is amenable relative to $\cL(G_{I-\{i\}}) \otimes \mathcal M_n({\mathbb{C}})$ and thus $G_i$ is amenable, a contradiction.

\noindent Then by Theorem \ref{convweakcomp} we get that $1 \otimes \alpha_t^{i} \to id$ uniformly on $(\cL(\Omega \cap \Sigma)z)_1$. Here, $1 \otimes \alpha_t^{i}$ is defined on $(\cL(G_{I-\{i\}} \otimes \mathcal M_n(\mathbb{C})) \otimes \cL(G_i)$ where $\alpha_t^i$ is the core length deformation on the wreath product algebra $\cL(G_i)=\cL(A^i \wr_{G_0/K_{0}} G_{0}^i)$. Thus, using \cite[Theorem 4.2]{IPV10} one of the following must hold:
\begin{enumerate}
    \item $\cL(\Omega \cap \Sigma) \preceq \cL(G_{I -\{i\}}) \otimes \mathcal M_n(\mathbb{C})$;
    \item $\cL(\Sigma)z \preceq \cL(G_{I-\{i\}}) \otimes \mathcal M_n(\mathbb{C}) \bar \otimes \cL(A^i \rtimes {\rm stab}_{G_{0}^i}(hK_{0}))$;
    \item there exists $v$ partial isometry such that $vv^*=z$ and $v\cL(\Sigma)zv^* \subset \cL(G_{i}) \otimes \mathcal M_{n}(\mathbb{C}) \bar \otimes \cL(G_{0}^i)$.
\end{enumerate}
Notice that $(\cL(\Sigma)z)' \cap \cL(G_{I})^t = \mathbb{C}z$. Then by \cite[Lemma 2.4]{DHI16} all intertwinings in 1) and 2) are strong. Next, we argue that 2) and 3) do not hold.

Assume 2) holds. Since $\cL(\Sigma)z \subseteq \cL(G_{I})^t$ is finite index we have by \cite[Lemma 2.9(2)]{Dr19b} that
\[
\cL(G_{I})^t \preceq \cL(G_{I-\{i\}}) \otimes \mathcal M_{n}(\mathbb{C}) \bar \otimes \cL(A_i \rtimes {\rm stab}_{G_{0}^i}(hK_{0}))
\]
but this implies that the inclusion $A_i \rtimes {\rm stab}_{G_{0}^i}(hK_{0}) \leqslant G_{i}$ has finite index which contradicts the fact that ${\rm stab}_{G_{0}^i}(hK_{0})=hK_0h^{-1}$ is amenable and $G_{i}$ is non-amenable.

\noindent Now, assume that 3) holds. Reasoning the same way, we have that 
\[
\cL(G_I)^t \preceq \cL(G_{I-\{i\}}) \otimes \mathcal M_n(\mathbb{C}) \bar \otimes \cL(G_0^i)
\]
which further implies that the inclusion $G_0^i \leqslant  G_i$ has finite index, a contradiction.

\noindent In conclusion, we have obtained that for all $i \in I$ we have  \[
\cL(\Omega \cap \Sigma) z \preceq^s \cL(G_{I-\{i\}}) \otimes \mathcal M_n(\mathbb{C}).
\]
Combining this with \cite[Lemma 2.8]{DHI16} inductively we get that \[
\cL(\Omega \cap \Sigma)z \preceq^{s} \bigcap_{i \in I}{\cL(G_{I-\{i\}}) \otimes \mathcal M_n(\mathbb{C})} = 1 \otimes \mathcal M_n(\mathbb{C}) . 
\]

\noindent This implies that a corner of $\cL(\Omega \cap \Sigma)z$ is atomic and hence, there exists a non-zero projection $z_0 \in \mathscr Z(\cL(\Omega \cap \Sigma)z)$ such that $\cL(\Omega \cap \Sigma)z_0 = \mathbb{C}z_0$. Thus, by applying \cite[Corollary 2.7]{CdSS15} we get that $\Omega \cap \Sigma$ is finite, contradiction. \hfill$\blacksquare$

\end{document}